\documentclass[reqno,11pt]{amsart}
\usepackage{amsfonts}
\usepackage{graphicx}
\usepackage{amsfonts,amsmath, amssymb}
\usepackage{bbm,mathrsfs,mathscinet}
\usepackage{hyperref}
\usepackage{marginnote}
\usepackage{color}

\numberwithin{equation}{section}

\newcommand{\fN}{\mathfrak{N}}
\newcommand{\fu}{\mathfrak{u}}
\newcommand{\fb}{\mathfrak{b}}
\newcommand{\fa}{\mathfrak{a}}
\newcommand{\fw}{\mathfrak{w}}
\newcommand{\R}{\mathbb{R}}
\newcommand{\MH}{\mathcal{H}}
\newcommand{\bp}{\partial_{t,\sp}}
\newcommand{\tw}{\tilde{w}}

\newtheorem{theorem}{Theorem}[section]

\newtheorem{lemma}[theorem]{Lemma}
\newtheorem{proposition}[theorem]{Proposition}
\newtheorem{remark}[theorem]{Remark}

\def\v{\varepsilon}

\def\x{\xi}
\def\t{\theta}
\def\k{\kappa}

\def\a{\alpha}
\def\b{\beta}
\def\g{\gamma}
\def\d{\delta}

\def\r{\rho}
\def\s{\sigma}

\def\o{\omega}

\def\i{\infty}
\def\f{\frac}
\def\pa{\partial}
\def\fep{f^{\varepsilon}}
\def\hep{h^{\varepsilon}}

\def\FRep{F^{\varepsilon}_{R}}
\def\pt{\partial_{t}}
\def\nax{\nabla_{x}}
\def\dis{\displaystyle}
\def\la{\langle}
\def\ra{\rangle}
\def\intr{\int_{\mathbb{R}^{3}}}
\def\ints{\int_{\mathbb{S}^{2}}}
\def\smum{\sqrt{\mu_{M}}}
\def\sp{\shortparallel}

\begin{document}
	
	\title[Hilbert expansion of Boltzmann equation]{Hilbert expansion of the Boltzmann equation with specular boundary condition in  half-space}

	\author[Y. Guo]{Yan Guo}
	\address[Y. Guo]{The Division of Applied Mathematics, Brown University, Providence,
		USA.}
	\email{yan\_guo@brown.edu}
	
	\author[F.M. Huang]{Feimin Huang}
    \address[F.M. Huang]{Academy of Mathematics and Systems Science, Chinese Academy of Sciences, Beijing 100190, China; School of Mathematical Sciences, University of Chinese Academy of Sciences, Beijing 100049, China}
    \email{fhuang@amt.ac.cn}
	
	\author[Y. Wang]{Yong Wang}
	\address[Y. Wang]{Academy of Mathematics and Systems Science, Chinese Academy of Sciences, Beijing 100190, China; School of Mathematical Sciences, University of Chinese Academy of Sciences, Beijing 100049, China}
	\email{yongwang@amss.ac.cn}

\begin{abstract}
	
Boundary effects play an important role in the study of hydrodynamic limits in the Boltzmann theory. Based on a systematic derivation and study of the viscous layer equations and the $L^2$ to $L^\infty$ framework, we establish the validity of the Hilbert expansion for the Boltzmann equation with specular reflection boundary conditions, which leads to derivations of compressible Euler equations and acoustic equations.

\end{abstract}

%\subjclass[2000]{35Q20, 35B20, 35B35, 35B45}
\date{\today}
%\thanks{}
\maketitle
	
\setcounter{tocdepth}{2}
\tableofcontents
	
\thispagestyle{empty}

	%%%%%%%%%%%%%%%%%%%%%%%%%%%%%%%%%%%%%%%%%%%%%%%%%%%%%%%%%%%%%%%%%%
	
	\section{Introduction and Main Results}
	
\subsection{Introduction}

In  the founding work of Maxwell \cite{Maxwell} and Boltzmann \cite{Boltzmann},  it was shown that the Boltzmann equation is closely related to the fluid dynamical systems for both compressible and incompressible flows. A great effort has been devoted to the study of the hydrodynamic limit from the Boltzmann equation to the fluid systems. In 1912,   Hilbert proposed a systematic formal asymptotic expansion for Boltzmann equation  with respect to Knudsen number $\mathscr{K}_n \ll1$. A bit later  Enskog and Chapman independently proposed a somewhat different  formal expansion in 1916 and 1917, respectively. Either Hilbert or the Chapman-Enskog expansions yield the compressible and incompressible fluid equations, for instance: the compressible Euler and Navier-Stokes systems, the incompressible Euler and Navier-Stokes (Fourier) systems,  and the acoustic system, {\it et.al.} It is a challenging problem to rigorously justify these formal approximation. In fact, the purpose of Hilbert's sixth problem \cite{Hilbert} is to establish  the laws of motion of continua from more microscopic physical models, such as Boltzmann theory, from a rigorous mathematical standpoint.

\smallskip

Based on the truncated Hilbert expansion, Caflisch \cite{Caflish} rigorously established  the hydrodynamic limit from the Boltzmann equation to the  compressible Euler equations when solution is smooth, see also   \cite{Grad,Lachowicz,Nishida,Ukai-Asano} and \cite{Guo Jang Jiang-1,Guo Jang Jiang} via a recent $L^2$-$L^\infty$ framework.  As is well known, solutions of the compressible Euler equations in general develop singularities, such as shock waves. Generally, there are three basic wave patterns for compressible Euler equations, that is, shock wave, rarefaction wave, and contact discontinuity.
The hydrodynamic limit of  Boltzmann to such wave patterns have been proved   \cite{Huang-Jiang-Wang,Huang-Wang-Yang,Huang-Wang-Yang-1,Huang-Wang-Wang-Yang,Xin-Zeng,Yu} in one  dimensional case. For multi-dimensional case, the only result is \cite{WWZ} for planar rarefaction wave.

\smallskip

The acoustic equations are the linearization of the compressible Euler equations about a spatially homogeneous fluid state. And it is somehow the simplest system of fluid dynamics, being essentially the wave equation. Bardos, Golse and Levermore \cite{Bardos-3} established the convergence from the DiPerna-Lions  \cite{Diperna-Lions} solutions of Boltzmann equation to the solution of acoustic system over a periodic spatial domain with  restriction on the size of fluctuation. Recently, the restriction has been relaxed in \cite{JLM,Jang-Jiang}, and finally  removed in \cite{Guo Jang Jiang} via a $L^2$-$L^\infty$ framework.

\smallskip

There have been extensive research efforts and literature to derive the incompressible Navier-Stokes system, see  \cite{Bardos,Bardos-2,Bardos-Ukai,CEMP,Diperna-Lions,ELM-1994,E-Guo-M,E-Guo-K-M,Golse-Saint-Raymond,Guo2006,Guo-Liu,Jang-Kim,Jiang-Masmoudi,Masmoudi-Raymond,Saint-Raymond2009,Wu} and the references cited therein.
%For the case of  incompressible flows, the program was initiated by Bardos, Golse, and Levermore \cite{Bardos,Bardos-2} to justify the global weak solution of incompressible flows in the frame work of global renormalized solution of DiPerna-Lions \cite{Diperna-Lions}. In particular, Golse and Saint-Raymond \cite{Golse-Saint-Raymond} proved that the limits of the DiPerna-Lions renormalized solutions of the Boltzmann equation are the Leray solutions to the incompressible Navier-Stokes equations. There are also many  progresses on this topic such as \cite{Bardos-Ukai,Guo2006,Guo-Liu} without physical boundary and \cite{E-Guo-M,E-Guo-K-M,Jiang-Masmoudi,Masmoudi-Raymond,Jang-Kim,Wu} with physical boundary, see also the lecture note \cite{Saint-Raymond2009} and the references therein. We will not go into details about the incompressible limits since we will concentrate on the compressible Euler limit in this paper.

\smallskip

All of the above-mentioned works on the compressible Euler limit and the acoustic limit  were carried out in either the periodic spatial domain or the whole space. However,  in many important physical applications, boundaries occur naturally, and boundary effects are crucial in the hydrodynamic limit of dilute gases governed by the Boltzmann equation. Hence it is  important  to study  the hydrodynamic limit from the  Boltzmann equation to the compressible Euler equations in the appearance of  physical boundary.
The purpose of this paper is to justify the compressible Euler limit and the acoustic limit of Boltzmann equation with specular reflection boundary conditions by the Hilbert expansion method. The main difficulty, due to the appearance of physical boundary, is the possible appearance of  both viscous and Knudsen boundary layers,  and the interplay between these layers is  very complicate.

\smallskip

More precisely, we consider  the scaled Boltzmann equation
	\begin{equation}\label{1.0}
		F_t+v\cdot\nabla_x F=\frac1{\mathscr{K}_n}Q(F, F),
	\end{equation}
where $F(t,x,v)\geq 0$ is the density distribution function for the gas particles with position $x\in\mathbb{R}^3_+=\{x\in\mathbb{R}^3: x_3>0\}$  and velocity $v\in\mathbb{R}^3$ at  time $t>0$, and $\mathscr{K}_n>0$ is  Knudsen number which is proportional to the mean free path. The Boltzmann collision term $Q(F_1,F_2)$ on the right is defined in terms of the following bilinear form
	\begin{align}\label{1.2}
		Q(F_1,F_2)&\equiv\int_{\mathbb{R}^3}\int_{\mathbb{S}^2} B(v-u,\t)F_1(u')F_2(v')\,{d\omega du}
		-\int_{\mathbb{R}^3}\int_{\mathbb{S}^2} B(v-u,\t)F_1(u)F_2(v)\,{d\omega du}\nonumber\\
		&:=Q_+(F_1,F_2)-Q_-(F_1,F_2),
	\end{align}
where the relationship between the post-collision velocity $(v',u')$ of two particles with the pre-collision velocity $(v,u)$ is given by
	\begin{equation*}%\label{1.3}
		u'=u+[(v-u)\cdot\omega]\omega,\quad v'=v-[(v-u)\cdot\omega]\omega,
	\end{equation*}
	for $\omega\in \mathbb{S}^2$, which can be determined by conservation laws of momentum and energy
	\begin{equation*}%\label{1.3-1}
		u'+v'=u+v,\quad |u'|^2+|v'|^2=|u|^2+|v|^2.
	\end{equation*}
The Boltzmann collision kernel $B=B(v-u,\theta)$ in \eqref{1.2} depends only on $|v-u|$ and $\theta$ with $\cos\theta=(v-u)\cdot \omega/|v-u|$.   Throughout this paper,  we  consider the hard sphere model, i.e.,
\begin{equation*}%\label{1.4}
B(v-u,\t)=|(v-u)\cdot \omega|
\end{equation*}

We denote $\vec{n}=(0,0,1)$ to be  the outward normal of  $\mathbb{R}^3_+$. We denote the phase boundary in the space $\mathbb{R}^3_+\times\mathbb{R}^3$ as $\gamma:=\partial\mathbb{R}^3_+\times\mathbb{R}^3$ and split it into outgoing boundary $\gamma_+$, incoming boundary $\gamma_-$, and grazing boundary $\gamma_0$:
\begin{align}
\begin{split}\nonumber
\gamma_+=\{(x,v) : x\in \partial\mathbb{R}^3_+, v\cdot \vec{n}=-v_3>0\},\\
\gamma_-=\{(x,v) : x\in \partial\mathbb{R}^3_+, v\cdot \vec{n}=-v_3<0\},\\
\gamma_0=\{(x,v) : x\in \partial\mathbb{R}^3_+, v\cdot \vec{n}=-v_3=0\}.
\end{split}
\end{align}
In the present paper we consider the Boltzmann equation with specular reflection boundary conditions, i.e.,
\begin{equation}\label{1.5}
F(t,x,v)|_{\gamma_-}=F(t,x,R_xv),
\end{equation}
where
\begin{equation}\label{1.6-1}
R_xv=v-2\{v\cdot \vec{n}\} \vec{n}=(v_1,v_2, -v_3)^{\top}.
\end{equation}

%%%%%%%%%%%%%%%%%%%%%%%%%%%%%%%%%%%%%%%%%%%%%%%%%%%%%%%%%%%%%%%%%%%%%%%%%%%%%%%%%%%%%%

\subsection{Asymptotic Expansion} Since the thickness of viscous boundary layer is $\sqrt{\mathscr{K}_n}$, for simplicity, we  use the new parameter $\v=\sqrt{\mathscr{K}_n}$, and denote the Boltzmann solution to be $F^\v$, then the Boltzmann equation \eqref{1.0} is rewritten as
\begin{equation}\label{1.1}
\partial_t F^\v+v\cdot \nabla_x F^\v =\frac{1}{\v^2} Q(F^\v, F^\v).
\end{equation}

\subsubsection{Interior expansion:} We define the interior  expansion
	\begin{equation}\label{1.8}
		F^\v(t,x,v)\sim \sum_{k=0}^{\infty}\v^{k}F_{k}(t,x,v),
	\end{equation}
Substituting \eqref{1.8} into  \eqref{1.1} and comparing the order of $\v$, one obtains
\begin{align}\label{1.7-1}
\begin{split}
\dis \v^{-2}: \quad\quad\quad\quad\quad\qquad 0&=Q(F_{0},F_{0}),\\[2mm]
\dis \v^{-1}:  \quad\quad\quad\quad\quad\qquad 0&=Q(F_{0},F_{1})+Q(F_1,F_0),\\[2mm]
\dis \v^0:\quad \{\pt+v\cdot\nax\}F_{0}&=Q(F_{0},F_{2})+Q(F_{2},F_{0})+Q(F_1,F_1),\\[2mm]
\dis \v: \quad\{\pt+v\cdot\nax\}F_{1}&=Q(F_{0},F_{3})+Q(F_{3},F_{0})+Q(F_{1},F_{2})+Q(F_2,F_1),\\[2mm]
\dis \quad\quad\quad\quad\quad\quad\quad\ \vdots\\
\dis \v^{k}:\quad  \{\pt+v\cdot\nax\}F_{k}&=Q(F_{0},F_{k+2})+Q(F_{k+2},F_{0})+\sum_{\substack{i+j=k+2\\  i,j\geq1}}Q(F_{i},F_{j}).
%\dis \quad\quad\quad\quad\quad\quad\quad\ \vdots\\
\end{split}
\end{align}
%	We can construct smooth function $F_{1}(t,x,v),...,F_{5}(t,x,v)$ for $ 0\leq t\leq \tau$. For more detailed discussion, see \cite{C,R,E}.

It follows from  $\eqref{1.7-1}_1$ and the celebrated H-theorem that $F_0$ should be a local Maxwellian
\begin{equation*}%\label{1.6}
\mu(t,x,v):=F_0(t,x,v)\equiv \f{\rho(t,x)}{[2\pi T(t,x)]^{3/2}}\exp{\left\{-\f{|v-\fu(t,x)|^2}{2T(t,x)}\right\}},
\end{equation*}
where $\rho(t,x)$, $\mathfrak{u}(t,x)=(\fu_1,\fu_2,\fu_3)(t,x),$ and $T(t,x)$ are defined by
\begin{align}\nonumber
\intr F_0 dv =\rho,\quad \intr vF_0dv =\rho \fu, \quad \intr|v|^2F_0 dv =\rho|\fu|^2+3\rho T,
\end{align}
which represent the macroscopic density, velocity and temperature, respectively.  Projecting the equation  $\eqref{1.7-1}_3$ onto $1$, $v$, $\f{|v|^2}{2}$, which are five collision invariants for the Boltzmann collision operator $Q(\cdot,\cdot)$, one obtains that $(\rho, \fu, T)$ satisfies the compressible Euler system
	\begin{equation}\label{1.7}
\begin{cases}
\dis \pt\r+\mbox{\rm div} (\r \fu)=0,\\[2mm]
\dis \pt(\r\fu)+\mbox{\rm div}(\r \fu\otimes \fu)+\nabla p=0,\\[2mm]
\dis \pt[\r(\f{3T}{2}+\f{|\fu|^2}{2})]+\mbox{\rm div} [\r \fu(\f{3T}{2}+\f{|\fu|^2}{2})]+\mbox{\rm div}(p\fu)=0,
\end{cases}
\end{equation}
where $x\in\mathbb{R}^3_+, t>0$ and  $p=\rho T$ is the pressure function. For the compressible Euler equations \eqref{1.7}, we impose the slip boundary condition
\begin{equation}\label{1.12}
\fu\cdot \vec{n}|_{x_3=0}=\fu_{3}|_{x_3=0}=0.
\end{equation}
and the initial data
\begin{equation}\label{1.12-2}
(\rho,\fu,T)(0,x)=(1+\delta\varphi_0,\delta \Phi_0, 1+\delta \vartheta_0)(x),
\end{equation}
with $\|(\varphi_0, \Phi_0,\vartheta_0)\|_{H^{s_0}}\leq 1$ where $\delta>0$ is a  parameter and $s_0\geq 3$ is some given positive number.
Choose $\delta_1>0$ so that for any $\delta\in(0,\delta_1]$, the positivity of $1+\delta\varphi_0$ and $1+\delta \vartheta_0$ is guaranteed. Then for each $\delta\in(0,\delta_1]$, there is a family of classical solutions $(\rho^{\delta},\fu^{\delta}, T^{\delta})\in C([0,\tau^{\delta}]; H^{s_0}(\mathbb{R}^3_+)) \cap C^1([0,\tau^{\delta}]; H^{s_0-1}(\mathbb{R}^3_+))$  of the compressible Euler equations \eqref{1.7}-\eqref{1.12-2} such that $\rho^{\delta}>0$ and $T^{\delta}>0$, see Lemma \ref{lem2.1-1} for  details.

Generally, the solution of  interior expansion $F_i , i=1,2,\cdots$ does not  satisfy the specular reflection boundary conditions. %We will give an explicit example to explain it.
Then it fails to prove the summation $ \sum_{k=0}^{\infty}\v^{k}F_{k}(t,x,v)$ satisfies the specular boundary condition. To overcome the difficulty coming from physical boundary,  the boundary layer expansions is needed, see \cite{Sone-2002,Sone-2007}.

\smallskip

For later use, we define the linearized collision operator $\mathbf{L}$ by
\begin{equation}%\label{2.1}
\mathbf{L}g=-\frac{1}{\sqrt{\mu}}\Big\{Q(\mu,\sqrt{\mu} g)+Q(\sqrt{\mu} g,\mu)\Big\}.\nonumber
\end{equation}
%and the nonlinear operator
%\begin{equation*}
%\Gamma(g_1,g_2)=\frac{1}{\sqrt{\mu}} Q(\sqrt{\mu}g_1,\sqrt{\mu}g_2).
%\end{equation*}
The null space $\mathcal{N}$ of $\mathbf{L}$ is generated by
\begin{align}\nonumber%\label{2.2}
\begin{split}
\chi_0(v)&\equiv\f{1}{\sqrt{\rho}}\sqrt{\mu},\\
\chi_i(v)&\equiv\f{v_i-\fu_{i}}{\sqrt{\rho T}}\sqrt{\mu},\quad i=1,2,3,\\
\chi_4(v)&\equiv\f{1}{\sqrt{6\rho}}\left\{\f{|v-\fu|^2}{T}-3\right\}\sqrt{\mu}.
\end{split}
\end{align}
It is easy to check that $\displaystyle\int_{\mathbb{R}^3}\chi_i\cdot \chi_j dv=\delta_{ij}$ for $0\leq i,j\leq 4$.
We also define the collision frequency $\nu$:
\begin{equation}\label{1.19}
\nu(t,x,v)\equiv\nu(\mu):=\intr\ints B(v-u,\t)\mu(u)d\o du.
\end{equation}
It is direct to know that
\begin{equation}\nonumber
\frac1{C} (1+|v|)^{\gamma}\leq \nu(t,x,v) \leq  C(1+|v|)^{\gamma},
\end{equation}
where  $C>0$ is some given positive constant.
%Throughout this paper, we shall  $\|\cdot\|_\nu$ to donate the space-velocity inner product $\la\cdot,\cdot\ra_\nu$.
Let $\mathbf{P}g$ be the $L_{v}^{2}$ projection with respect to $[\chi_0,...,\chi_4]$. It is well-known that there exists a positive number $c_{0}>0$ such that for any function $g$%\in\mathcal{N}^{\perp}
\begin{equation}\nonumber%\label{coercivity}
\la \mathbf{L}g,g\ra\geq c_{0}\|\{\mathbf{I-P}\}g\|_{\nu}^{2},
\end{equation}
where the weighted $L^2$-norm $\|\cdot\|_{\nu}$ is defined as
\begin{equation}\nonumber
\|g\|_{\nu}^2:=\int_{\mathbb{R}^3_x\times\mathbb{R}^3_v}  g^2(x,v) \nu(v)dxdv.
\end{equation}

For each $k\geq 1$, we define the macroscopic and microscopic part of  $\displaystyle\frac{F_k}{\sqrt{\mu}}$ as
\begin{align}\label{1.25}
\frac{F_k}{\sqrt{\mu}}
&=\mathbf{P}\left(\frac{F_k}{\sqrt{\mu}}\right)+\{\mathbf{I-P}\}\left(\frac{F_k}{\sqrt{\mu}}\right)\nonumber\\
&\equiv \left\{\frac{\rho_k}{\sqrt{\rho}} \chi_0+\sum_{j=1}^3\sqrt{\frac{\rho}{T}} u_{k,j}\cdot \chi_j+\sqrt{\frac{\rho}{6}}\frac{\theta_k}{T} \chi_4 \right\}+\{\mathbf{I-P}\}\left(\frac{F_k}{\sqrt{\mu}}\right)\nonumber\\
&\equiv \left\{ \frac{\rho_k}{\rho}+u_{k}\cdot \frac{v-\fu}{T}+\frac{\theta_k}{6T}(\frac{|v-\fu|^2}{T}-3)\right\}\sqrt{\mu}+\{\mathbf{I-P}\}\left(\frac{F_k}{\sqrt{\mu}}\right).
\end{align}

\subsubsection{ Viscous boundary layer expansion:}
We  define the scaled normal coordinate:
\begin{equation}\label{1.12-1}
y:=\frac{x_3}{\v}.
\end{equation}
For simplicity of presentation, we denote
\begin{equation}\label{1.12-3}
x_{\shortparallel}=(x_1,x_2),\quad  \nabla_{\sp}=(\partial_{x_1},\partial_{x_2})\quad\mbox{and}\quad  v_{\sp}=(v_1,v_2).
\end{equation}

Noting \eqref{1.12}, it is direct to know that  the local Maxwellian $\mu$ satisfies the specular reflection boundary conditions. However, in general, $F_1$ may not satisfy the specular reflection boundary conditions, therefore we need to construct viscous boundary layer to compensate the boundary condition starting from the first order of $\v$.

\smallskip

Motivated by \cite{Sone-2007}, we define the viscous boundary layer expansion as
\begin{align*}%\label{1.13}
\bar{F}^\v(t,x_\sp, y)\sim \sum_{k=1}^\infty \v^k \bar{F}_k(t,x_\sp, y,v).
\end{align*}
Plugging $F^\v+\bar{F}^\v$  into the Boltzmann equation \eqref{1.1} and comparing the order of $\v$, then using \eqref{1.7-1},  in the neighborhood of  physical boundary, we have
\begin{equation}\label{1.14}
\begin{split}
\dis& \v^{-1}:  \qquad\quad 0=Q(\mu_0,\bar{F}_{1})+Q(\bar{F}_1,\mu_0),\\[2mm]
\dis & \v^0:\quad v_3 \frac{\partial\bar{F}_1}{\partial y }=[Q(\mu_0,\bar{F}_{2})+Q(\bar{F}_{2},\mu_0)]+y [Q(\partial_3\mu_0, \bar{F}_1)+Q( \bar{F}_1,\partial_3\mu_0)]\\
&\qquad\qquad\qquad\qquad+Q(F_1^0,\bar{F}_1)+Q(\bar{F}_1,F_1^0)+Q(\bar{F}_1,\bar{F}_1),\\
&\dis \quad\quad\quad\quad\quad\quad\quad\ \vdots\\
\dis & \v^{k}: \quad  \{\partial_t+v_\sp\cdot\nabla_\sp\}\bar{F}_{k}+v_3\frac{\partial\bar{F}_{k+1}}{\partial y }\\
&\qquad=Q(\mu_0,\bar{F}_{k+2})+Q(\bar{F}_{k+2},\mu_0) +\sum_{\substack{l+j=k+2\\    1\leq l\leq \fb,\, j\geq1}}\frac{y^l}{l!} \big[Q(\partial_3^l\mu_0, \bar{F}_{j})+Q( \bar{F}_{j},\partial_3^l\mu_0) \big] \\
&\qquad\quad+\sum_{\substack{i+j=k+2\\   i, j\geq1}}\big[Q(F_i^0,\bar{F}_j)+Q(\bar{F}_j,F_i^0)+Q(\bar{F}_i,\bar{F}_j)\big] \\
&\qquad\quad+\sum_{\substack{i+j+l=k+2\\  1\leq l\leq \fb,\, i, j\geq1}} \frac{y^l}{l!} \big[Q(\partial_3^lF_i^0, \bar{F}_{j})+Q( \bar{F}_{j},\partial_3^lF_i^0)\big], \quad \mbox{for}\  k\geq 1,
\end{split}
\end{equation}
where we have used  the Taylor expansions of $\mu$ and $F_i$ at $x_3=0$, i.e.,
\begin{align}\label{1.14-1}
\mu(t,x_1,x_2,x_3,v)
&=\mu_0+\sum_{l=1}^{\fb} \frac{1}{l!} \partial_3^l\mu_0\cdot x_3^l + \frac{x_3^{\fb+1}}{(\fb+1)!} \partial_3^{\fb+1}\tilde{\mu},
\end{align}
and for $i\geq 1$
\begin{align}\label{1.14-2}
F_i(t,x_1, x_2, x_3,v)
&=F_i^0+\sum_{l=1}^{\fb} \frac{1}{l!} \partial_3^l F_i^0\cdot x_3^l + \frac{x_3^{\fb+1}}{(\fb+1)!} \partial_3^{\fb+1}\mathfrak{F}_i.
\end{align}
Here we have used the simplified notations
\begin{align}\label{1.14-3}
\begin{split}
\partial_3^l\mu_0:&=(\partial_3^l\mu)(t,x_1, x_2,0,v),
\quad \partial_3^{\fb+1}\tilde{\mu}:=(\partial_3^{\fb+1}\mu)(t,x_1, x_2, \xi_0,v),\\
 \partial_3^lF_i^0:&=(\partial_3^lF_i)(t,x_1, x_2, 0, v),
\quad \partial_3^{\fb+1}\mathfrak{F}_i:=(\partial_3^{\fb+1}F_i)(t,x_1, x_2, \xi_i,v),
\end{split}
\end{align}
for some $\xi_i \in (0,x_3)$ with $i\geq0$.  The number $\fb\in \mathbb{N}_+$ will be chosen later.

\smallskip

The main reason to use \eqref{1.14-1}  is to make the linearized operator of \eqref{1.14} be independent of $y$. And the reason for \eqref{1.14-2} is to make the coefficients of viscous boundary layer system be independent of $\v$, see \eqref{2.19-2}-\eqref{2.20} for details. Noting the polynomial growth coefficients $y^l$ in  $\eqref{1.14}$, it is imperative  to  prove that they decay with enough polynomial  rate as $y\rightarrow+\infty$ provided the initial data decay sufficiently fast. For later use, we define $\bar{F}_0=0$.

\smallskip

For the macro-micro decomposition of viscous and Knudsen boundary layers, we define the corresponding linearized collision operator,  macroscopic projection,  and null space as
\begin{align*}%\label{2.8}
\mathbf{L}_0=\mathbf{L}(t,x_\sp,0,v),\qquad \mathbf{P}_0=\mathbf{P}(t,x_\sp,0,v), \qquad \mathcal{N}_0=\mathcal{N}(t,x_\sp,0,v).
\end{align*}
It is noted that $\mathbf{L}_0, \mathbf{P}_0$ and $\mathcal{N}_0$ are independent of normal variable. We define
\begin{equation}\label{2.8-1}
\bar{f}_k:=\frac{\bar{F}_k}{\sqrt{\mu_0}},
\end{equation}
 then it holds that
\begin{align*}%\label{2.10}
\bar{f}_k&=\mathbf{P}_0\bar{f}_k+\{\mathbf{I-P_0}\}\bar{f}_k\nonumber\\
&= \left\{ \frac{\bar\rho_k}{\rho^0}+\bar{u}_{k}\cdot \frac{v-\fu^0}{T^0}+\frac{\bar{\theta}_k}{6T^0}(\frac{|v-\fu^0|^2}{T^0}-3)\right\}\sqrt{\mu_0}+\{\mathbf{I-P_0}\}\bar{f}_k.
\end{align*}
where and whereafter we always use the notation $(\rho^0,\fu^0,T^0):=(\rho,\fu,T)(t,x_\sp,0)$.

\smallskip

Throughout the present paper, we always assume the far field condition
\begin{equation}\label{2.11}
\bar{f}_k(t,x_\sp,y,v)\rightarrow 0, \quad \mbox{as}\  y\rightarrow+\infty.
\end{equation}
In fact, it follows from $\eqref{1.14}_1$ that $\bar{f}_1\in \mathcal{N}_0$, i.e.,
\begin{equation}\label{2.12}
\bar{f}_1\equiv\mathbf{P}_0\bar{f}_1= \left\{ \frac{\bar{\rho}_1}{\rho^0}+\bar{u}_{1}\cdot \frac{v-\fu^0}{T^0}+\frac{\bar{\theta}_1}{6T^0}(\frac{|v-\fu^0|^2}{T^0}-3)\right\}\sqrt{\mu_0}.
\end{equation}
We denote
\begin{equation}\label{Boussinesq}
\bar{p}_k=\frac{\rho^0 \bar{\theta}_k+3T^0 \bar{\rho}_k}{3}.
\end{equation}
Multiplying $\eqref{1.14}_2$ by $\sqrt{\mu_0}, v_3\sqrt{\mu_0}$ and integrating over $\R^3$ with respect to $v$, one obtains
\begin{equation}\label{2.13}
\partial_y\bar{u}_{1,3}=0,\quad
\partial_y\bar{p}_1=0,
\end{equation}
where $\bar{p}_1$ is the one defined in \eqref{Boussinesq} with $k=1$. Noting from \eqref{2.11} and \eqref{2.13}, we have
\begin{equation}\label{2.14}
\bar{u}_{1,3}(t,x_\sp,y)\equiv0\quad \mbox{and}\quad
\bar{p}_1(t,x_\sp,y)\equiv0,\quad \forall\,  (t,x_\sp,y)\in [0,\tau]\times\R^2\times\R_+.
\end{equation}
It is noted that $\eqref{2.14}_2$ is similar to the Boussinesq relation in the diffusive limit of Boltzmann equation.

\smallskip

For later use, we define the Burnett functions
$\mathcal{A}_{ij}$ and $\mathcal{B}_{i}$
\begin{align}\label{2.4}
\begin{split}
\mathcal{A}_{ij}&:=\left\{\f{(v_i-\fu_i)(v_j-\fu_j)}{T}-\d_{ij}\f{|v-\fu|^2}{3T}\right\}\sqrt{\mu},\\
\mathcal{B}_{i}&:=\f{v_i-\fu_i}{2\sqrt{T}}\left(\f{|v-\fu|^2}{T}-5\right)\sqrt{\mu}.
\end{split}
\end{align}
and denote
\begin{equation*}
\mathcal{A}_{i,j}^0:=\mathcal{A}_{i,j}(t,x_\sp,0,v) \quad \mbox{and}\quad \mathcal{B}_{i}^0:=\mathcal{B}_{i}(t,x_\sp,0,v).
\end{equation*}
We define the viscosity and thermal conductivity coefficients  $\mu(T^0),\, \kappa(T^0)$ of viscous boundary layer
\begin{align}\label{2.15-2}
\begin{split}
\mu(T^0)&:=T^0\langle \mathcal{A}_{31}^0,\  \mathbf{L}_0^{-1} \mathcal{A}_{31}^0 \rangle\equiv T^0\langle\mathcal{A}_{ij}^0,\  \mathbf{L}_0^{-1} \mathcal{A}_{ij}^0 \rangle,\quad \forall i\neq j,\\[1.5mm]
\kappa(T^0)&:=\frac23T^0\langle \mathcal{B}_{3}^0,\  \mathbf{L}_0^{-1} \mathcal{B}_{3}^0\rangle \equiv \frac23T^0 \langle\mathcal{B}_{i}^0,\  \mathbf{L}_0^{-1} \mathcal{B}_{i}^0 \rangle,
\end{split}
\qquad i,j=1,2,3.
\end{align}
By Lemma 4.4 in \cite{Bardos-2}, it holds that $\langle T^0\mathcal{A}_{33}^0, \mathbf{L}_0^{-1} \mathcal{A}_{33}^0\rangle=\frac43\mu(T^0) $.
We denote
\begin{align*}%\label{2.15-3}
f_1^0:&=f_1(t,x_\sp,0,v)\equiv \left\{ \frac{\rho_1^0}{\rho^0}+u_{1}^0\cdot \frac{v-\fu^0}{T^0}+\frac{\theta_1^0}{6T^0}(\frac{|v-\fu^0|^2}{T^0}-3)\right\}\sqrt{\mu_0},
\end{align*}
where $(\rho_1^0,u_1^0,\theta_1^0):=(\rho_1,u_1,\theta_1)(t,x_\sp,0)$.

\smallskip

$\bar{f}_k$ will be constructed inductively as follows:
\begin{theorem}\label{thm2.5}
	Let $\bar{f}_k\  (k\geq1)$ be the solution of \eqref{1.14}. For  $k\geq1$, $(\bar{u}_{k,\sp},\bar{\theta}_k)$ satisfies
	\begin{align}\label{2.19-2}
	&\rho^0 \partial_t \bar{u}_{k,i}+\rho^0 (\fu^0_{\sp} \cdot\nabla_\sp) \bar{u}_{k,i}+\rho^0 (\partial_3\fu_3^0\cdot y+u^0_{1,3}) \partial_y \bar{u}_{k,i}\nonumber\\
	&\qquad\qquad\qquad\qquad\qquad+\rho^0 \bar{u}_{k,\sp}\cdot \nabla_\sp \fu^0_i-\frac{\partial_i p^0}{3T^0} \bar{\theta}_k-\mu(T^0) \partial_{yy} \bar{u}_{k,i}\nonumber\\
	&=\bar{\mathfrak{f}}_{k-1}=:-\rho^0 \partial_y[(\partial_3\fu_i^0\cdot y+u^0_{1,i}+\bar{u}_{1,i})   \bar{u}_{k,3}]-[\partial_i-\frac{\partial_ip^0}{p^0}]\bar{p}_k\\
	&\qquad\qquad\qquad\qquad\qquad\qquad\qquad+\bar{W}_{k-1,i}-T^0\partial_y\langle\bar{J}_{k-1},\  \mathcal{A}_{3i}^0\rangle,\quad i=1,2,\nonumber\\[2mm]
	&\rho^0 \partial_t \bar{\theta}_{k}+\rho^0 (\fu^0_{\sp} \cdot\nabla_\sp) \bar{\theta}_{k}+\rho^0 (\partial_3\fu_3^0\cdot y+u^0_{1,3}) \partial_y\bar{\theta}_k+\frac23\rho^0 \mbox{\rm div}\fu^0 \bar{\theta}_k-\frac35\kappa(T^0) \partial_{yy}\bar{\theta}_k\nonumber\\
	&=\bar{\mathfrak{g}}_{k-1}=:-\rho^0 \partial_y[(3\partial_3T^0\cdot y+\theta^0_{1}+\bar{\theta}_{1})   \bar{u}_{k,3}]+\frac35\bar{H}_{k-1}-\frac65(T^0)^{\frac32}\partial_y\langle \bar{J}_{k-1},\  \mathcal{B}_{3}^0\rangle \label{2.20}\\
	&\qquad\qquad\qquad+\frac35\Big\{2\partial_t+2\fu_\sp^0\cdot\nabla_\sp+\frac{10}{3}\mbox{\rm div}\fu^0\Big\} \bar{p}_k,\nonumber
	\end{align}
	where $\mbox{\rm div}\fu^0:=(\mbox{\rm div}\fu)(t,x_\sp,0)$.
	Once one solves $\bar{u}_{k,\sp}$, then  $(\mathbf{I-P_0})\bar{f}_{k+1}, \bar{u}_{k+1,3}, \bar{p}_{k+1}$ can be determined by the following equations
	\begin{align}\label{2.21}
	(\mathbf{I-P_0})\bar{f}_{k+1}&=\mathbf{L}_0^{-1}\Bigg\{-(\mathbf{I-P_0})(v_3\partial_y\mathbf{P}_0 \bar{f}_k)\nonumber\\
	&\quad+\frac{y}{\sqrt{\mu_0}}\Big[Q(\partial_3\mu_0,\sqrt{\mu_0}\mathbf{P_0}\bar{f}_{k})+Q(\sqrt{\mu_0}\mathbf{P_0}\bar{f}_{k},\partial_3\mu_0)\Big]\nonumber\\
	&\quad+ \frac{1}{\sqrt{\mu_0}} \Big[Q(F_1^0,\sqrt{\mu_0}\mathbf{P_0}\bar{f}_k)+Q(\sqrt{\mu_0}\mathbf{P_0}\bar{f}_k, F_1^0)
	\Big]\nonumber\\
	&\quad+\frac1{\sqrt{\mu_0}} \Big[Q(\sqrt{\mu_0}\bar{f}_1,\sqrt{\mu_0}\mathbf{P_0}\bar{f}_k)+Q(\sqrt{\mu_0}\mathbf{P_0}\bar{f}_k,\sqrt{\mu_0}\bar{f}_1)\Big]\Bigg\}+\bar{J}_{k-1},
	\end{align}
	
	\begin{align}
	\partial_y\bar{u}_{k+1,3}&=-\frac{1}{\rho^0}\Big\{\partial_t\bar{\rho}_{k}+\mbox{\rm div}_\sp(\rho^0\bar{u}_{k,\sp}+\bar{\rho}_{k} \fu^0_\sp)\Big\},\label{2.24}\\
	\partial_y \bar{p}_{k+1}
	&=-\rho^0\partial_t \bar{u}_{k,3}-\rho^0 (\fu^0_\sp\cdot \nabla_\sp)\bar{u}_{k,3}+\rho^0 \partial_3\fu_3^0\bar{u}_{k,3}-\frac43\rho^0\partial_y[(\partial_3\fu_3^0\cdot y+u^0_{1,3})\bar{u}_{k,3}]\nonumber\\
	&\quad+\frac43\mu(T^0) \partial_{yy}\bar{u}_{k,3}-T^0\partial_y\langle\bar{J}_{k-1}
	,\,  \mathcal{A}_{33}^0\rangle+\bar{W}_{k-1,3},\label{2.25}
	\end{align}
	and
	\begin{align}
	\bar{W}_{k-1,i}&=-\sum_{j=1}^2\partial_j \int_{\R^3} T^0 \{\mathbf{I-P_0}\}\bar{f}_k\cdot  \mathcal{A}^0_{i,j} dv,\quad\mbox{for}\  i=1,2,3,\label{2.24-1}\\
	\bar{H}_{k-1}&=-\sum_{j=1}^2\partial_j \Big\{(T^0)^{\frac32}\int_{\R^3} \{\mathbf{I-P_0}\}\bar{f}_k\cdot  \mathcal{B}^0_{j} dv+\sum_{l=1}^22T^0\fu_l^0\int_{\R^3}  \{\mathbf{I-P_0}\}\bar{f}_k \cdot  \mathcal{A}^0_{l,j} dv\Big\}\nonumber\\
	&\qquad\quad -2\fu_\sp^0\cdot \bar{W}_{k-1,\sp}.\label{2.24-2}
	\end{align}
	and
	\begin{align}\label{2.26}
	\bar{J}_{k-1}&=\mathbf{L}_0^{-1}\Bigg\{-(\mathbf{I-P_0})\Big(\frac{(\partial_t+v_\sp\cdot \nabla_\sp)\bar{F}_{k-1}}{\sqrt{\mu_0}} \Big)-(\mathbf{I-P_0})[v_3\partial_y(\mathbf{I-P_0})\bar{f}_k]\nonumber\\
	&\quad+\sum_{\substack{j+l=k+1\\2\leq l\leq \fb,\,  j\geq 1}}\frac{y^l}{l!}\frac{1}{\sqrt{\mu_0}}\Big[Q(\partial_3^l\mu_0,\sqrt{\mu_0}\bar{f}_{j})+Q(\sqrt{\mu_0}\bar{f}_{j},\partial_3^l\mu_0)\Big]\nonumber\\
	&\quad+\sum_{\substack{i+j=k+1\\i\geq 2, j\geq1}} \frac{1}{\sqrt{\mu_0}} \Big[Q(F_i^0,\sqrt{\mu_0}\bar{f}_j)+Q(\sqrt{\mu_0}\bar{f}_j,F_i^0)
	\Big]+\sum_{\substack{i+j=k+1\\i,j\geq2}} \frac{Q(\sqrt{\mu_0}\bar{f}_i,\sqrt{\mu_0}\bar{f}_j)}{\sqrt{\mu_0}}\nonumber\\
	&\quad+\sum_{\substack{i+j+l=k+1\\1\leq l\leq \fb,\, i,j\geq1}}\frac{y^l}{l!}\frac{1}{\sqrt{\mu_0}}\Big[Q(\partial_3^lF^0_i,\sqrt{\mu_0}\bar{f}_{j})+Q(\sqrt{\mu_0}\bar{f}_{j},\partial_3^lF^0_i)\Big]\nonumber\\
	&\quad+\frac{y}{\sqrt{\mu_0}}\Big[Q(\partial_3\mu_0,\sqrt{\mu_0}(\mathbf{I-P_0})\bar{f}_{k})+Q(\sqrt{\mu_0}(\mathbf{I-P_0})\bar{f}_{k},\partial_3\mu_0)\Big]\nonumber\\
	&\quad+ \frac{1}{\sqrt{\mu_0}} \Big[Q(F_1^0,\sqrt{\mu_0}(\mathbf{I-P_0})\bar{f}_k)+Q(\sqrt{\mu_0}(\mathbf{I-P_0})\bar{f}_k, F_1^0)
	\Big]\nonumber\\
	&\quad+\frac1{\sqrt{\mu_0}} \Big[Q(\sqrt{\mu_0}\bar{f}_1,\sqrt{\mu_0}(\mathbf{I-P_0})\bar{f}_k)+Q(\sqrt{\mu_0}(\mathbf{I-P_0})\bar{f}_k,\sqrt{\mu_0}\bar{f}_1)\Big]\Bigg\},
	\end{align}
	We point out that $\bar{W}_{k-1}, \  \bar{H}_{k-1}$ and $\bar{J}_{k-1}$ depend on $\bar{f}_{j}, 1\leq j\leq k-1$.
\end{theorem}

\begin{remark}
We remark that the Taylor expansion of $\mu$ in the derivation of the viscous boundary layer equations is crucial for the control of large velocity $v$ in our $L^2$-$L^\infty$ framework. On the other hand, such an expansion creates a factor of y, which leads to only algebraic decay in $y$ in general. This is in stark contrast to the other standard boundary layer theories which typically have exponential decay in the normal direction.
\end{remark}

\begin{remark}
	For $k=1$, noting $\bar{J}_{0}=\bar{W}_{0}=\bar{H}_{0}=0$ and \eqref{2.14},  then
	the system \eqref{2.19-2}-\eqref{2.20} for $[\bar{u}_{1,\sp}, \bar\theta_1]$ becomes
	\begin{align}
	&\rho^0 \partial_t \bar{u}_{1,i}+\rho^0 (\fu^0_{\sp} \cdot\nabla_\sp) \bar{u}_{1,i}+\rho^0 (\partial_3\fu_3^0\cdot y+u^0_{1,3}) \partial_y \bar{u}_{1,i}\nonumber\\
	&\qquad\qquad\qquad\qquad\qquad+\rho^0 \bar{u}_{1,\sp}\cdot \nabla_\sp \fu^0_i%+\rho^0\partial_y u_{1,3} \cdot \bar{u}^\v_{1,i}
	-\frac{\partial_i p^0}{3T^0} \bar{\theta}_1-\mu(T^0) \partial_{yy} \bar{u}_{1,i}=0,\ i=1,2, \label{2.15}\\
	&\rho^0 \partial_t \bar{\theta}_1+\rho^0 (\fu^0_{\sp} \cdot\nabla_\sp) \bar{\theta}_1+\rho^0 (\partial_3\fu_3^0\cdot y+u^0_{1,3}) \partial_y\bar{\theta}_1+\frac23\rho^0 \mbox{\rm div}\fu^0 \bar{\theta}_1-\frac35\kappa(T^0) \partial_{yy}\bar{\theta}_1=0,\label{2.15-1}
	\end{align}
	which is indeed a linear system for $(\bar{u}_{1,\sp}, \bar{\theta}_1)$.
%Moreover $\bar{u}_{2,3}, \bar{p}_2$ and $(\mathbf{I-P_0})\bar{f}_{2}$ satisfy
%	\begin{align*}%\label{2.18}
%	\partial_y\bar{u}_{2,3}=-\frac{1}{\rho^0}\Big\{\partial_t\bar{\rho}_{1}+\mbox{\rm div}_\sp(\rho^0\bar{u}_{1,\sp}+\bar{\rho}_{1} \fu^0_\sp)\Big\} \quad\mbox{and}\quad
%	\partial_y \bar{p}_{2}=0,
%	\end{align*}
%	and
%	\begin{align*}%\label{2.19-1}
%	&(\mathbf{I-P_0})\bar{f}_{2}=\mathbf{L}_0^{-1}\Bigg\{-v_3\partial_y \bar{f}_1+\frac{y}{\sqrt{\mu_0}}\Big[Q(\partial_3\mu_0,\sqrt{\mu_0}\bar{f}_{1})+Q(\sqrt{\mu_0}\bar{f}_{1}, \partial_3\mu_0)\Big]\nonumber\\
%	&\qquad\quad+ \frac{1}{\sqrt{\mu_0}} \Big[Q(\sqrt{\mu_0}f_1^0,\sqrt{\mu_0}\bar{f}_1)+Q(\sqrt{\mu_0}\bar{f}_1,\sqrt{\mu_0}f_1^0)
%	\Big] +\frac1{\sqrt{\mu_0}} \Big[Q(\sqrt{\mu_0}\bar{f}_1,\sqrt{\mu_0}\bar{f}_1)\Big]\Bigg\}.
%	\end{align*}
As indicated later in Remark \ref{rem1.4},  the $\v$-order viscous boundary layer $\bar{F}_1$ will appear  if one of $\partial_3 \fu^0_{1}(t,x_\sp,0), \ \partial_3 \fu^0_{2}(t,x_\sp,0)$ and $\partial_3 T^0(t,x_\sp,0)$ is nonzero. That means whether the main viscous boundary layer $\bar{F}_1$ appears depends only   on the boundary  properties of compressible Euler solution, and has no relation with the interior expansion $F_1$.
\end{remark}

\subsubsection{ Knudsen boundary layer expansion:}
To construct the solution that satisfies the boundary condition at higher orders, we have to introduce the Knudsen boundary layer. Firstly, we define the new scaled normal coordinate:
\begin{equation}
\eta:=\frac{x_3}{\v^2}.
\end{equation}
The  Knudsen boundary layer expansion is defined as
\begin{equation}\nonumber
\hat{F}^\v(t,x_\sp, \eta)\sim \sum_{k=1}^\infty \v^k \hat{F}_k(t,x_\sp, \eta,v).
\end{equation}
Plugging $F^\v+\bar{F}^\v+\hat{F}^\v$  into \eqref{1.1} and comparing the order of $\v$,  then using \eqref{1.7-1}, \eqref{1.14}, one obtains
\begin{equation}\label{1.19-1}
\begin{split}
\dis& \v^{-1}:  \quad v_3 \frac{\partial\hat{F}_1}{\partial \eta }-\big[Q(\mu_0,\hat{F}_{1})+Q(\hat{F}_1,\mu_0)\big]=0,\\[2mm]
\dis & \v^0:\quad \  v_3 \frac{\partial\hat{F}_2}{\partial \eta }-\big[ Q(\mu_0,\hat{F}_{2})+Q(\hat{F}_{2},\mu_0) \big]\\
&\qquad\quad= Q(F_1^0+\bar{F}^0_1,\hat{F}_1)+Q(\hat{F}_1,F_1^0+\bar{F}^0_1)
+Q(\hat{F}_1,\hat{F}_1),\\[2mm]
&\dis \quad\quad\quad\quad\quad\quad\quad\ \vdots \\
\dis & \v^{k}:\quad   v_3\frac{\partial\hat{F}_{k+2}}{\partial \eta }-\big[ Q(\mu_0,\hat{F}_{k+2})+Q(\hat{F}_{k+2},\mu_0) \big] \\
&\qquad\quad=-\{\partial_t+v_\sp\cdot\nabla_\sp\}\hat{F}_{k}+\sum_{\substack{j+2l=k+2\\  1\leq l\leq \fb, j\geq1}} \frac{\eta^l}{l!} \big[Q(\partial_3^l\mu_0, \hat{F}_{j})+Q( \hat{F}_{j},\partial_3^l\mu_0) \big] \\
&\qquad\qquad+\sum_{\substack{i+j=k+2\\ i,j\geq1}}\big[Q(F_i^0+\bar{F}_i^0,\hat{F}_j)+Q(\hat{F}_j,F_i^0+\bar{F}_i^0)+Q(\hat{F}_i,\hat{F}_j)\big] \\
&\qquad\qquad+\sum_{\substack{i+j+2l=k+2\\ i,j\geq1, 1\leq l\leq \fb}}\frac{\eta^l}{l!} \big[Q(\partial_3^l F_i^0, \hat{F}_{j})+Q( \hat{F}_{j},\partial_3^l F_i^0) \big]   \\
&\qquad\qquad+\sum_{\substack{i+j+l=k+2\\ i,j\geq1, 1\leq l\leq \fb}}\frac{\eta^l}{l!} \big[Q(\partial_y^l \bar{F}_i^0, \hat{F}_{j})+Q( \hat{F}_{j},\partial_y^l \bar{F}_i^0) \big],
\quad \mbox{for}\  k\geq 1,
\end{split}
\end{equation}
where we have used \eqref{1.14-1}-\eqref{1.14-2} and  the Taylor expansion of $\bar{F}_i$
\begin{equation*}%\label{1.14-4}
\bar{F}_i(t, x_1, x_2, y,v)
=F_i^0+\sum_{l=1}^{\fb} \frac{1}{l!} \partial_y^l \bar{F}_i^0\cdot y^l + \frac{y^{\fb+1}}{(\fb+1)!} \partial_3^{\fb+1}\bar{\mathfrak{F}}_i,
\end{equation*}
with
\begin{align}\label{1.14-5}
\partial_3^l\bar{F}_i^0:&=(\partial_3^l\bar{F}_i)(t,x_1, x_2, 0, v),
\quad \partial_3^{\fb+1}\bar{\mathfrak{F}}_i:=(\partial_3^{\fb+1}\bar{F}_i)(t,x_1, x_2, \bar{\xi}_i,v),
\, \mbox{for}\,\, 0\leq l\leq \fb,
\end{align}
for some $\bar{\xi}_i\in [0,y]$.

It is noted that the Knudsen boundary layer  \eqref{1.19-1} is in fact a steady problem with $(t,x_\sp)$ as parameters, and the well-posedness has already been obtained in \cite{GPS-1988,Jiang-Wang} under some conditions on the source term and boundary condition.  However, we shall use the existence results in  \cite{Jiang-Wang} since the continuity and uniform estimate in  $L^\infty_{\eta,v}$ is needed in the present paper.
%As in \cite{GPS-1988},one can prove that the Knudsen layers decay exponentially as $\eta\rightarrow\infty$, but the solution of \cite{GPS-1988} is not in the space $L^\infty_{y,v}$.

\subsection{Hilbert Expansion}
Now we consider the Boltzmann solution with the following multiscale Hilbert expansion
\begin{align}\label{1.22}
F^\v=\mu+\sum_{i=1}^{N} \v^i F_i(t,x,v)+\sum_{i=1}^{N} \v^i\bar{F}_i(t,x_\sp,\frac{x_3}{\v},v)+\sum_{i=1}^{N} \v^i\hat{F}_i(t,x_\sp,\frac{x_3}{\v^2},v)+\v^5 F^\v_{R},
\end{align}
which, together with \eqref{1.7-1}, \eqref{1.14} and \eqref{1.19-1}, yields that
the equation of reminder  $F_R^{\v}$ is given by
\begin{align}\label{1.9-2}
&\pt F^{\v}_{R}+v\cdot \nabla_x F^{\v}_R-\frac{1}{\v^2}\{Q(\mu, \FRep)+Q(\FRep, \mu)\}\nonumber\\
&=\v^3Q(\FRep, \FRep)+\sum_{i=1}^N \v^{i-2}\{Q(F_i+\bar{F}_i+\hat{F}_i,F^\v_R)+Q(F^\v_R, F_i+\bar{F}_i+\hat{F}_i)\}\nonumber\\
&\quad+ R^\v+\bar{R}^\v+\hat{R}^\v,
\end{align}
where
\begin{align}\label{1.9-3}
R^\v=&-\v^{N-6}\{\partial_t+v\cdot \nabla_x\} (F_{N-1}+\v F_n)+\v^{N-6}\sum_{\substack{i+j\geq N+1\\1\leq i,j\leq N}}\v^{i+j-N-1}Q(F_{i},F_{j}),
\end{align}
\begin{align}\label{1.9-4}
\bar{R}^\v&=-\v^{N-6}\{\partial_t+v_\sp\cdot \nabla_\sp\} (\bar{F}_{N-1}+\v \bar{F}_N)-\v^{N-6} v_3\partial_y \bar{F}_N\nonumber\\
&\quad+\v^{N-6}\sum_{\substack{j+l\geq N+1\\ 1\leq j\leq N, \,  1\leq l\leq \fb}} \v^{l+j-N-1}\cdot\frac{y^l}{l!} \big[Q(\partial_3^l\mu_0, \bar{F}_{j})+Q( \bar{F}_{j},\partial_3^l\mu_0) \big]\nonumber\\
&\quad+\v^{N-6}\sum_{\substack{i+j\geq N+1\\   1\leq i,j\leq N}} \v^{i+j-N-1} \big[ Q(F_i^0,\bar{F}_j)+Q(\bar{F}_j,F_i^0)+Q(\bar{F}_i,\bar{F}_j) \big]\nonumber\\
&\quad+\v^{N-6}\sum_{\substack{i+j+l\geq N+1\\ 1\leq i, j\leq N, \,  1\leq l\leq \fb}} \v^{i+j+l-N-1}\cdot \frac{y^l}{l!} \big[Q(\partial_3^l F_i^0, \bar{F}_{j})+Q( \bar{F}_{j},\partial_3^l F_i^0) \big]\nonumber\\
&\quad +  \v^{\fb-5}\frac{y^{\fb+1}}{(\fb+1)!} \sum_{j=1}^N \v^{j-1} [Q(\partial_3^{\fb+1}\tilde{\mu}, \bar{F}_{j})+Q( \bar{F}_{j},\partial_3^{\fb+1}\tilde{\mu})]\nonumber\\
&\quad +\v^{\fb-4}\frac{y^{\fb+1}}{(\fb+1)!} \sum_{i,j=1}^N \v^{i+j-2} \big[Q(\partial_3^{\fb+1}\mathfrak{F}_i, \bar{F}_{j})+Q( \bar{F}_{j},\partial_3^{\fb+1}\mathfrak{F}_i)\big],
\end{align}
and
\begin{align}\label{1.9-5}
\hat{R}^\v&=-\v^{N-6}\{\partial_t+v_\sp\cdot \nabla_\sp\} (\hat{F}_{N-1}+\v \hat{F}_N)\nonumber\\
&\quad+\v^{N-6}\sum_{\substack{j+2l\geq N+1\\  1\leq j\leq N, 1\leq l\leq \fb}} \v^{j+2l-N-1}\cdot \frac{\eta^l}{l!} \big[ Q(\partial_3^l\mu_0, \hat{F}_{j})+Q( \hat{F}_{j},\partial_3^l\mu_0) \big]\nonumber\\
&\quad+\v^{N-6}\sum_{\substack{i+j\geq N+1\\  1\leq i,j\leq N}} \v^{i+j-N-1} \big[Q(F_i^0+\bar{F}_i^0,\hat{F}_j)+Q(\hat{F}_j,F_i^0+\bar{F}_i^0)+Q(\hat{F}_i,\hat{F}_j)\big]\nonumber\\
&\quad+\v^{N-6}\sum_{\substack{i+j+2l\geq N+1\\  1\leq i,j\leq N, 1\leq l\leq \fb}} \v^{i+j+2l-N-1}\cdot \frac{\eta^l}{l!} \big[ Q(\partial_3^l F_i^0, \hat{F}_{j})+Q( \hat{F}_{j},\partial_3^l F_i^0) \big]\nonumber\\
&\quad+\v^{N-6}\sum_{\substack{i+j+l\geq N+1\\  1\leq i,j\leq N, 1\leq l\leq \fb}} \v^{i+j+l-N-1}\cdot \frac{\eta^l}{l!} \big[ Q(\partial_y^l \bar{F}_i^0, \hat{F}_{j})+Q( \hat{F}_{j},\partial_y^l \bar{F}_i^0) \big]\nonumber\\
&\quad + \v^{2\fb -4} \frac{\eta^{\fb+1}}{(\fb+1)!} \sum_{j=1}^N \v^{j-1} \big[ Q(\partial_3^{\fb+1}\tilde{\mu}, \hat{F}_{j})+Q( \hat{F}_{j},\partial_3^{\fb+1}\tilde{\mu}) \big]\nonumber\\
&\quad + \v^{2\fb-3}\frac{\eta^{\fb+1}}{(\fb+1)!} \sum_{i,j=1}^N \v^{i+j-2} \big[ Q(\partial_3^{\fb+1}\mathfrak{F}_i, \hat{F}_{j})+Q( \hat{F}_{j},\partial_3^{\fb+1}\mathfrak{F}_i) \big]\nonumber\\
&\quad + \v^{\fb-4}\frac{\eta^{\fb+1}}{(\fb+1)!} \sum_{i,j=1}^N \v^{i+j-2} \big[ Q(\partial_3^{\fb+1}\bar{\mathfrak{F}}_i, \hat{F}_{j})+Q( \hat{F}_{j},\partial_3^{\fb+1}\bar{\mathfrak{F}}_i) \big],
\end{align}
where $\pa_3^l\mu_0, \pa_3^{\fb+1}\tilde{\mu}$, $\pa_3^lF_i^0, \pa_3^{\fb+1}\mathfrak{F}_i$ and $\pa_y^l\bar{F}_i^0, \pa_y^{\fb+1}\bar{\mathfrak{F}}_i$ are the ones defined in \eqref{1.14-3}, \eqref{1.14-5}.

\

The main aim of the present paper is to  establish the validity of the Hilbert expansion for the Boltzmann equation around  the local Maxwellian $\mu$ determined by the compressible Euler equations \eqref{1.7}, so it is natural to rewrite the remainder as
\begin{equation}\label{1.20}
F^{\v}_{R}=\sqrt{\mu} f^\v_R.
\end{equation}
To use the $L^2$-$L^\infty$ framework \cite{Guo Jang Jiang-1}, we  introduce a global Maxwellian
	\begin{equation*}%\label{1.28}
		\mu_{M}:=\frac{1}{(2\pi T_{M})^{3/2}}\exp{\left\{-\frac{|v|^{2}}{2T_{M}}\right\}},
	\end{equation*}
	where $T_{M}>0$ satisfies the condition
	\begin{equation}\label{1.29}
		T_{M}<\min_{x\in\mathbb{R}^3_+} T(t,x)\leq \max_{x\in\mathbb{R}^3_+} T(t,x)<2T_{M}.
	\end{equation}
By the assumption \eqref{1.29}, one can easily deduce that there exists positive constant $C>0$ such that for some $\frac{1}{2}<\alpha<1$ , the following holds:
\begin{equation}\label{relation of mu and muM}
\frac1C\mu_{M}\leq \mu(t,x,v)\leq C\mu_{M}^{\alpha}.
\end{equation}
We further define
\begin{equation}\label{def of h}
F_{R}^{\v}=\{1+|v|^{2}\}^{-\frac{\k}{2}}\sqrt{\mu_{M}} h^\v_R\equiv\frac{1}{w_{\k}(v)}\sqrt{\mu_{M}}\hep_R,
\end{equation}
with the velocity weight function
\begin{equation*}
w_{\k}(v):=\{1+|v|^{2}\}^{\frac{\k}{2}},\quad \mbox{for} \quad \k\geq0.
\end{equation*}
 %We also define the exponential growth weight function
% \begin{align}\label{ewf}
%\tw_{\kappa}(v)= w_{\k}(v) \, e^{\fa |v|^2},
% \end{align}
%for any $\k\geq0$ and $0\leq \fa<\frac{\alpha}{2}$, where $\frac12<\alpha<1$ is the one fixed in \eqref{relation of mu and muM}.

\begin{theorem}\label{theorem}
Let $\tau>0$ be the life-span of smooth solution of compressible Euler equations \eqref{1.7}. Let  $\k\geq 7$,   $N\geq 6$ and $\fb\geq 5$.  We assume the initial data
\begin{equation}%\label{1.35}
F^{\v}(0,x,v)=\mu(0,x,v)+\sum_{i=1}^{N}\v^i \left\{F_i(0,x,v)+ \bar{F}_i(0,x_\sp, \frac{x_3}{\v},v)+ \hat{F}_i(0,x_\sp,\frac{x_3}{\v^2},v)\right\}+\v^5F^\v_R(0,x,v)\geq 0,\nonumber
\end{equation}
and $F_i(0), \bar{F}_i(0), \hat{F}_i(0),  i=1,\cdots, N$  satisfy the regularity and compatibility condition described in  Proposition \ref{prop5.1},  and
\begin{align*}%\label{1.36}
\Big\|(\frac{F^\v_R}{\sqrt{\mu}})(0)\Big\|_{L^2_{x,v}}+\v^3\Big\|(w_\k \frac{F^\v_R}{\sqrt{\mu_M}})(0)\Big\|_{L^\infty_{x,v}}
<\infty.
 \end{align*}
Then there exists a small positive constants $\v_0>0$  such that IBVP of Boltzmann equation \eqref{1.1}, \eqref{1.5}  has a unique solution for  $\v\in(0,\v_0]$ over the time interval $t\in [0,\tau]$ in the following form of expansion
\begin{equation}\label{1.36-1}
F^{\v}(t,x,v)=\mu(t,x,v)+\sum_{i=1}^{N}\v^i \left\{F_i(t,x,v)+ \bar{F}_i^\v(t,x_\sp, \frac{x_3}{\v},v)+ \hat{F}_i^\v(t,x_\sp,\frac{x_3}{\v^2},v)\right\}+\v^5F^\v_R(t,x,v)\geq 0
\end{equation}
with
\begin{align}
&\sup_{t\in[0,\tau]} \left\{\Big\|\frac{F^\v_{R}(t)}{\sqrt{\mu}}\Big\|_{L^2_{x,v}}+\v^3\Big\|w_\k(v)\frac{F^\v_{R}(t)}{\sqrt{\mu_M}}\Big\|_{L^\infty_{x,v}}\right\}\leq C(\tau)<\infty.\label{1.39}%\\
%& \sup_{t\in[0,\tau]} \sum_{i=1}^n\left\{\Big\|\tw_{\k_i}\frac{F_i(t)}{\sqrt{\mu_M}}\Big\|_{L^\infty_{x,v}}+\Big\|\tw_{\bar{\k}_i}\frac{\bar{F}_i(t)}{\sqrt{\mu_M}}\Big\|_{L^\infty_{x,v}}+\Big\|\tw_{\hat{\k}_i}\frac{\hat{F}_i(t)}{\sqrt{\mu_M}}\Big\|_{L^\infty_{x,v}}\right\}\nonumber\\
%&+\|??\|_{L^2}\leq C(\tau)<\infty.\label{1.40}
\end{align}
Here the functions $F_i(t,x,v), \bar{F}_i(t,x_{\sp},y,v)$ and $\hat{F}_i(t,x_{\sp},\eta,v)$ are  respectively the interior expansion, viscous and Knudsen boundary layers  constructed in Proposition \ref{prop5.1}.
%The positive constants $\k_i, \bar{\k}_i, \hat{\k}_i$ are chosen such that  $\k_i\gg \bar{\k}_i\gg\hat{\k}_i\gg \k_{i+1}\gg \bar{\k}_{i+1}\gg\hat{\k}_{i+1}\gg1$.
\end{theorem}

\begin{remark}
From  \eqref{1.36-1}-\eqref{1.39} and the uniform estimates in  Proposition \ref{prop5.1}, it is direct to check that
	\begin{align}\nonumber
	\sup_{t\in[0,\tau]}\left\{ \Big\|\Big(\frac{F^\v-\mu}{\sqrt{\mu}}\Big)(t)\Big\|_{L^2(\R_+^3\times\R^3)}+ \Big\|w_\k \, \Big(\frac{F^\v-\mu}{\sqrt{\mu_M}}\Big)(t)\Big\|_{L^\infty(\R_+^3\times\R^3)}
	\right\}\leq C \v\rightarrow 0.
	\end{align}
Hence we have established  the hydrodynamic limit from the Boltzmann equation to the compressible Euler system for the half-space  problem.
\end{remark}

\begin{remark}
For the initial data $F_i(0,x,v), \bar{F}_i(0,x_\sp,y,v)$ and $\hat{F}_i(0,x_\sp,\eta,v)$,
we only need to  impose data on the macroscopic part of  $F_i(0,x,v)$, and the part of the macroscopic part of viscous boundary layer $\bar{F}_i(0,x_\sp,y,v)$, and no conditions are needed on $\hat{F}_i(0,x_\sp,\eta,v)$, see Proposition \ref{prop5.1} for more details. Also, from Proposition \ref{prop5.1}, we know that $ \bar{F}_i$ decay algebraically with respect to  $y$, and $\hat{F}_i$ decay exponentially with respect to $\eta$; and the decay estimates are crucial for us to close the estimate for $\bar{R}^\v, \hat{R}^\v$.
\end{remark}

\begin{remark}\label{rem1.4}
 For the first order of viscous boundary layer $\bar{F}_1$, its boundary condition is close related to the boundary value of compressible Euler solution, i.e., (see \eqref{7.9} for details)
\begin{align}\label{1.32}
\begin{cases}
\dis\partial_y\bar{u}_{1,i}(t,x_\sp,y)|_{y=0}=-\partial_3 \fu^0_{i}(t,x_\sp,0),\quad i=1,2,\\[2mm]
\dis\partial_y\bar{\theta}(t,x_\sp,y)|_{y=0}=-3\partial_3 T^0(t,x_\sp,0).
\end{cases}
\end{align}
It is also noted that $\bar{F}_1$ does not interplay with $F_1$, see \eqref{2.14}, \eqref{2.15}-\eqref{2.15-1} and \eqref{2.51}.
Hence, if one of $\partial_3 \fu^0_{1}(t,x_\sp,0), \ \partial_3 \fu^0_{2}(t,x_\sp,0)$ and $\partial_3 T^0(t,x_\sp,0)$ is nonzero, then the viscous boundary layer $\bar{F}_1$ must be  nonzero. That means the  $\v$-order viscous boundary layer  $\bar{F}_1$ is imperative, and must be included in the Hilbert expansion.  On the other hand, the first order of Knudsen boundary layer $\hat{F}_1$ does not appear, i.e.,  $\hat{F}_1\equiv0$ (see  \eqref{7.11} for  details), and this is reasonable since the Knudsen boundary layer is  used to mend the boundary condition at higher orders.
Therefore the interplay of interior expansion, viscous and Knudsen boundary layers start from $\v^2$-order.
\end{remark}

%\begin{remark}
%We do not describe the precise relations between $\k_i, \bar{\k}_i$ and $\hat{\k}_i$ because the functions  $F_i, \bar{F}_i$ and $\hat{F}_i$ indeed decay exponentially with respect to particle velocity.
%\end{remark}

\begin{remark}
For non-flat domain, as pointed out in \cite{Wu},  one  needs to modify the equation of expansion for boundary layers due to the non-trivial geometry of domain. Indeed it is still a very important open problem for general smooth domain.
\end{remark}

\subsection{Acoustic Limit}
The acoustic system is the linearization of  compressible Euler equations around a uniform
fluid state, for instance,   $(1,0,1)$. After a suitable choice of units, the fluid fluctuations $(\varphi, \Phi, \vartheta)=(\varphi, \Phi_1,\Phi_2,\Phi_3, \vartheta)$ satisfies
\begin{align}\label{1.9}
\begin{cases}
\partial_t \varphi+\mbox{\rm div} \Phi=0,\\
\partial_t \Phi+\nabla(\varphi+\vartheta)=0,\\
\frac32\partial_t\vartheta+\mbox{\rm div} \Phi=0,
\end{cases}
(t, x)\in\R_+\times\R^3_+.
\end{align}
We impose \eqref{1.9} with the following initial and boundary data
\begin{align}\label{1.9-1}
(\varphi, \Phi, \vartheta)(0,x)=(\varphi_0, \Phi_0, \vartheta_0)(x)\in H^{s_0}(\R^3_+),\quad\mbox{and}\quad
\Phi_3(t,x)|_{x_3=0}=0.
\end{align}
It is direct to know that the IBVP \eqref{1.9}-\eqref{1.9-1} is a linear hyperbolic system with constant coefficients and  characteristic boundary, and there exists a unique global smooth solution $(\varphi, \Phi, \vartheta)(t)\in H^{s_0}(\R^3_+)$. In fact, we can also use Lemma \ref{lem3.1} (taking the Euler solution $(\rho,\fu, T)=(1,0,1)$ in \eqref{5.1},  even though the coefficients are slightly different, but Lemma \ref{lem3.1} is still valid) to obtain the global existence of smooth solution to IBVP \eqref{1.9}-\eqref{1.9-1}. Moreover it holds that
\begin{align*}%\label{1.9-10}
\sup_{s\in[0,t]}\|(\varphi, \Phi, \vartheta)(s)\|_{H^{s_0}(\R^3_+)}\leq C\big(t, \|(\varphi_0, \Phi_0, \vartheta_0)\|_{H^{s_0}(\R^3_+)}\big),\quad \forall \, t>0.
\end{align*}

On the other hand, the acoustic system \eqref{1.9} can also be formally derived from the Boltzmann equation \eqref{1.1} by letting
\begin{equation}\label{1.11}
F^\v=\tilde{\mu}_M+\d G^\v,
\end{equation}
where $\tilde{\mu}_M$ is the global Maxwellian determined by the uniform state $(1,0,1)$, i.e.,
\begin{equation*}
\tilde{\mu}_M=\frac{1}{(2\pi)^{\frac32}} \exp(-\frac{|v|^2}{2}).
\end{equation*}
The fluctuation amplitude $\delta$ is a function of $\v$ satisfying:
\begin{equation}\nonumber
\delta\rightarrow 0  \quad \mbox{as}\quad \v\rightarrow0.
\end{equation}
For instance we can take
\begin{equation}\nonumber
\delta=\v^{\varpi} \quad \mbox{for}\quad \varpi>0.
\end{equation}
With the above scalings, $G^\v$ formally converges to
\begin{align}\label{1.15}
G:= \left\{\varphi+v\cdot \Phi+\frac{|v|^2-3}{2}\vartheta  \right\} \tilde{\mu}_M, \quad \mbox{as}\ \v\rightarrow0,
\end{align}
where
$(\varphi, \Phi, \vartheta)$ is the solution of acoustic system \eqref{1.9},  see \cite{Bardos,Bardos-3} for detailed formal derivation.

One of the purpose of present paper is to establish the acoustic limit for initial boundary value problem of Boltzmann equation over half-space $\mathbb{R}^3_+$. %via a recent $L^2-L^\infty$ framework.
We use $\delta$ to denote the fluctuation amplitude and assume that
\begin{equation*}%\label{1.21}
\frac{\v}{\delta}\rightarrow0\quad \mbox{as}\quad  \v\rightarrow0.
\end{equation*}
\begin{theorem}\label{thm5.1}
Let $\tau>0$ be any given time.   Let $\mu^{\delta}(0,x,v)$ be the local Maxwellian with initial datum $1+\delta \varphi_0, \delta \Phi_0$ and $1+\delta \vartheta_0$:
	\begin{align}\nonumber
	\mu^{\delta}(0,x,v)=\frac{1+\delta \varphi_0(x)}{[2\pi(1+\delta \vartheta_0(x))]^{\frac32}}
	\exp{\left\{-\frac{|v-\delta \Phi_0(x)|^2}{2(1+\delta \vartheta_0(x))}\right\}},
	\end{align}
where	$(\varphi_0,\Phi_0,\vartheta_0)$ is the initial data given in \eqref{1.9-1}.
We assume that the conditions in Theorem \ref{theorem} hold, and  rewrite the corresponding Hilbert expansion established in Theorem \ref{theorem} as
\begin{align}
F^\v(t,x,v)&=\mu^{\delta}(t,x,v)+\sum_{i=1}^{N}\v^i\{F_i(t,x,v)+\bar{F}_i(t,x_\sp,\frac{x_3}{\v},v)+\hat{F}_i(t,x_\sp,\frac{x_3}{\v^2},v) \}\nonumber\\
&\qquad+\v^5 F^\v_R(t,x,v)\geq0.\nonumber
\end{align}
Then there exist $\v_0>0$ and $\delta_0>0$ such that for each $\v\in(0,\v_0]$ and $\d\in(0,\delta_0]$, there exists a constant $C>0$ so that
\begin{equation*}
\sup_{0\leq t\leq \tau}\Big\{\|G^\v(t)-G(t)\|_{L^\infty(\R_+^3\times\R^3)}+\|G^\v(t)-G(t)\|_{L^2(\R_+^3\times\R^3)}\Big\}\leq C\{\delta+\frac{\v}{\delta}\},
\end{equation*}
where $\dis\frac{\v}{\delta}\rightarrow0\,\, \mbox{as}\  \v\rightarrow0$, and $G^\v$ and $G$ are defined in \eqref{1.11} and \eqref{1.15} respectively. The constant $C>0$ here depends only on $\tau$ and initial data $\|(\varphi_0,\Phi_0,\vartheta_0)\|_{H^{s_0}(\R^3_+)}$.
\end{theorem}

%\smallskip

We now  briefly comment on the analysis of the present  paper. For  the Hilbert expansion of Boltzmann equation over the half-space $x\in \R_+^3$ with specular reflection boundary conditions, in general, the viscous and Knudsen  boundary layers will appear. To solve the interior expansion, viscous and Knudsen boundary layers, we need to determine the boundary conditions so that each of them is well-posed. We notice that the Knudsen boundary layer  \eqref{1.19-1} is indeed a steady problem with $(t,x_\sp)\in [0,\tau]\times \R^2$ as parameters. From \cite{Jiang-Wang}, we know that Knudsen boundary layer is well-posed in weighted  $L^\infty_{\eta,v}$-space  under the conditions \eqref{2.30}-\eqref{2.31}, see Lemma \ref{lem2.6} for details. Especially, the weight with respect to normal variable $\eta$ can grow exponentially, which is important for us to close the remainder estimate.  In general, the source term on the right hand side of \eqref{1.19-1}   does not satisfy \eqref{2.30}.  %and the normal regularity can not be obtained.
And we use the idea of \cite{BCN-1986} to introduce a correct function $f_{i,1}$ to  overcome this difficulty, see Lemma \ref{lem2.7} for details. To determine the boundary conditions for $F_i, \bar{F}_i, \hat{F}_i,\, i=1,\cdots,N$, we hope each $F_i+\bar{F}_i+\hat{F}_i$ satisfies the specular reflection boundary conditions, which  together with \eqref{2.31}, we can finally obtain the boundary conditions for interior expansion and viscous boundary layer, see section \ref{sec2.4} for details. Here we point out that $F_i, \bar{F}_i$ and $\hat{F}_i$ may not satisfy the specular reflection boundary conditions alone.

For the existence of  interior expansion, we have to consider a linear  hyperbolic system with characteristic boundary, see Lemma \ref{lem2.1}. We can indeed find some local existence of smooth solution for the linear hyperbolic problem, but without good enough energy estimate, so we resolve this hyperbolic system in Lemma \ref{lem3.1} by energy estimate and careful boundary analysis.
For the existence of viscous boundary  layer, we note that it is involved to a linear hyperbolic system with partial viscosity (only in the normal direction) and linear growth coefficients, which is not a standard linear parabolic system, see Theorem \ref{thm2.5}. By using the energy estimate and  several cut-off approximate arguments,  we establish its well-posedness in a weighted Sobolev space with algebraically growth weight of $y$, see Lemma \ref{lem4.1} for details.

With above preparations, then we can establish the well-posedness of $F_i, \bar{F}_i$ and $\hat{F}_i, i=1,\cdots, N$  and obtain the uniform estimate, see Proposition \ref{prop5.1}.
Now with the help of uniform estimates in Proposition \ref{prop5.1},  we can use the $L^2$-$L^\infty$ framework \cite{Guo Jang Jiang,Guo2010} to obtain the uniform estimate for remainder term $F^\v_R$, and hence obtain the solution of Boltzmann equation in the form of \eqref{1.36-1}.

\smallskip

The paper is organized as follows. In Section \ref{section2},  we reformulate  the interior expansion  and Knudsen boundary layers, and give the proof of Theorem \ref{thm2.5}.  Also we derive the corresponding boundary conditions so that the formulations of  interior expansion, the viscous and Knudsen boundary layers are all well-posed. Section \ref{section3} is devoted to an existence theory for a linear hyperbolic system with characteristic boundary, which are used to construct interior expansion $F_i$. In Section \ref{section4}, to construct the existence of viscous boundary layer $\bar{F}_i$, we establish an existence theory of  IBVP  for a  linear parabolic system with degenerate viscosity and linear growth coefficients in a weighted Sobolev space. In Section \ref{section5}, we construct solutions of  interior expansion, the viscous and Knudsen boundary layers. Finally we prove Theorems \ref{theorem} and \ref{thm5.1} in Sections \ref{section6} and \ref{section7} respectively. In the appendix, we present some anisotropic trace estimates.

\

\noindent{\bf Notations.}  Throughout this paper, $C$ denotes a generic positive constant  and  vary from line to line. And $C(a),C(b),\cdots$ denote the generic positive constants depending on $a,~b,\cdots$, respectively, which also may vary from line to line. We use $\langle\cdot ,\cdot \rangle$ to denote the standard $L^2$ inner product in $\R^3_v$.
% while we use $(\cdot ,\cdot )$ to denote the $L^2$ inner product
$\|\cdot\|_{L^2}$ denotes the standard $L^2(\mathbb{R}^3_+\times\mathbb{R}^3_v)$-norm, and $\|\cdot\|_{L^\infty}$ denotes the $L^\infty(\mathbb{R}^3_+\times\mathbb{R}^3_v)$-norm. \vspace{1.5mm}

%%%%%%%%%%%%%%%%%%%%%%%%%%%%%%%%%%%%%%%%%%%%%%%%%%%%%%%%%%%%%%%%%%%%%%%%%%%%%%%%%%%%%	
\section{Reformulations of  Expansions and Boundary Conditions}\label{section2}

\subsection{Reformulation of Interior Expansion} Firstly we introduce the existence result on the compressible Euler equations.
\begin{lemma}\label{lem2.1-1}
Let $s_0\geq 3$ be some positive integer.  Consider the IBVP of compressible Euler equations  \eqref{1.7}-\eqref{1.12-2}.
Choose $\delta_1>0$ so that for any $\delta\in(0,\delta_1]$, the positivity of $\rho_0$ and $T_0$ is guaranteed. Then for  $\delta\in(0,\delta_1]$, there is a family of classical solutions $(\rho^{\delta},\fu^{\delta}, T^{\delta})\in C([0,\tau^{\delta}]; H^{s_0}(\mathbb{R}^3_+)) \cap C^1([0,\tau^{\delta}]; H^{s_0-1}(\mathbb{R}^3_+))$  of IBVP \eqref{1.7}-\eqref{1.12-2} such that $\rho^{\delta}>0$ and $T^{\delta}>0$, and the following estimate holds:
\begin{equation}\label{2.2-1}
\|(\rho^{\delta}-1,\fu^{\delta}, T^{\delta}-1)\|_{C([0,\tau^{\delta}]; H^{s_0}(\mathbb{R}^3_+)) \cap C^1([0,\tau^{\delta}]; H^{s_0-1}(\mathbb{R}^3_+))}\leq C_0.
\end{equation}
The life-span $\tau^{\delta}$ have the following lower bound
\begin{equation}\label{2.2-2}
\tau^\delta\geq \frac{C_1}{\delta}.
\end{equation}
%Moreover, for any given time $\tau \in(0,\tau^\delta)$, if $\delta$ can be further small (the smallness depends on $\tau_0$), then it holds
%\begin{equation}\label{2.2-2}
%\|(\rho^{\delta}-1,\fu^{\delta}, T^{\delta}-1)\|_{C([0,\tau]; H^{s_0}(\mathbb{R}^3_+)) \cap C^1([0,\tau]; H^{s_0-1}(\mathbb{R}^3_+))} \ll 1.
%\end{equation}
The constant $C_0, C_1$ are independent of $\delta$, depending only on the $H^{s_0}$-norm of $(\varphi_0,\Phi_0,\vartheta_0)$.
\end{lemma}

We refer \cite{Schochet} for the local existence of the IBVP of compressible Euler equation \eqref{1.7}-\eqref{1.12-2}, see also \cite{Chenshuxing} and the references cited therein. We point out that the local existence result in \cite{Schochet} is for smooth bounded domain with $C^\infty$ boundary, but the method can also be applied to our half-space problem. On the other hand,  we can obtain \eqref{2.2-1}-\eqref{2.2-2} by using similar  arguments as in Lemma \ref{lem3.1} below.

Throughout this paper, we will drop the superscript of $(\rho^{\delta},\fu^{\delta}, T^{\delta})$ when no confusion arises. To  derive the estimates of interior expansion, i.e., $F_1(t,x,v),\cdots, F_N(t,x,v)$,
 we  firstly  present a useful lemma which will be used to estimate the bound of linear terms.
\begin{lemma}[\cite{Guo Jang}]\label{lem2.1}
Let $(\rho,\fu, T)(t)$ be some smooth solution of compressible Euler equations \eqref{1.7}.  For each given nonnegative integer $k$, assume $F_k$'s are found. Then the microscopic part of $F_{k+1}$ is determined through the equation for $F_k$ in  \eqref{1.7-1}:
\begin{equation}\label{2.1}
\dis \{\mathbf{I-P}\}(\frac{F_{k+1}}{\sqrt{\mu}})=\mathbf{L}^{-1}\left(-\f{\{\partial_t+v\cdot\nabla_{x}\}F_{k-1}-\sum_{\substack{i+j=k+1\\i,j\geq 1}}Q(F_i,F_j)}{\sqrt{\mu}}\right), \,\,\mbox{for}\,\, k\geq0.
\end{equation}
where we define $F_{-1}=0$ for the consistency of notation.
For the macroscopic part, $\r_{k+1},u_{k+1},\t_{k+1}$ satisfy the following:
\begin{align}\label{2.2}
\begin{split}
\dis \pa_t\r_{k+1}+\mbox{\rm div}_x(\rho u_{k+1}+\r_{k+1}\fu)&=0,\\[1.5mm]
\dis \rho \Big\{\partial_t u_{k+1}+u_{k+1}\cdot \nabla_x \fu+\fu\cdot\nabla_x u_{k+1}\Big\}-\f{\r_{k+1}}{\rho}\nabla_x(\rho T)+\nabla_x(\f{\r\t_{k+1}+3T\r_{k+1}}{3})&=\mathfrak{f}_{k},\\[3mm]
\dis \rho\Big\{\pa_t\t_{k+1}+\f23(\t_{k+1}\mbox{\rm div}_x\fu+3T\mbox{\rm div}_xu_{k+1})+\fu\cdot\nabla_x\t_{k+1}+3u_{k+1}\cdot\nabla_xT\Big\}&=\mathfrak{g}_k,
\end{split}
\end{align}
with
	\begin{align}\label{2.3}
	\begin{split}
	\mathfrak{f}_{k,i}&=-\sum_{j=1}^3\partial_{x_j}\left(T\int_{\mathbb{R}^3}\mathcal{A}_{i,j} \frac{F_{k+1}}{\sqrt{\mu}}dv\right),\\
	 \mathfrak{g}_{k}&=-\sum_{i=1}^3\partial_{x_i}\left(2T^{\f32}\int_{\mathbb{R}^3}\mathcal{B}_i\frac{F_{k+1}}{\sqrt{\mu}}dv+\sum_{j=1}^32\fu_{j}T\int_{\mathbb{R}^3}\mathcal{A}_{i,j}\frac{F_{k+1}}{\sqrt{\mu}}dv\right)-2\fu\cdot \mathfrak{f}_{k},
	\end{split}
	\end{align}
where the $\mathcal{A}_{ij}$ and $\mathcal{B}_{i}$ are the Burnett functions defined \eqref{2.4},
and  we have used  the subscript $k$ for the source terms $\mathfrak{f}_{k}$ and $\mathfrak{g}_{k}$ in order to emphasize that the right hand side depends only on $F_i$'s for $0\leq i\leq k$.
\end{lemma}

\begin{remark}\label{rem2.3}
To solve \eqref{1.7-1}, it is equivalent to solve the linear  hyperbolic system \eqref{2.2}. Since it is an initial boundary value problem, we still need to impose a suitable boundary condition to \eqref{2.2}. In fact, to ensure each $F_k+\bar{F}_k+\hat{F}_k$ satisfies the specular reflection boundary conditions and Knudsen layer are solvable, we can not impose boundary data of \eqref{2.2} arbitrarily.
 It is  a very technical process to determine the boundary condition and we will show the details  in  section \ref{sec2.4} below.
\end{remark}

\begin{remark}
The original version of Lemma \ref{lem2.1} in \cite{Guo Jang} is for the Hilbert expansion of Vlasov-Poisson-Boltzmann equations,  and one can obtain Lemma \ref{lem2.1} by dropping the electric field. Lemma \ref{lem2.1} is slightly different from the original version in \cite{Guo Jang} because we also consider the orders of $\v^{2k-1}$, but the proof is very similar, so we omit the details for simplicity of presentation. Noting  $\mathcal{A}_{i,j}$ and $\mathcal{B}_{i}$ are microscopic functions,  the source term $\mathfrak{f}_{k}$ and $\mathfrak{g}_{k}$ depend only on the microscopic part $\mathbf{(I-P)}(\frac{F_{k+1}}{\sqrt{\mu}})$ and hence depend only on $F_i$'s for $0\leq i\leq k$.
\end{remark}

%Since we consider the initial boundary value problem, hence for the linearized system  \eqref{2.2} we impose it with slip boundary condition
%\begin{equation}\label{2.5}
%u_{k+1,3}|_{x_3=0}=0.
%\end{equation}

\smallskip

\subsection{Proof of Theorem \ref{thm2.5}}
%Recall the definition of $\bar{f}_k$ in \eqref{2.8-1}.
%\noindent{\bf Proof.}
For the macroscopic variables $\bar{\rho}_k, \bar{u}_{k}$ and $\bar{\theta}_k$ of $\bar{F}_k$, a direct calculation shows that
\begin{align}\label{2.27}
\begin{split}
&\int_{\R^3} \bar{F}_k dv= \bar{\rho}_k,\quad
\int_{\R^3} (v_i-\fu^0_i) \bar{F}_k  dv=\rho^0 \bar{u}_{k,i},
\quad \int_{\R^3} v_i \bar{F}_k  dv=\rho^0 \bar{u}_{k,i}+\bar{\rho}_k \fu^0_{i},\\
&\int_{\R^3} |v|^2 \bar{F}_k dv =\rho^0 \bar{\theta}_k+3T^0\bar{\rho}_k+2\rho^0\fu^0\cdot\bar{u}_k+\bar{\rho}_k |\fu^0|^2,\\
&\int_{\R^3} |v-\fu^0|^2 \bar{F}_k  dv =\rho^0 \bar{\theta}_k+3T^0 \bar{\rho}_k,\\
&\int_{\R^3} v_i^2 \bar{F}_k  dv=2\rho^0 \fu^0_i \bar{u}_{k,i}+\bar{\rho}_{k} |\fu^0_i|^2+\frac{\rho^0 \bar{\theta}_k+3T^0 \bar{\rho}_k}{3}  + \int_{\R^3} T^0  (\mathbf{I-P}_0) \bar{f}_k\cdot \mathcal{A}_{ii}^0 dv,\\
&\int_{\R^3} v_iv_j \bar{F}_k  dv=\rho^0 \fu^0_i \bar{u}_{k,j}+\rho^0 \fu^0_j \bar{u}_{k,i}+\bar{\rho}_{k} \fu^0_i\fu^0_j+ \int_{\R^3} T^0  (\mathbf{I-P}_0) \bar{f}_k\cdot \mathcal{A}_{ij}^0 dv,\  i\neq j,
\end{split}
\end{align}
\begin{align*}%\label{2.28}
\int_{\R^3} v_i|v|^2 \bar{F}_k  dv&=(5T^0+|\fu^0|^2) \rho^0 \bar{u}_{k,i}+\frac{5}{3}\fu^0_i  (\rho^0 \bar{\theta}_k+3T^0 \bar{\rho}_k)+\fu^0_i (2\rho^0 \fu^0\cdot \bar{u}_k+\bar{\rho}_k |\fu^0|^2)\nonumber\\
&\quad+ \sum_{l=1}^2 \int_{\R^3} 2T^0\fu^0_l \mathcal{A}_{il}^0 \cdot (\mathbf{I-P}_0) \bar{f}_k dv+\int_{\R^3} (T^0)^{\f32} \mathcal{B}_i^0\cdot (\mathbf{I-P}_0) \bar{f}_k dv.
\end{align*}
and
\begin{align}\label{2.28-1}
\begin{split}
&\int_{\R^3} v_3\sqrt{\mu_0} \cdot\partial_y \mathbf{P}_0\bar{f}_k dv=\rho^0 \partial_y\bar{u}_{k,3},\\
&\int_{\R^3} v_iv_3\sqrt{\mu_0} \cdot\partial_y \mathbf{P}_0\bar{f}_k dv=\rho^0 \fu^0_i \partial_y\bar{u}_{k,3}, \  i=1,2,\\
&\int_{\R^3} v_3^2\sqrt{\mu_0}\cdot \partial_y \mathbf{P}_0\bar{f}_k dv= \partial_y(\frac{\rho^0 \bar{\theta}_{k}+3T^0 \bar{\rho}_{k}}{3}),\\
&\int_{\R^3} v_3|v|^2\sqrt{\mu_0}\cdot \partial_y \mathbf{P}_0\bar{f}_k dv=\rho^0 (5T^0+|\fu^0|^2) \partial_y \bar{u}_{k,3},\\
&\int_{\R^3} v_3(|v-\fu^0|^2-5T^0)\sqrt{\mu_0}\cdot \partial_y \mathbf{P}_0\bar{f}_k \, dv=0.
\end{split}
\end{align}

Multiplying $\eqref{1.14}_3$ by $1,v,|v|^2$, integrating over $\R^3$ and using \eqref{2.27}-\eqref{2.28-1}, we obtain
\begin{align}
&\partial_t \bar{\rho}_k+\sum_{j=1}^2 \pa_j (\rho^0 \bar{u}_{k,j}+\bar{\rho}_k \fu^0_j)+\rho^0 \partial_y\bar{u}_{k+1,3}=0,\label{2.28-2}\\
&\partial_t(\rho^0 \bar{u}_{k,i}+\bar{\rho}_k \fu^0_i)
+\sum_{j=1}^2 \partial_j\Big(\rho^0\fu^0_j \bar{u}_{k,i}+\rho^0\fu^0_i \bar{u}_{k,j}+\bar{\rho}_k \fu^0_i \fu^0_j +\delta_{ij} \frac{\rho^0 \bar{\theta}_{k}+3T^0\bar{\rho}_{k}}{3}\Big)+ \rho^0 \fu^0_i \partial_y\bar{u}_{k+1,3}\nonumber\\
&\quad+\delta_{3i} \,\partial_y (\frac{\rho^0 \bar{\theta}_{k+1}+3T^0 \bar{\rho}_{k+1}}{3})+\int_{\R^3} v_iv_3 \sqrt{\mu_0}\cdot (\mathbf{I-P}_0)\partial_y\bar{f}_{k+1} dv=\bar{W}_{k-1,i},\label{2.28-3}
\end{align}
and
\begin{align}
&\partial_t(\rho^0\bar{\theta}_k+3T^0\bar{\rho}_k+2\rho^0\fu^0\cdot\bar{u}_k+\bar{\rho}_k |\fu^0|^2)\nonumber\\
&+\sum_{j=0}^2 \partial_j \{(5T^0+|\fu^0|^2) \rho^0 \bar{u}_{k,i}+\frac{5}{3}\fu^0_i  (\rho^0 \bar{\theta}_k+3T^0 \bar{\rho}_k)+\fu^0_i (2\rho^0 \fu^0\cdot \bar{u}_k+\bar{\rho}_k |\fu^0|^2)\}\nonumber\\
&+\rho^0 (5T^0+|\fu^0|^2) \partial_y \bar{u}_{k+1,3} +\int_{\R^3} v_3 |v|^2 \sqrt{\mu_0}\cdot (\mathbf{I-P}_0)\partial_y\bar{f}_{k+1} dv\nonumber\\
&=-\sum_{j=1}^2 \partial_j\left\{\sum_{l=1}^2 2T^0\fu^0_l\int_{\R^3}  \mathcal{A}_{jl}^0 \cdot (\mathbf{I-P}_0) \bar{f}_k dv+(T^0)^{\f32}\int_{\R^3}  \mathcal{B}_j^0\cdot (\mathbf{I-P}_0) \bar{f}_k dv \right\}.\label{2.28-4}
\end{align}
Then \eqref{2.24} follow directly from \eqref{2.28-2}.

Substituting \eqref{2.28-2} into \eqref{2.28-3} and by tedious calculations, one can obtain
\begin{align}
&\rho^0 \partial_t \bar{u}_{k,i}+\rho^0 (\fu^0_{\sp}\cdot\nabla_\sp) \bar{u}_{k,i}-\rho^0 \partial_3\fu_3^0\, \bar{u}_{k,i}-\bar{\rho}_k \frac{\partial_i p^0}{\rho^0}+\rho^0 \bar{u}_{k,\sp}\cdot \nabla_\sp \fu^0_i\nonumber\\
&\qquad\qquad\qquad\qquad\qquad\qquad
+\partial_y\langle T^0\mathcal{A}_{3i}^0,\  (\mathbf{I-P}_0)\bar{f}_{k+1} \rangle
=\bar{W}_{k-1,i}-\partial_i\bar{p}_k,%(\frac{\rho^0 \bar{\theta_{k}^\v}+3T^0\rho^\v_{k}}{3}),
\ i=1,2,\label{2.28-5}\\[2mm]
&\rho^0 \partial_t \bar{u}_{k,3}+\rho^0 (\fu^0_{\sp} \cdot\nabla_\sp) \bar{u}_{k,3}-\rho^0 \partial_3\fu_3^0\, \bar{u}_{k,3}+\partial_y\bar{p}_{k+1}%\Big(\frac{\rho^0 \bar{\theta_{k+1}^\v}+3T^0\rho^\v_{k+1}}{3}\Big)
+\partial_y \langle T^0  \mathcal{A}_{33}^0, \  (\mathbf{I-P}_0)\bar{f}_{k+1} \rangle
=\bar{W}_{k-1,3}.\label{2.28-6}
\end{align}
Similarly, by using \eqref{2.28-2}-\eqref{2.28-3}, the equation \eqref{2.28-4} for $\bar{\theta}_k$ can be reduced to be
\begin{align}\label{2.28-7}
&\frac53\rho^0 \partial_t \bar{\theta}_{k}+\frac53\rho^0 \fu^0_{\sp} \cdot\nabla_\sp \bar{\theta}_{k} +\Big(\frac{10}{9}\rho^0 \mbox{\rm div} \fu^0-\frac53\rho^0 \partial_3\fu^0_3\Big) \bar{\theta}_k\nonumber\\
&+(3\rho^0\nabla_\sp T^0-2T^0\nabla_\sp\rho^0) \bar{u}_{k,\sp}+\partial_y\langle2(T^0)^{\frac32}  \mathcal{B}_3^0,\  (\mathbf{I-P}_0)\bar{f}_{k+1} \rangle \nonumber\\
&=\bar{H}_{k-1}+\{2\partial_t+2\fu^0_\sp\cdot \nabla_\sp+\frac{10}{3}\mbox{\rm div}\fu^0\}\bar{p}_k.
%\Big(\frac{\rho^0 \bar{\theta_{k}^\v}+3T^0\rho^\v_{k}}{3}\Big)
\end{align}

We still need to deal with the microscopic parts in \eqref{2.28-5}-\eqref{2.28-7}. In fact, by using \eqref{1.14}, it is direct to obtain \eqref{2.21}.  Then by using \eqref{2.21}, we have
\begin{align}\label{2.28-8}
&\Big\langle(\mathbf{I-P_0})\bar{f}_{k+1},\mathcal{A}_{3i}^0\Big\rangle
=\Big\langle -\mathbf{L}_0^{-1}\{(\mathbf{I-P_0})(v_3\partial_y\mathbf{P}_0 \bar{f}_k)\},\ \mathcal{A}_{3i}^0\Big\rangle\nonumber\\
&\qquad+\left\langle\mathbf{L}_0^{-1}\Big\{\frac{y}{\sqrt{\mu_0}}[Q(\partial_3\mu_0,\sqrt{\mu_0}\mathbf{P_0}\bar{f}_{k})+Q(\sqrt{\mu_0}\mathbf{P_0}\bar{f}_{k}, \partial_3\mu_0)]\right.\nonumber\\
&\quad\quad \left.+ \frac{1}{\sqrt{\mu_0}} [Q(\sqrt{\mu_0} f_1^0,\sqrt{\mu_0}\mathbf{P_0}\bar{f}_k)+Q(\sqrt{\mu_0}\mathbf{P_0}\bar{f}_k, \sqrt{\mu_0} f_1^0)]\right.\nonumber\\
&\quad\quad+\left.\frac1{\sqrt{\mu_0}} [Q(\sqrt{\mu_0}\bar{f}_1,\sqrt{\mu_0}\mathbf{P_0}\bar{f}_k)+Q(\sqrt{\mu_0}\mathbf{P_0}\bar{f}_k,\sqrt{\mu_0}\bar{f}_1)]\Big\},\  \mathcal{A}_{3i}^0\right\rangle+\Big\langle\bar{J}_{k-1},\  \mathcal{A}_{3i}^0\Big\rangle.
\end{align}
A direct calculation shows that
\begin{align}\label{2.28-9}
&\Big\langle -\mathbf{L}_0^{-1}\{(\mathbf{I-P_0})(v_3\partial_y\mathbf{P}_0 \bar{f}_k)\},\ \mathcal{A}_{3i}^0\Big\rangle=\Big\langle -v_3\partial_y\mathbf{P}_0 \bar{f}_k,\ \mathbf{L}_0^{-1}\mathcal{A}_{3i}^0\Big\rangle\nonumber\\
&=-\Big\langle\Big\{ \frac{\partial_y\bar\rho_k}{\rho^0}+\partial_y\bar{u}_{k}\cdot \frac{v-\fu^0}{T^0}+\frac{\partial_y\bar{\theta}_k}{6T^0}(\frac{|v-\fu^0|^2}{T^0}-3)\Big\}v_3\sqrt{\mu_0},\  \mathbf{L}_0^{-1}\mathcal{A}_{3i}^0\Big\rangle\nonumber\\
&=-\partial_y\bar{u}_{k,i} \Big\langle \mathcal{A}_{3i}^0 , \ \mathbf{L}_0^{-1}\mathcal{A}_{3i}^0\Big\rangle
=
\begin{cases}
\displaystyle -\frac{\mu(T^0)}{T^0} \partial_y\bar{u}_{k,i},\  i=1,2,\\[2mm]
\displaystyle -\frac43\frac{\mu(T^0)}{T^0} \partial_y\bar{u}_{k,3},\ i=3,
\end{cases}
\end{align}
where we have used \eqref{2.15-2} in the last equality.

From \cite[p.649]{Guo2006}, we know that
\begin{align}\label{0.1}
&\mathbf{L}_0^{-1}\Big\{\frac{1}{\sqrt{\mu_0}} Q(\sqrt{\mu_0} \mathbf{P}_0g,\sqrt{\mu_0} \mathbf{P}_0\bar{f}_k)+\frac{1}{\sqrt{\mu_0}} Q(\sqrt{\mu_0}\mathbf{P}_0\bar{f}_k, \sqrt{\mu_0} \mathbf{P}_0g)\Big\}\nonumber\\
&=(\mathbf{I-P}_0)\left\{\frac{\mathbf{P}_0g\cdot\mathbf{P}_0\bar{f}_k}{\sqrt{\mu_0}}\right\}.
\end{align}
We assume
\begin{equation}\nonumber
 \mathbf{P}_0g=\left\{ \frac{a}{\rho^0}+b\cdot \frac{v-\fu^0}{T^0}+\frac{c}{6T^0}(\frac{|v-\fu^0|^2}{T^0}-3)\right\}\sqrt{\mu_0},
\end{equation}
then a direct calculation show that
\begin{align}\label{0.2}
(\mathbf{I-P}_0)\left\{\frac{\mathbf{P}_0g\cdot\mathbf{P}_0\bar{f}_k}{\sqrt{\mu_0}}\right\}
&=\sum_{l,j=1}^3 \frac{b_l \, \bar{u}_{k,j}}{T^0} \mathcal{A}_{lj}^0
+\frac1{\sqrt{T^0}}\Big(\frac{\bar{\theta}_k}{3T^0} b+\frac{c}{3T^0}\bar{u}_{k}\Big)\cdot \mathcal{B}^0\nonumber\\
&\quad+\frac{c\cdot \bar{\theta}_k}{36(T^0)^2}(\mathbf{I-P}_0)\left\{ \Big(\frac{|v-\fu^0|^2}{T^0}-5\Big)^2\sqrt{\mu_0}\right\}.
\end{align}
Noting
\begin{equation}\label{0.8}
\frac{\partial\mu}{\sqrt{\mu}}=\left\{ \frac{\partial\rho}{\rho}+\partial\fu\cdot \frac{v-\fu}{T}+\frac{3\partial T}{6T}(\frac{|v-\fu|^2}{T}-3)\right\}\sqrt{\mu},
\end{equation}
which, together with \eqref{0.1} and \eqref{0.2},  yields that the second term on RHS of \eqref{2.28-8} is expressed as
\begin{align}\label{0.3}
&\sum_{l,j=1}^3\frac{\bar{u}_{k,j}}{T^0} [y \partial_3 \fu^0_{l}+u^0_{1,l}+\bar{u}_{1,l} ] \Big\langle  \mathcal{A}_{lj}^0,\  \mathcal{A}_{3i}^0\Big\rangle\nonumber\\
&=
\begin{cases}
\displaystyle \frac{\rho^0}{T^0} \{[\partial_3\fu^0_3\cdot y+u^0_{1,3}+\bar{u}_{1,3}]\bar{u}_{k,i}+[\partial_3\fu^0_i\cdot y+u^0_{1,i}+\bar{u}_{1,i}]\bar{u}_{k,3}\},\  i=1,2,\\[2mm]
\displaystyle \frac{4}{3} \frac{\rho^0}{T^0} [\partial_3\fu^0_3\cdot y+u^0_{1,3}+\bar{u}_{1,3}]\bar{u}_{k,3}, \  i=3,
\end{cases}
\end{align}
where we have used
\begin{align*}%\label{0.4}
\Big\langle \mathcal{A}_{ij},\ \mathcal{A}_{ij} \Big\rangle=\rho,\  \mbox{for}\ i\neq j,\quad \mbox{and} \  \Big\langle \mathcal{A}_{ii},\ \mathcal{A}_{ii} \Big\rangle=\frac43\rho.
\end{align*}
Combining \eqref{2.28-8}, \eqref{2.28-9} and \eqref{0.3}, then we obtain
\begin{align}\label{0.5}
&\partial_y\langle T^0\mathcal{A}_{3i}^0, \   (\mathbf{I-P}_0)\bar{f}_{k+1} \rangle\nonumber\\
&=
\begin{cases}
\displaystyle - \mu(T^0)\partial_{yy}\bar{u}_{k,i}+\rho^0\partial_y\{[\partial_3\fu^0_3\cdot y+u^0_{1,3}+\bar{u}_{1,3}]\bar{u}_{k,i}\}\\[1mm]
\displaystyle\qquad\qquad\qquad\quad+\rho^0\partial_y\{[\partial_3\fu^0_i\cdot y+u^0_{1,i}+\bar{u}_{1,i}]\bar{u}_{k,3}\} +T^0\partial_y \langle\bar{J}_{k-1},\  \mathcal{A}_{3i}^0\rangle,\  i=1,2,\\[2mm]
\displaystyle -\frac43\mu(T^0) \partial_{yy}\bar{u}_{k,3}+\frac43\rho^0\partial_y\{[\partial_3\fu^0_3\cdot y+u^0_{1,3}+\bar{u}_{1,3}]\bar{u}_{k,3}\}+T^0\partial_y \langle\bar{J}_{k-1},\  \mathcal{A}_{33}^0\rangle,\ i=3.
\end{cases}
\end{align}
Substituting \eqref{0.5} into \eqref{2.28-5} and \eqref{2.28-6}, then by  tedious calculations, one can get  \eqref{2.19-2} and \eqref{2.25}.

By similar arguments as in \eqref{2.28-8}-\eqref{0.5}, we can obtain
\begin{align}\label{0.6}
 \partial_y\langle2 (T^0)^{\frac32} \mathcal{B}_{3}^0, \ (\mathbf{I-P}_0)\bar{f}_{k+1} \rangle
&=-\kappa(T^0) \partial_{yy} \bar{\theta}_k +\frac53 \rho^0 \partial_y \{(\partial_3\fu^0_3\cdot y+u^0_{1,3}+\bar{u}_{1,3}) \bar{\theta}_k\}\nonumber\\
&\hspace{-6mm}+\frac53 \rho^0 \partial_y \{(3\partial_3T^0_3\cdot y+\theta^0_{1}+\bar{\theta}_{1}) \bar{u}_{k,3}\}+2 (T^0)^{\frac32} \partial_y\langle\bar{J}_{k-1},\ \mathcal{B}_{3}^0\rangle,
\end{align}
which, together with \eqref{2.28-7} and \eqref{2.14}, yields \eqref{2.20}. Therefore the proof of Theorem \ref{thm2.5} is completed.  $\hfill\Box$

\subsection{Reformulation of Knudsen Boundary Layer}
Define $\displaystyle \hat{f}_k:=\frac{\hat{F}_k}{\sqrt{\mu_0}}$, then we can rewrite \eqref{1.19-1} as
\begin{align}\label{2.33}
 v_3\frac{\partial\hat{f}_{k}}{\partial \eta }+\mathbf{L}_0\hat{f}_{k}=\hat{S}_{k}, \quad k\geq1,
 \end{align}
 where $\hat{S}_k:=\hat{S}_{k,1}+\hat{S}_{k,2} (k\geq1)$ with
\begin{align}
\hat{S}_{k,1}&=-\mathbf{P_0}\left\{\frac{\{\partial_t+v_\sp\cdot\nabla_\sp\}\hat{F}_{k-2}}{\sqrt{\mu_0}}\right\},\label{2.35}\\
\hat{S}_{k,2}&=\sum_{\substack{i+j=k\\ i,j\geq1}} \frac{1}{\sqrt{\mu_0}} \big[Q(F_i^0+\bar{F}_i^0,\sqrt{\mu_0}\hat{f}_j)+Q(\sqrt{\mu_0}\hat{f}_j,F_i^0+\bar{F}_i^0)+Q(\sqrt{\mu_0}\hat{f}_i, \sqrt{\mu_0}\hat{f}_j) \big] \nonumber\\
&\qquad+\sum_{\substack{j+2l=k\\  1\leq l\leq \fb, j\geq1}} \frac{\eta^l}{l!} \frac{1}{\sqrt{\mu_0}}\big[Q(\partial_3^l\mu_0, \sqrt{\mu_0}\hat{f}_{j})+Q( \sqrt{\mu_0}\hat{f}_{j},\partial_3^l\mu_0) \big]\nonumber \\
&\qquad+\sum_{\substack{i+j+2l=k\\ i,j\geq1, 1\leq l\leq \fb}}\frac{\eta^l}{l!} \frac{1}{\sqrt{\mu_0}}\big[Q(\partial_3^l F_i^0, \sqrt{\mu_0}\hat{f}_{j})+Q( \sqrt{\mu_0}\hat{f}_{j},\partial_3^l F_i^0) \big] \label{2.34} \\
&\qquad+\sum_{\substack{i+j+l=k\\ i,j\geq1, 1\leq l\leq \fb}}\frac{\eta^l}{l!} \frac{1}{\sqrt{\mu_0}} \big[Q(\partial_y^l \bar{F}_i^0, \sqrt{\mu_0}\hat{f}_{j})+Q( \sqrt{\mu_0}\hat{f}_{j},\partial_y^l \bar{F}_i^0) \big]\nonumber\\
&\qquad-(\mathbf{I-P_0})\left\{\frac{\{\partial_t+v_\sp\cdot\nabla_\sp\}\hat{F}_{k-2}}{\sqrt{\mu_0}}\right\}.\nonumber
\end{align}
Here we have used the notation $\hat{F}_{-1}=\hat{F}_{0}=0$ for simplicity of presentation. It is easy to notice that $\hat{S}_{k,1}\in \mathcal{N}_0$, $\hat{S}_{k,2}\in \mathcal{N}_0^{\perp}$ and
\begin{equation}\label{2.36}
\hat{S}_1=\hat{S}_{1,1}=\hat{S}_{1,2}=\hat{S}_{2,1}=0.
\end{equation}

For later use, we introduce a result  on the existence of solution to the Knudsen boundary layer problem with a perturbed specular reflection boundary conditions.  Consider the following half-space linear  problem
\begin{align}\label{2.29}
\begin{cases}
&\dis v_3 \partial_{\eta} f+\mathbf{L}_0f=S(t,x_\sp,\eta,v),\\[2mm]
&\dis f(t,x_\sp, 0,v_\sp,v_3)|_{v_3>0}=f(t,x_\sp,0,v_\sp,-v_3)+f_b(t,x_\sp,v_\sp,-v_3),\\[2mm]
&\dis \lim_{\eta\rightarrow\infty} f(t,x_\sp,\eta,v)=0,
\end{cases}
\end{align}
where $\eta\in\R_+$ and $(t,x_\sp)\in [0,\tau]\times \R^2$. In fact, we think $(t,x_\sp)\in [0,\tau]\times \R^2$ to be parameters in \eqref{2.29}. The function $f_b(t,x_\sp,v)$ is defined only for $v_3< 0$, and we always assume that it is extended to be $0$ for $v_3>0$.
\begin{lemma}[\cite{Jiang-Wang}]\label{lem2.6}
Let $0\leq \fa<\frac12$ and $\k\geq 3$. For each $(t,x_\sp)\in[0,\tau]\times \R^2$, we assume  that
\begin{align}\label{2.30}
S\in \mathcal{N}_0^{\perp}\quad \mbox{and}\quad \|w_\k \mu_0^{-\fa} f_b(t,x_\sp,0,\cdot)\|_{L^\infty_v}+
\|\nu^{-1} w_\k \mu_0^{-\fa} e^{\zeta_0\eta}S(t,x_\sp,\cdot,\cdot)\|_{L^\infty_{\eta,v}}<\infty,
\end{align}
for some positive constant $\zeta_0>0$, and
\begin{align}\label{2.31}
\begin{cases}
\dis \int_{\R^3}v_3\,  f_b(t,x_\sp,v)\,  \sqrt{\mu_0}\, dv&\equiv0,\\[2mm]
\dis \int_{\R^3}(v_1-\fu_1^0) v_3\,  f_b(t,x_\sp,v) \, \sqrt{\mu_0}\, dv&\equiv0,\\[2mm]
\dis \int_{\R^3}(v_2-\fu_2^0)v_3\,  f_b(t,x_\sp,v) \, \sqrt{\mu_0}\, dv&\equiv0,\\[2mm]
\dis \int_{\R^3} |v-\fu^0|^2v_3 \, f_b(t,x_\sp,v) \, \sqrt{\mu_0}\, dv&\equiv0.
\end{cases}
%\qquad\quad\forall (t,x_\sp)\in \R_+\times \R^2.
\end{align}
Then the boundary value problem \eqref{2.29}-\eqref{2.31} has a unique solution $f$  satisfying
\begin{align*}%\label{2.52-2}
& \|w_{\k} \mu_0^{-\fa}e^{\zeta\eta}f(t,x_\sp,\cdot,\cdot)\|_{L^\infty_{\eta,v}}
+\|w_{\k} \mu_0^{-\fa} f(t,x_\sp,0,\cdot)\|_{L^\infty_v}\nonumber\\
&\leq \frac{C}{\zeta_0-\zeta} \, \Big(\|w_\k \mu_0^{-\fa}f_b(t,x_\sp,0,\cdot)\|_{L^\infty_v}+
\|\nu^{-1}w_\k \mu_0^{-\fa} e^{\zeta_0 \eta}S(t,x_\sp,\cdot,\cdot)\|_{L^\infty_{\eta,v}}\Big),
\end{align*}
for all  $ (t,x_\sp)\in [0,\tau]\times \R^2$, where $C>0$ is a positive constant independent of $(t,x_\sp)$, and  $\zeta$ is any positive constant such that $0<\zeta< \zeta_0$. Moreover, if $S$ is continuous in $(t,x_\sp,\eta,v)\in [0,\tau]\times \R^2\times\R_+\times\R^3$ and $f_b$ is continuous in $(t,x_\sp,v_\sp,-v_3)\in [0,\tau]\times \R^2\times\R^2\times\R_+$, then the solution $f$ is continuous away from the grazing set $[0,\tau]\times\gamma_0$.%$\{(x_\sp,0, v)\, | \, x_\sp\in\R^2, \, v_\sp\in\R^2, v_3\neq0 \}$.
\end{lemma}

\begin{remark}
As indicated in \cite{GKTT}, in general, it is hard to obtain  the normal derivatives estimates for the boundary value problem \eqref{2.29}. Fortunately, it is easy to obtain the tangential and time derivatives estimates for the solution of \eqref{2.29}, i.e.,
\begin{align}\label{2.53-1}
& \sum_{i+j\leq r}\|w_{\k} \mu_0^{-\fa} e^{\zeta\eta}\partial_t^{i}\nabla_\sp^{j}f(t,x_\sp,\cdot,\cdot)\|_{L^\infty_{\eta,v}}
+\|w_{\k} \mu_0^{-\fa} \partial_t^{i}\nabla_\sp^{j}f(t,x_\sp,0,\cdot)\|_{L^\infty_v}\nonumber\\
%\|w_{\beta}e^{\zeta\eta}\partial_t^{i}\nabla^{j}f\|_{L^\infty_{\eta,v}}\nonumber\\
&\leq \frac{C}{\zeta_0-\zeta}  \sum_{i+j\leq r}\Big\{ \|w_\k \mu_0^{-\fa} \partial_t^{i}\nabla_\sp^{j}f_b(t,x_\sp,\cdot)\|_{L^\infty_v}+
\|\nu^{-1}w_\k \mu_0^{-\fa} e^{\zeta_0 \eta} \partial_t^{i}\nabla_\sp^{j}S\|_{L^\infty_{\eta,v}}\Big\},
\end{align}
 provided the right hand side of \eqref{2.53-1} is bounded. And such estimate \eqref{2.53-1} is enough for us to establish the Hilbert expansion. Moreover, taking $L^\infty_{x_\sp}\cap L^2_{x_\sp}$ over \eqref{2.53-1}, one obtains
 \begin{align}\label{2.53-2}
 & \sum_{i+j\leq r} \sup_{t\in[0,\tau]} \Big\{\|w_{\k} \mu_0^{-\fa}e^{\zeta\eta}\partial_t^{i}\nabla_\sp^{j}f(t)\|_{L^\infty_{x_\sp,\eta,v}\cap L^2_{x_\sp}L^\infty_{\eta,v}}
 +\|w_{\k} \mu_0^{-\fa} \partial_t^{i}\nabla_\sp^{j}f(t,\cdot,0,\cdot)\|_{L^\infty_{x_\sp,v}\cap L^2_{x_\sp}L^\infty_{v}}\Big\}\nonumber\\
 &\leq \frac{C}{\zeta_0-\zeta}  \sup_{t\in[0,\tau]}\left\{ \sum_{i+j\leq r} \Big\{ \|w_\k \mu_0^{-\fa} \partial_t^{i}\nabla_\sp^{j}f_b(t)\|_{L^\infty_{x_\sp,v}\cap L^2_{x_\sp}L^\infty_{v}}\right.\nonumber\\
 &\qquad\qquad\qquad\qquad\qquad\left.+
\sum_{i+j\leq r} \|\nu^{-1}w_\k \mu_0^{-\fa} e^{\zeta_0 \eta} \partial_t^{i}\nabla_\sp^{j}S(t)\|_{L^\infty_{x_\sp,\eta,v}\cap L^2_{x_\sp}L^\infty_{\eta,v}}\right\}.
 \end{align}
\end{remark}

\begin{remark}
Golse, Perthame and Sulem \cite{GPS-1988} have proved an existence result for \eqref{2.29} in the space  $\dis \int_{\R_+\times\R^3}(1+|v|) e^{2\zeta \eta} f^2 dvd\eta+\int_{\R^3}\|e^{\zeta \eta} f \|^2_{L^\infty_{\eta}} dv$. In the present paper, since  the continuity and  weighted  $L^\infty_{\eta,v}$ estimate are need, so the second and third authors of present paper  proved  the well-posedness of \eqref{2.29} in  the new functional space, see \cite{Jiang-Wang}.
\end{remark}

Since  the source term $S\in \mathcal{N}_0^{\perp}$ in Lemma \ref{lem2.6} is demanded, but  $\hat{S}_{k}\notin \mathcal{N}_0^{\perp} (k\geq3)$ in general, i.e., $\hat{S}_{k,1}\neq0$. Hence, to solve \eqref{2.33}, we need to cancel the term $\hat{S}_{k,1}$.
We assume that
\begin{equation}\label{2.37}
\hat{S}_{k,1}=\left\{ \hat{a}_k+\hat{b}_{k}\cdot (v-\fu^0)+\hat{c}_k |v-\fu^0|^2\right\}\sqrt{\mu_0}.
\end{equation}
where $(\hat{a}_k,\hat{b}_k,\hat{c}_k)=(\hat{a}_k,\hat{b}_k,\hat{c}_k)(t,x_\sp,\eta)$.
By similar arguments as in \cite{BCN-1986}, we have the following lemma. The details of proof are omitted for simplicity of presentation.
\begin{lemma}\label{lem2.7}
For $(\hat{a}_k,\hat{b}_k,\hat{c}_k)$ defined in \eqref{2.37}, we assume that
\begin{align*}%\label{2.38-1}
\lim_{\eta\rightarrow\infty} e^{\zeta\eta}|(\hat{a}_k,\hat{b}_k,\hat{c}_k)(t,x_\sp,\eta)|=0,
\end{align*}
for some positive constant $\zeta>0$. Then there exists a function
\begin{align}\label{2.38}
\hat{f}_{k,1}=\Big\{\hat{A}_k v_3+ \hat{B}_{k,1} v_3(v_1-\fu_1^0)+\hat{B}_{k,2} v_3(v_2-\fu_2^0)+\hat{B}_{k,3}+ \hat{C}_{k} v_3 |v-\fu^0|^2 \Big\}\sqrt{\mu_0},
\end{align}
such that
\begin{equation*}%\label{2.39}
v_3\partial_{\eta} \hat{f}_{k,1}-\hat{S}_{k,1}\in \mathcal{N}_0^{\perp},
\end{equation*}
where
\begin{align}\label{2.40}
\begin{split}
\hat{A}_k(t,x_\sp,\eta)&=-\int_{\eta}^{\infty}\Big(\frac2{T^0} \hat{a}_k+3\hat{c}_k\Big)(t,x_\sp,s)ds,\\
\hat{B}_{k,i}(t,x_\sp,\eta)&=-\int_{\eta}^{\infty}\frac1{T^0} \hat{b}_{k,i}(t,x_\sp,s)ds,\quad i=1,2,\\
\hat{B}_{k,3}(t,x_\sp,\eta)&=-\int_{\eta}^{\infty}\hat{b}_{k,3}(t,x_\sp,s)ds,\\
\hat{C}_k(t,x_\sp,\eta)&=\frac1{5(T^0)^2}\int_{\eta}^{\infty}\hat{a}_{k}(t,x_\sp,s)ds.
\end{split}
\end{align}
Moreover it holds that
\begin{align*}%\label{2.42-1}
|v_3\partial_{\eta} \hat{f}_{k,1}-\hat{S}_{k,1}|\leq C|(\hat{a}_k,\hat{b}_k,\hat{c}_k)(t,x_\sp,\eta)| (1+|v|)^4\sqrt{\mu_0},
\end{align*}
and
\begin{equation*}%\label{2.43}
|\hat{f}_{k,1}(t,x_\sp,\eta,v)| \leq C(1+|v|)^3\sqrt{\mu_0} \int_{\eta}^{\infty} |(\hat{a}_k,\hat{b}_k,\hat{c}_k)| \rightarrow 0\, \, \mbox{as}\, \,  \eta\rightarrow\infty.
\end{equation*}
\end{lemma}

\begin{remark}
It is very important to note that $\hat{S}_{k,1}$ depends only on $\hat{f}_{k-2}$, which is  already known function when we consider the existence of $\hat{f}_{k}$. Thus  $\hat{f}_{k,1}$  (or equivalent $(\hat{A}_k,\hat{B}_k,\hat{C}_k)$) is determined by $\hat{f}_{k-2}$. On the other hand, since $\hat{S}_{1,1}=\hat{S}_{2,1}=0$, one has that
\begin{align}\label{2.43-2}
\hat{f}_{1,1}=\hat{f}_{2,1}\equiv0,
\end{align}
which yields that $(\hat{A}_1,\hat{B}_1,\hat{C}_1)\equiv(0,0,0)$ and $(\hat{A}_2,\hat{B}_2,\hat{C}_2)\equiv(0,0,0)$.
\end{remark}

\

Now we consider $\hat{f}_{k,2}$ satisfying
\begin{align}\label{2.42}
 v_3\frac{\partial\hat{f}_{k,2}}{\partial \eta }+\mathbf{L}_0\hat{f}_{k,2}=\hat{S}_{k,2}-\mathbf{L}_0\hat{f}_{k,1}-\left(v_3\partial_{\eta} \hat{f}_{k,1}-\hat{S}_{k,1}\right)\in \mathcal{N}_0^{\perp},
\end{align}
which yields that we can use Lemma \ref{lem2.6} to solve $\hat{f}_{k,2}$. Then it is easy to check that
\begin{equation}\label{2.43-1}
\hat{f}_{k}:=\hat{f}_{k,1}+\hat{f}_{k,2},
\end{equation}
is a solution of \eqref{2.33}.

%\section{Existence of solution to expansions}\label{sec4}
%In this section we will construct solution the interior expansion, viscous boundary layer and Knudsen layer. To solve these expansions, we need to determine the corresponding boundary conditions.

\subsection{Boundary Conditions}\label{sec2.4}
To construct the solutions for  interior expansion, viscous and Knudsen boundary layers,  the remain problem is to give  suitable boundary conditions so that they are all well-posed. As mentioned in Remark \ref{rem2.3}, we hope that each $F_k+\bar{F}_k+\hat{F}_k$ satisfies the specular reflection boundary conditions, i.e.,
\begin{align}
(f_k+\bar{f}_k+\hat{f}_{k})(t,x_\sp,0,v_\sp,v_3)=(f_k+\bar{f}_k+\hat{f}_{k})(t,x_\sp,0,v_\sp,-v_3),\nonumber
\end{align}
which, together with \eqref{2.43-1}, yields
\begin{align}\label{2.45}
\hat{f}_{k,2}(t,x_\sp,0,v_\sp,v_3)|_{v_3>0}&=\hat{f}_{k,2}(t,x_\sp,0,v_\sp,-v_3)+[f_k+\bar{f}_k+\hat{f}_{k,1}](t,x_\sp,0,v_\sp,-v_3)\nonumber\\
&\qquad-[f_k+\bar{f}_k+\hat{f}_{k,1}](t,x_\sp,0,v_\sp,v_3).
\end{align}
For simplicity of presentation, we denote
\begin{align}\label{2.46}
\hat{g}_k(t,x_\sp,v_\sp,v_3)=
\begin{cases}
0,\qquad v_3>0,\\[3mm]
[f_k+\bar{f}_k+\hat{f}_{k,1}](t,x_\sp,0,v_\sp,v_3)\\
\qquad-[f_k+\bar{f}_k+\hat{f}_{k,1}](t,x_\sp,0,v_\sp,-v_3),\quad  v_3<0.
\end{cases}
\end{align}
On the other hand, we impose the far field boundary condition
\begin{equation}\label{2.48}
\lim_{\eta\rightarrow\infty} \hat{f}_{k,2}(t,x_\sp,\eta,v)=0.
\end{equation}
Noting from Lemma \ref{lem2.6}, to solve \eqref{2.42}, \eqref{2.45}-\eqref{2.46} and \eqref{2.48}, we need $\hat{g}_k$ to satisfy \eqref{2.31}, i.e.,
\begin{align}
\int_{\R^3}v_3 \hat{g}_k\sqrt{\mu_0}dv
&=\int_{\R^3}(v_1-\fu_1^0)v_3 \hat{g}_k \sqrt{\mu_0}dv
=\int_{\R^3}(v_2-\fu_2^0)v_3 \hat{g}_k \sqrt{\mu_0}dv\nonumber\\
&=\int_{\R^3}v_3 |v-\fu^0|^2\hat{g}_k \sqrt{\mu_0}dv=0,\nonumber
\end{align}
which is equivalent to
\begin{align}\label{2.49}
\begin{split}
&\int_{\R^3} v_3 \sqrt{\mu_0} \mathbf{P}_0[f_k+\bar{f}_k](t,x_\sp,0,v) dv=-\rho^0T^0 (\hat{A}_k+5T^0 \hat{C}_k)(t,x_\sp,0), \\
&\int_{\R^3} v_3 (v_i-\fu_i^0) \sqrt{\mu_0} (\mathbf{I-P_0})[f_k+\bar{f}_k](t,x_\sp,0,v) dv=-\rho^0(T^0)^2 \hat{B}_{k,i}(t,x_\sp,0),\quad i=1,2,\\
&\int_{\R^3} v_3 (|v-\fu^0|^2-5T^0) \sqrt{\mu_0} (\mathbf{I-P_0})[f_k+\bar{f}_k](t,x_\sp,0,v) dv=-10\rho^0 (T^0)^3 \hat{C}_k(t,x_\sp,0),
\end{split}
\end{align}
where we have used \eqref{2.38}. \vspace{1mm}

For the case $k=1$, noting from \eqref{2.12}, \eqref{2.14} and \eqref{2.43-2}, it is direct to know that $\eqref{2.49}_1$ is equivalent to
\begin{align}\label{2.51}
u_{1,3}(t,x_\sp,0)=0.% \quad \forall t\geq0, x_\sp\in\R^2,
\end{align}
Noting $ f_1(t,x_\sp,0,v) , \bar{f}_1(t,x_\sp,0,v) \in \mathcal{N}_0$, then it is direct to know $\eqref{2.49}_{2,3}$ holds naturally for $k=1$.

 \

Now we consider \eqref{2.49} for the case $k\geq 2$. From  $\eqref{2.49}_1$, a direct calculation shows  that
\begin{align}\label{2.50}
 u_{k,3}(t,x_\sp,0)&=-\bar{u}_{k,3}(t,x_\sp,0)-T^0 (\hat{A}_k+5T^0 \hat{C}_k)(t,x_\sp,0) \nonumber\\
 &=-\int_0^\infty\frac{1}{\rho^0}\Big\{\partial_t\bar{\rho}_{k-1}+\mbox{\rm div}_\sp(\rho^0\bar{u}_{k-1,\sp}+\bar{\rho}_{k-1} \fu^0_\sp)\Big\}(t,x_\sp,y)dy\nonumber\\
 &\quad-T^0 (\hat{A}_k+5T^0 \hat{C}_k)(t,x_\sp,0),
\end{align}
where we have used \eqref{2.24}. It is direct to know that the right hand side terms of \eqref{2.50} can be determined by $f_i, \bar{f}_i \, ( i\leq k-1)$ and $\hat{f}_j \, ( j\leq k-2) $.

Now we consider the rest terms of \eqref{2.49}. By similar arguments as in  \eqref{2.28-8}-\eqref{0.6},  and utilizing \eqref{2.14} and \eqref{2.51}, one can obtain
\begin{align}
&\int_{\R^3} v_3 (v_1-\fu_1^0) \sqrt{\mu_0} (\mathbf{I-P_0})\bar{f}_k(t,x_\sp,0,v) dv\nonumber\\
&=-\Big[\mu(T^0)\partial_{y} \bar{u}_{k-1,i}-\rho^0 (u_{1,i}+\bar{u}_{1,i}) \bar{u}_{k-1,3}-\langle T^0\mathcal{A}_{3i}^0,\    \bar{J}_{k-2}\rangle\Big](t,x_\sp,0),\, \, i=1,2,
\label{2.53}\\[2mm]
&\int_{\R^3} v_3 (|v-\fu^0|^2-5T^0) \sqrt{\mu_0} (\mathbf{I-P_0})\bar{f}_k(t,x_\sp,0,v) dv\nonumber\\
&=-\Big[\kappa(T^0) \partial_y\bar{\theta}_{k-1}-\frac53 \rho^0 (\theta_1+\bar{\theta}_1) \bar{u}_{k-1,3}-\langle 2(T^0)^{\frac32} \mathcal{B}_{3}^0,\,   \bar{J}_{k-2}\rangle\Big](t,x_\sp,0).\label{2.54}
\end{align}
Using \eqref{2.53} and \eqref{2.54}, we can rewrite $\eqref{2.49}_{2,3}$ as
\begin{align}
\partial_{y} \bar{u}_{k-1,i}(t,x_\sp,0)
&=\frac{1}{\mu(T^0)}\Big\{[\rho^0 (u_{1,i}+\bar{u}_{1,i}) \bar{u}_{k-1,3}]+ \langle T^0  \mathcal{A}_{3i}^0,\   \bar{J}_{k-2} \rangle\nonumber\\
&\qquad+\langle T^0 \mathcal{A}_{3i}^0, \  (\mathbf{I-P_0})f_k\rangle+\rho^0(T^0)^2 \hat{B}_{k,i}\Big\}(t,x_\sp,0), \  i=1,2, \label{2.55}\\
\partial_y\bar{\theta}_{k-1}(t,x_\sp,0)
&=\frac{1}{\kappa(T^0)}\Big\{\frac53 \rho^0 [(\theta_1+\bar{\theta}_1) \bar{u}_{k-1,3}]+\langle 2(T^0)^{\frac32}\mathcal{B}_{3}^0,\ \bar{J}_{k-2}\rangle\nonumber\\
&\qquad + \langle 2(T^0)^{\frac32}\mathcal{B}_3^0,\   (\mathbf{I-P_0})f_k\rangle+10\rho^0 (T^0)^3 \hat{C}_k\Big\}(t,x_\sp,0).\label{2.56}
\end{align}
\begin{remark}
It is direct to check that the terms on RHS of \eqref{2.55}-\eqref{2.56} depends on $f_i \,  (i\leq k-1)$, $\bar{f}_j\, (j\leq k-2)$ and $\hat{f}_l\, (l\leq k-2)$. Once we solved $(\bar{u}_{k-1}, \bar{\theta}_{k-1})$ with the boundary condition \eqref{2.55}-\eqref{2.56}, then  $\hat{f}_k$ will be solvable by using Lemma \ref{lem2.6}.
\end{remark}

\section{Existence of Solution for a Linear Hyperbolic System}\label{section3}
To study existence of interior expansion, we  first need to consider the following linear problem for $(\tilde{\rho}, \tilde{u}, \tilde{\theta})(t,x)$:
\begin{equation}\label{5.1}
\begin{cases}
\dis \pa_t\tilde{\rho}+\mbox{\rm div}_x(\rho \tilde{u}+\tilde{\rho}\fu)=0,\\[1.5mm]
\dis \rho \{\partial_t \tilde{u}+\tilde{u}\cdot \nabla_x \fu+\fu\cdot\nabla_x \tilde{u}\}-\f{\nabla_xp}{\rho}\tilde{\rho}+\nabla_x(\f{\rho\tilde{\theta}+3T\tilde{\rho}}{3})=\mathfrak{f},\\[3mm]
\dis \rho \{\pa_t\tilde{\theta}+\f23(\tilde{\theta}\mbox{\rm div}_x\fu+3T\mbox{\rm div}_x\tilde{u})+\fu\cdot\nabla_x\tilde{\theta}+3\tilde{u}\cdot\nabla_xT\}=\mathfrak{g},
\end{cases}
\end{equation}
with $(t,x)\in (0,\tau)\times \R_+^3.$ We impose \eqref{5.1} with a given boundary condition
\begin{equation}\label{5.2}
\tilde{u}_3(t,x_\sp,0)=d(t,x_\sp), \quad \forall (t,x)\in (0,\tau)\times \R^2,
\end{equation}
and the initial condition
\begin{equation}\label{5.3}
(\tilde{\rho}, \tilde{u}, \tilde{\theta})(0,x)=(\tilde{\rho}_0, \tilde{u}_0, \tilde{\theta}_0)(x).
%\in H^s(\R_+^3).
\end{equation}

\vspace{1.5mm}

For later use we define
\begin{equation}\nonumber
\bp^\alpha=\partial_t^{\a_0}\partial_{x_1}^{\a_1}\partial_{x_2}^{\a_2},
\end{equation}
where $\alpha$ is a vector index  which is different from the one defined in \eqref{relation of mu and muM}.
Define the nations $\|\cdot\|_{\mathcal{H}^k(\R_+^3)}$ and $\|\cdot\|_{\mathcal{H}^k(\R^2)}$:
\begin{equation}\label{5.4}
\|f(t)\|_{\mathcal{H}^k(\R_+^3)}^2=\sum_{|\alpha|+i\leq k} \|\bp^\alpha \pa_3^if(t)\|_{L^2(\R^3_+)}^2, \quad
\|g(t)\|_{\mathcal{H}^k(\R^2)}^2=\sum_{|\alpha|\leq k} \|\bp^\alpha g(t)\|_{L^2(\R^2)}^2.
\end{equation}

\begin{lemma}\label{lem3.1}
Let  $(\rho,\fu, T)$ be the smooth solution of compressible Euler system obtained in Lemma \ref{lem2.1-1}, and $\tau>0$ be its lifespan. We assume that
\begin{align}\label{5.4-1}
&\|(\tilde{\rho}_0, \tilde{u}_0, \tilde{\theta}_0)\|^2_{\MH^k(\R^3_+)}
+\sup_{t\in(0,\tau)} \Big[\|(\mathfrak{f},\mathfrak{g})(t)\|^2_{\MH^{k+1}(\R^3_+)} +\|d(t)\|^2_{\MH^{k+2}}\Big]<\infty,
\end{align}
with $k\geq 3$, and the compatibility condition is satisfied for the initial data.  Then there exists a  unique smooth solution to \eqref{5.1}-\eqref{5.3} for $t\in[0,\tau]$, and satisfies
\begin{align}\label{5.5}
\sup_{t\in[0,\tau]}\|(\tilde{\rho}, \tilde{u}, \tilde{\theta})(t)\|^2_{\MH^k(\R^3_+)}
&\leq C(\tau,E_{k+2}) \Big\{\|(\tilde{\rho}_0, \tilde{u}_0, \tilde{\theta}_0)\|^2_{\MH^k(\R^3_+)} \nonumber\\
&\quad+\sup_{t\in[0,\tau]} \Big[\|(\mathfrak{f},\mathfrak{g})(t)\|^2_{\MH^{k+1}(\R^3_+)} +\|d(t)\|^2_{\MH^{k+2}(\R^2)}\Big]\Big\},
\end{align}
where $E_{k}:=\sup_{t\in[0,\tau]}\|(\rho-1,\fu,T-1)(t)\|_{H^{k}}$.
\end{lemma}

\noindent{\bf Proof.}  We define
\begin{equation*}
\tilde{p}:=\f{\rho\tilde{\theta}+3T\tilde{\rho}}{3}.
\end{equation*}
To deal with the boundary terms, it is more convenient to use the variables $(\tilde{p},\tilde{u},\tilde{\theta})$. Then \eqref{5.1} is equivalent to
\begin{equation}\label{5.6}
\begin{cases}
\dis \pa_t\tilde{p}+\fu\cdot \nabla_x\tilde{p}+\frac53 p\, \mbox{\rm div}_x \tilde{u}+\frac53\mbox{\rm div}_x \fu \, \tilde{p}+\nabla_xp\cdot \tilde{u}=\frac13 \mathfrak{g},\\[3mm]
\dis  \rho\partial_t \tilde{u}+\rho\fu\cdot\nabla_x \tilde{u}+\nabla_x\tilde{p}-\frac{\nabla_xp}{p} \tilde{p}+\rho\tilde{u}\cdot \nabla_x \fu-\f{\nabla_xp}{3T}\tilde{\theta}=\mathfrak{f},\\[3mm]
\dis  \rho\pa_t\tilde{\theta}+\rho\fu\cdot\nabla_x\tilde{\theta}+2p\,\mbox{\rm div}_x\tilde{u}+3\rho\tilde{u}\cdot\nabla_xT+\f23\rho\mbox{\rm div}_x\fu \, \tilde{\theta}=\mathfrak{g}.
\end{cases}
\end{equation}
Let $\chi$ be a smooth monotonic cut-off function  such that
\begin{align}\label{5.7}
\chi(s)=
\begin{cases}
1,\quad s\in[0,1],\\
0,\quad s\in[2,\infty).
\end{cases}
\end{align}
Then we define
\begin{equation}\label{5.7-1}
u_d(t,x):=\big(0,0,d(t,x_\sp)\chi(x_3)\big)^\top\quad\mbox{and}\quad
\tilde{w}:=\tilde{u}-u_d.
\end{equation}
Now we can rewrite \eqref{5.6} to be
\begin{equation}\label{5.8}
\begin{cases}
\dis \pa_t\tilde{p}+\fu\cdot \nabla_x\tilde{p}+\frac53 p \, \mbox{\rm div}_x \tilde{w}+\frac53\mbox{\rm div}_x \fu \, \tilde{p}+\nabla_xp\cdot \tilde{w}\\[3mm]
\dis \qquad\qquad\qquad\qquad\qquad\qquad=G_0:=\frac13 \mathfrak{g}-\frac53 p\, \mbox{\rm div} u_d-\nabla_xp\cdot u_d,\\[3mm]
\dis  \rho\partial_t \tilde{w}+\rho\fu\cdot\nabla_x \tilde{w}+\nabla_x\tilde{p}-\frac{\nabla_xp}{p} \tilde{p}+\rho\tilde{w}\cdot \nabla_x \fu-\f{\nabla_xp}{3T}\tilde{\theta}\\[3mm]
\dis \qquad\qquad\qquad\qquad\qquad\qquad=G_1:=\mathfrak{f}-\rho \partial_t u_d-\rho\fu\cdot\nabla_x u_d-\rho u_d\cdot \nabla_x \fu,\\[3mm]
\dis  \rho\pa_t\tilde{\theta}+\rho\fu\cdot\nabla_x\tilde{\theta}+2p\, \mbox{\rm div}_x\tilde{w}+3\rho\tilde{w}\cdot\nabla_xT+\f23\rho\mbox{\rm div}_x\fu\, \tilde{\theta}\\[3mm]
\dis \qquad\qquad\qquad\qquad\qquad\qquad=G_2:=\mathfrak{g}-2p\, \mbox{\rm div}_x u_d-3\rho u_d\cdot\nabla_xT.
\end{cases}
\end{equation}
From \eqref{5.2}, the boundary condition now becomes
\begin{equation}\label{5.9}
\tilde{w}_3(t,x_\sp,0)\equiv0.
\end{equation}
We can write the linear system \eqref{5.8} as a symmetric hyperbolic equations
\begin{align}\label{5.10}
A_0 \partial_t U+\sum_{i=1}^3 A_i \partial_i U+A_4U=G,
\end{align}
where
\begin{equation*}%\label{5.11}
U=\left(
\begin{array}{c}
\tilde{p}\\
\tilde{w}\\
\tilde{\theta}\\
\end{array}
\right),
\quad
A_0= \left(
\begin{array}{ccccc}
\frac95 & 0 & 0 & 0 & -\rho \\
0 & \rho p & 0 & 0 & 0\\
0 & 0 &\rho p & 0 & 0 \\
0 & 0 & 0 & \rho p& 0 \\
-\rho & 0 & 0 & 0 & \f{5}{6} \\
\end{array}
\right),
\quad
A_1= \left(
\begin{array}{ccccc}
\frac95 \fu_1 & p & 0 & 0 & -\rho \fu_1 \\
p & \rho p\fu_1 & 0 & 0 & 0\\
0 & 0 &\rho p\fu_1 & 0 & 0 \\
0 & 0 & 0 & \rho p\fu_1& 0 \\
-\rho \fu_1 & 0 & 0 & 0 & \f{5}{6}\rho^2 \fu_1 \\
\end{array}
\right),
\end{equation*}
and
\begin{equation}\label{5.12}
A_2= \left(
\begin{array}{ccccc}
\frac95 \fu_2 & 0 & p & 0 & -\rho \fu_2 \\
0 & \rho p\fu_1 & 0 & 0 & 0\\
p & 0 &\rho p\fu_2 & 0 & 0 \\
0 & 0 & 0 & \rho p\fu_2& 0 \\
-\rho \fu_2 & 0 & 0 & 0 & \f{5}{6}\rho^2 \fu_2 \\
\end{array}
\right),
\quad
A_3= \left(
\begin{array}{ccccc}
\frac95 \fu_3 & 0 & 0 & p & -\rho \fu_3 \\
0 & \rho p\fu_3 & 0 & 0 & 0\\
0 & 0 &\rho p\fu_3 & 0 & 0 \\
p & 0 & 0 & \rho p\fu_3& 0 \\
-\rho \fu_3 & 0 & 0 & 0 & \f{5}{6}\rho^2 \fu_3 \\
\end{array}
\right).
\end{equation}
The matrix $A_4$ and column vector $G$ can be easily write down, and we do not give the details here. It is easy to check that $A_0$ is positive.

Noting the matric $A_3$ is singular on the boundary, hence  the IBVP \eqref{5.8}-\eqref{5.9} is a linear hyperbolic system with characteristic boundary.  We refer \cite{Secchi,Chenshuxing-1980} for the  local existence of smooth solution. To close our lemma, one needs only to establish the {\it a priori} energy estimates.

It follows from Newtonian-Leibnitz formula that
\begin{align}\label{5.16-1}
\|U(t)\|_{\MH^{k-1}}^2&\leq \|U_0\|_{\MH^{k-1}}^2+2\int_0^t\|\pa_tU(s)\|_{\MH^{k-1}}\|U(s)\|_{\MH^{k-1}}ds\nonumber\\
&\leq \|U_0\|_{\MH^{k-1}}^2+\int_0^t\|U(s)\|_{\MH^{k}}^2ds.
\end{align}
Hence we need only to close the highest order derivatives estimates.  Let $|\alpha|+i=k$, and applying $\bp^\alpha \pa_3^i$ to \eqref{5.10}, we obtain
\begin{align}\label{5.14}
&A_0 \partial_t \bp^\alpha \pa_3^i U+\sum_{i=1}^3 A_i \partial_i \bp^\alpha \pa_3^iU\nonumber\\
&=\bp^\alpha \pa_3^iG-\bp^\alpha \pa_3^i(A_4U)-[\bp^\alpha \pa_3^i,A_0] \pa_t U
-\sum_{j=1}[\bp^\alpha \pa_3^i,A_j] \pa_j U.
\end{align}
where and whereafter the notation $[\cdot, \, \cdot]$ denote the commutator operator, i.e.,
\begin{align}\nonumber
[\partial^{\alpha},f]g=\pa^\alpha (fg)-f\pa^\alpha g=\sum_{\beta+\gamma=\alpha, |\beta|\geq1} C_{\beta,\gamma}\pa^\beta f\cdot \pa^{\gamma} g.
\end{align}
Multiplying \eqref{5.14} by $\bp^\alpha \pa_3^i U^{\top}$ and integrating the resultant equation over $[0,t]\times\R^3_+$,  we obtain
\begin{align}\label{5.16}
\|\bp^\alpha \pa_3^i  U(t)\|_{L^2}^2&\leq C\|\bp^\alpha \pa_3^i  U(0)\|_{L^2}^2+C\left|\int_0^t\int_{\R^2} (\bp^\alpha \pa_3^i U^{\top} A_3 \bp^\alpha \pa_3^i U)(s,x_\sp,0) dx_\sp ds\right|\nonumber\\
&\quad+C(E_{k+1})\int_0^t\|(U,G)(s)\|_{\MH^k}^2ds.
\end{align}
For the boundary term on RHS of \eqref{5.16}, noting \eqref{5.12} and $\fu_3(t,x_\sp,0)\equiv0$,  it holds
\begin{equation}\label{5.17}
\int_0^t\int_{\R^2} (\bp^\alpha \pa_3^i U^{\top} A_3 \bp^\alpha \pa_3^i U)(s,x_\sp,0)\, dx_\sp ds
=2\int_0^t\int_{\R^2} (p\, \bp^\alpha \pa_3^i \tilde{p}\cdot\bp^\alpha \pa_3^i \tilde{w}_3)(s,x_\sp,0)\, dx_\sp ds.
\end{equation}

\smallskip

To close the above estimates, we  use induction argument on the number of normal derivatives $\pa_3^i$.   For $i=0$, it follows from \eqref{5.9} that
\begin{equation}\label{5.18}
\bp^\alpha \tilde{w}_3(t,x_\sp,0)\equiv0,
\end{equation}
which, together with \eqref{5.16} and \eqref{5.17},  yields  that
%\begin{equation}
%\int_0^t\int_{\R^2} (U^{\top} A_3 U)(t,x_\sp,0) dx_\sp=2\int_0^t\int_{\R^2} (p\tilde{p} \tilde{w}_3)(t,x_\sp,0) dx_\sp\equiv0.
%\end{equation}
\begin{align*}%\label{5.19}
\sum_{|\alpha|=k}\|\bp^\alpha U(t)\|_{L^2}^2&\leq C\|U(0)\|_{\MH^k}^2+C(E_{k+1})\int_0^t\|(U,G)(s)\|_{\MH^k}^2ds.
\end{align*}
Assume that we have already obtained
\begin{align}\label{5.20}
\sum_{|\alpha|+i=k,i\leq l-1}\|\bp^\alpha\pa_3^i U(t)\|_{L^2}^2&\leq C\|(U,G)(0)\|_{\MH^k}^2+C(E_{k+1})\int_0^t\|U(s)\|_{\MH^k}^2+\|G(s)\|_{\MH^{k+1}}^2ds.
\end{align}
Next, we shall consider the case for $\|\bp^\alpha\pa_3^l U(t)\|_{L^2}^2$ with $|\alpha|+l=k$. Noting \eqref{5.16} and \eqref{5.17}, we need only to control the boundary term
\begin{equation}\label{5.21}
\int_0^t\int_{\R^2} (p\, \bp^\alpha \pa_3^l \tilde{p}\cdot\bp^\alpha \pa_3^l \tilde{w}_3)(s,x_\sp,0)\, dx_\sp ds.
\end{equation}
It follows from \eqref{5.8} that
\begin{align}
\pa_3\tilde{p}&=- \rho(\partial_t+\fu\cdot\nabla_x) \tilde{w}_3+\frac{\pa_3p}{p} \tilde{p}-\rho\tilde{w}\cdot \nabla_x \fu_3+\f{\partial_3p}{3T}\tilde{\theta}+G_{1,3},\label{5.22}\\[1mm]
p\partial_3\tilde{w}_3&=-\frac35(\partial_t+\fu\cdot\nabla_x) \tilde{p}-p \sum_{i=1}^2 \pa_i \tilde{w}_{i}-\mbox{\rm div}_x \fu \, \tilde{p}-\frac35\nabla_xp\cdot \tilde{w}+\frac35 G_0.\label{5.23}
\end{align}
By utilizing \eqref{1.7} and \eqref{1.12}, we have
\begin{equation}\label{5.24}
\pa_3 p(t,x_\sp,0)\equiv0\quad \mbox{and}\quad \fu_3(t,x_\sp,0)\equiv0.
\end{equation}
Substituting \eqref{5.24} into \eqref{5.22} and using \eqref{5.18}, one obtains
\begin{equation}\label{5.25}
\pa_3\tilde{p}(t,x_\sp,0)=G_{1,3}(t,x_\sp,0)
\end{equation}

Applying $\partial_3^{l-1}$ to \eqref{5.22} and \eqref{5.23}, and using \eqref{5.24}, then we  have
\begin{align}
\pa_3^l \tilde{p}(t,x_\sp,0)&\cong (\partial_t+\fu\cdot\nabla_x) \pa_3^{l-1}\tilde{w}_3+\partial_3^{l-1} \tilde{w}_3+\pa_t \partial_3^{l-2} \tilde{w}_3\nonumber\\
&\quad+\sum_{i=1}^2\sum_{|\beta|+j\leq l-2} \pa_i\bp^{\beta}\pa_3^j \tilde{w}_3+\sum_{|\beta|+j\leq l-2} \bp^{\beta}\pa_3^j (\tilde{p},\tilde{w},\tilde{\theta})+\partial_3^{l-1}G_{1,3},\label{5.26}\\
\partial_3^{l} \tilde{w}_3(t,x_\sp,0)&\cong (\partial_t+\fu\cdot\nabla_x) \pa_3^{l-1}\tilde{p}+\sum_{i=1}^2 \{\pa_i \pa_3^{l-1}\tilde{w}_{i}+\pa_3^{l-1}\tilde{w}_{i}\}+\partial_3^{l-1}\tilde{p}\nonumber\\
&\quad +\sum_{i=1}^2 \partial_i\partial_3^{l-2}( \tilde{p},\tilde{w})+\sum_{|\beta|+j\leq l-2} \bp^{\beta}\pa_3^j (\tilde{p},\tilde{w})+\pa_3^{l-1}G_0.\label{5.27}
\end{align}
where and whereafter  we  ignore the exact coefficients which depends only on the Euler solution. Substituting \eqref{5.26} into \eqref{5.27}, one can obtain
\begin{align}
\partial_3^{l} \tilde{w}_3(t,x_\sp,0)&\cong (\partial_t+\fu\cdot\nabla_x)\bp  \pa_3^{l-2}\tilde{w}_3+\sum_{i=1}^2 \{\pa_i \pa_3^{l-1}\tilde{w}_{i}+\pa_3^{l-1}\tilde{w}_{i}\}\nonumber\\
&\quad+\pa_t\pa_3^{l-2} \tilde{w}_3+\pa_t^2\pa_3^{l-3} \tilde{w}_3+\sum_{i=1}^2\sum_{|\beta|+j\leq l-2} \pa_i\bp^{\beta}\pa_3^j (\tilde{p},\tilde{w})\nonumber\\
&\quad +\sum_{|\beta|+j\leq l-2} \bp^{\beta}\pa_3^j (\tilde{p},\tilde{w},\tilde{\theta})+\pa^{l-1}G.\label{5.29}
\end{align}
If  $l$ is even, using \eqref{5.18} and \eqref{5.29},  step by step, one can get
\begin{align}\label{5.30}
\partial_3^{l} \tilde{w}_3(t,x_\sp,0)&\cong \sum_{i=1}^2\sum_{j=0}^{\frac{l}{2}-1} \{\pa_i \bp^{2j} \pa_3^{l-1-2j}\tilde{w}_{i}+\bp^{2j}\pa_3^{l-1-2j}\tilde{w}_{i}\} +\partial_t^{l-1}\tilde{p}\nonumber\\
&\hspace{-2mm}+\sum_{i=1}^2\sum_{|\beta|+j\leq l-2} \pa_i\bp^{\beta}\pa_3^j (\tilde{p},\tilde{w})+\sum_{|\beta|+j\leq l-2} \bp^{\beta}\pa_3^j (\tilde{p},\tilde{w},\tilde{\theta})+\sum_{j=0}^{l-1}\pa^{j}G.
\end{align}
Similarly, if $l$ is odd, step by step, we  have
\begin{align}\label{5.31-1}
\partial_3^{l} \tilde{w}_3(t,x_\sp,0)&\cong (\partial_t+\fu\cdot\nabla_x)\bp^{l-2} \pa_3\tilde{w}_3+ \sum_{i=1}^2\sum_{j=0}^{\frac{l-3}{2}} \{\pa_i \bp^{2j} \pa_3^{l-1-2j}\tilde{w}_{i}+\bp^{2j}\pa_3^{l-1-2j}\tilde{w}_{i}\} \nonumber\\
&\hspace{-3mm}+\partial_t^{l-1}\tilde{p}+\sum_{i=1}^2\sum_{|\beta|+j\leq l-2} \pa_i\bp^{\beta}\pa_3^j (\tilde{p},\tilde{w})+\sum_{|\beta|+j\leq l-2} \bp^{\beta}\pa_3^j (\tilde{p},\tilde{w},\tilde{\theta})+\sum_{j=0}^{l-1}\pa^{j}G.
\end{align}
Substituting \eqref{5.23} into \eqref{5.31-1},   one obtains,  for $l$ being odd, that
\begin{align}\label{5.31}
&\partial_3^{l} \tilde{w}_3(t,x_\sp,0)\cong (\partial_t+\fu\cdot\nabla_x)\bp^{l-1} \tilde{p}+ \sum_{i=1}^2\sum_{j=0}^{\frac{l-1}{2}} \{\pa_i \bp^{2j} \pa_3^{l-1-2j}\tilde{w}_{i}+\bp^{2j}\pa_3^{l-1-2j}\tilde{w}_{i}\} \nonumber\\
&\qquad\quad+\partial_t^{l-1}\tilde{p}+\sum_{i=1}^2\sum_{|\beta|+j\leq l-2} \pa_i\bp^{\beta}\pa_3^j (\tilde{p},\tilde{w})+\sum_{|\beta|+j\leq l-2} \bp^{\beta}\pa_3^j (\tilde{p},\tilde{w},\tilde{\theta})+\sum_{j=0}^{l-1}\pa^{j}G.
\end{align}

To estimate $\pa_3^l \tilde{p}$, we have to be careful. Let $l$ be even. Substituting \eqref{5.30} and \eqref{5.31} into \eqref{5.26}, we can get
\begin{align}\label{5.32-1}
\pa_3^l \tilde{p}(t,x_\sp,0)&\cong (\partial_t+\fu\cdot\nabla_x)\bp^{l-1}\tilde{p}+ (\partial_t+\fu\cdot\nabla_x) \sum_{i=1}^2\sum_{j=0}^{\frac{l}{2}-1} \{\pa_i \bp^{2j} \pa_3^{l-2-2j}\tilde{w}_{i}+\bp^{2j}\pa_3^{l-2-2j}\tilde{w}_{i}\}\nonumber\\
&\ +\pa_t^{l-1}\tilde{p}+\sum_{i=1}^2\sum_{|\beta|+j\leq l-2} \pa_i\bp^{\beta}\pa_3^j (\tilde{p},\tilde{w})+\sum_{|\beta|+j\leq l-2} \bp^{\beta}\pa_3^j (\tilde{p},\tilde{w},\tilde{\theta})+\sum_{j=0}^{l-1}\partial^{j}G\nonumber\\
&\cong (\partial_t+\fu\cdot\nabla_x)\bp^{l-1}\tilde{p}+ \sum_{i=1}^2\sum_{j=0}^{\frac{l}{2}-1} \pa_i^2 \bp^{2j} \pa_3^{l-2-2j}\tilde{p}%+\bp^{2j}\pa_3^{l-2-2j}\tilde{p}
+\pa_t^{l-1}\tilde{p}\nonumber\\
&+\sum_{i=1}^2\sum_{|\beta|+j\leq l-2} \pa_i\bp^{\beta}\pa_3^j (\tilde{p},\tilde{w}, \tilde{\theta})+\sum_{|\beta|+j\leq l-2} \bp^{\beta}\pa_3^j (\tilde{p},\tilde{w},\tilde{\theta})+\sum_{j=0}^{l-1}\partial^{j}G,
\end{align}
where we have used the facts
\begin{equation}\label{5.33}
(\partial_t+\fu\cdot\nabla_x) \tilde{w}_i\cong \partial_i \tilde{p}+O(1) (\tilde{p},\tilde{w},\tilde{\theta}) +G_{1,i}.
\end{equation}
which can be derived from \eqref{5.8}.

\smallskip

Iterating \eqref{5.32-1} again, step by step, we have, for $l$ being even, that
\begin{align*}%\label{5.32}
\pa_3^l \tilde{p}(t,x_\sp,0)&\cong  (\partial_t+\fu\cdot\nabla_x)\bp^{l-1} \tilde{p}+\pa_t^{l-1}\tilde{p}+\sum_{i=1}^2\sum_{|\beta|+j\leq l-2} \pa_i\bp^{\beta}\pa_3^j (\tilde{p},\tilde{w},\tilde{\theta})\nonumber\\
&\quad+\sum_{|\beta|+j\leq l-2} \bp^{\beta}\pa_3^j (\tilde{p},\tilde{w},\tilde{\theta})+\sum_{j=0}^{l-1}\partial^{j}G,
\end{align*}

For $l$ being odd, substituting   \eqref{5.30}  and \eqref{5.31} into \eqref{5.26}, and using \eqref{5.33}, we  get
\begin{align}\label{5.47-1}
\pa_3^l\tilde{p}(t,x_\sp,0)&\cong(\partial_t+\fu\cdot\nabla_x) \sum_{i=1}^2\sum_{j=0}^{\frac{l-3}{2}} \{\pa_i \bp^{2j} \pa_3^{l-2-2j}\tilde{w}_{i}+\bp^{2j}\pa_3^{l-2-2j}\tilde{w}_{i}\}+\pa_t^{l-1}\tilde{p}\nonumber\\
&\quad+\sum_{i=1}^2\sum_{|\beta|+j\leq l-2} \pa_i\bp^{\beta}\pa_3^j (\tilde{p},\tilde{w},\tilde{\theta})+\sum_{|\beta|+j\leq l-2} \bp^{\beta}\pa_3^j (\tilde{p},\tilde{w},\tilde{\theta})+\sum_{j=0}^{l-1}\partial^{j}G\nonumber\\
&\cong \sum_{i=1}^2\sum_{j=0}^{\frac{l-3}{2}} \pa_i^{2} \bp^{2j} \pa_3^{l-2-2j}\tilde{p} + \pa_t^{l-1}\tilde{p}+\sum_{i=1}^2\sum_{|\beta|+j\leq l-2} \pa_i\bp^{\beta}\pa_3^j (\tilde{p},\tilde{w},\tilde{\theta})\nonumber\\
&\quad+\sum_{|\beta|+j\leq l-2} \bp^{\beta}\pa_3^j (\tilde{p},\tilde{w},\tilde{\theta})+\sum_{j=0}^{l-1}\partial^{j}G.
\end{align}
Then, iterating \eqref{5.47-1} again, step by step, one obtains, for $l$ being odd, that
\begin{align}\label{5.47}
&\pa_3^l\tilde{p}(t,x_\sp,0)\cong \bp^{l-1}\pa_3 \tilde{p}+\pa_t^{l-1}\tilde{p}+\sum_{i=1}^2\sum_{|\beta|+j\leq l-2} \pa_i\bp^{\beta}\pa_3^j (\tilde{p},\tilde{w},\tilde{\theta})\nonumber\\
&\qquad\qquad\qquad+\sum_{|\beta|+j\leq l-2} \bp^{\beta}\pa_3^j (\tilde{p},\tilde{w},\tilde{\theta})+\sum_{j=0}^{l-1}\partial^{j}G\nonumber\\
&\cong \pa_t^{l-1}\tilde{p}+\sum_{i=1}^2\sum_{|\beta|+j\leq l-2} \pa_i\bp^{\beta}\pa_3^j (\tilde{p},\tilde{w},\tilde{\theta})+\sum_{|\beta|+j\leq l-2} \bp^{\beta}\pa_3^j (\tilde{p},\tilde{w},\tilde{\theta})+\sum_{j=0}^{l-1}\partial^{j}G,
\end{align}
where we have used \eqref{5.25} in the last step.

\smallskip

Now we estimate the boundary term \eqref{5.21} when  $l$ is  even. Using integration by parts, \eqref{5.26}, \eqref{5.30} and Lemma \ref{lemA.1},  it holds that
\begin{align}\label{5.48}
&\int_0^t\int_{\R^2} (p \bp^\alpha \pa_3^l \tilde{p}\cdot\bp^\alpha \pa_3^l \tilde{w}_3)(s,x_\sp,0) dx_\sp ds\nonumber\\
&\cong\int_0^t\int_{\R^2}  \Big\{\sum_{i=1}^2\sum_{j=0}^{\frac{l}{2}-1} \{\pa_i \bp^{\alpha+2j} \pa_3^{l-1-2j}\tilde{w}_{i}+\bp^{\alpha+2j}\pa_3^{l-1-2j}\tilde{w}_{i}\} +\bp^{\alpha}\partial_t^{l-1}\tilde{p}\nonumber\\
&\quad +\sum_{i=1}^2\sum_{|\beta|+j\leq l-2} \pa_i\bp^{\alpha+\beta}\pa_3^j (\tilde{p},\tilde{w})+\sum_{|\beta|+j\leq l-2} \bp^{\alpha+\beta}\pa_3^j (\tilde{p},\tilde{w},\tilde{\theta})+\sum_{j=0}^{l-1}\bp^{\alpha}\pa^{j}G\Big\}\nonumber\\
&\qquad \times \Big\{(\partial_t+\fu\cdot\nabla_x) \bp^{\alpha}\pa_3^{l-1}\tilde{w}_3+\bp^{\alpha}\partial_3^{l-1} \tilde{w}_3+\pa_t \bp^{\alpha}\partial_3^{l-2} \tilde{w}_3\nonumber\\
&\quad+\sum_{i=1}^2\sum_{|\beta|+j\leq l-2} \pa_i\bp^{\alpha+\beta}\pa_3^j \tilde{w}_3+\sum_{|\beta|+j\leq l-2} \bp^{\alpha+\beta}\pa_3^j (\tilde{p},\tilde{w},\tilde{\theta})+\bp^{\alpha}\partial_3^{l-1}G\Big\} dx_{\sp}  ds\nonumber\\
&\cong \int_0^t\int_{\R^2}  \bp^{\alpha}\pa_3^{l-1}\tilde{w}_3\times (\partial_t+\fu\cdot\nabla_x)\Big\{\sum_{i=1}^2\sum_{j=0}^{\frac{l}{2}-1} \{\pa_i \bp^{\alpha+2j} \pa_3^{l-1-2j}\tilde{w}_{i}+\bp^{\alpha+2j}\pa_3^{l-1-2j}\tilde{w}_{i}\} \nonumber\\
&\ +\bp^{\alpha}\partial_t^{l-1}\tilde{p}+\sum_{i=1}^2\sum_{|\beta|+j\leq l-2} \pa_i\bp^{\alpha+\beta}\pa_3^j (\tilde{p},\tilde{w})+\sum_{|\beta|+j\leq l-2} \bp^{\alpha+\beta}\pa_3^j (\tilde{p},\tilde{w},\tilde{\theta})+\sum_{j=0}^{l-1}\bp^{\alpha}\pa^{j}G_0\Big\}dx_\sp ds\nonumber\\
&\quad+\int_{\R^2}  \bp^{\alpha}\pa_3^{l-1}\tilde{w}_3\times \Big\{\sum_{i=1}^2\sum_{j=0}^{\frac{l}{2}-1} \{\pa_i \bp^{\alpha+2j} \pa_3^{l-1-2j}\tilde{w}_{i}+\bp^{\alpha+2j}\pa_3^{l-1-2j}\tilde{w}_{i}\} +\bp^{\alpha}\partial_t^{l-1}\tilde{p}\nonumber\\
&\quad +\sum_{i=1}^2\sum_{|\beta|+j\leq l-2} \pa_i\bp^{\alpha+\beta}\pa_3^j (\tilde{p},\tilde{w})+\sum_{|\beta|+j\leq l-2} \bp^{\alpha+\beta}\pa_3^j (\tilde{p},\tilde{w},\tilde{\theta})+\sum_{j=0}^{l-1}\bp^{\alpha}\pa^{j}G\Big\}dx_\sp\Big|_{0}^t\nonumber\\
&\quad + C(E_{k+1})\int_0^t\|(U,G)(s)\|_{\MH^k}^2ds.
\end{align}
By using \eqref{5.16-1}, \eqref{5.20} and Lemma \ref{lemA.1}, the second term on RHS of \eqref{5.48}  can be bounded by
\begin{align}\label{5.49}
&\lambda \sum_{|\alpha|=k-l}\|\bp^\alpha\pa_3^l U(t)\|_{L^2}^2+C_\lambda\sum_{|\alpha|+i=k,i\leq l-1}\|\bp^\alpha\pa_3^i U(t)\|_{L^2}^2\nonumber\\
&\qquad\qquad\qquad\quad+\|U(t)\|^2_{\MH^{k-1}}+\|(U,G)(0)\|_{\MH^k}^2+\|G(t)\|^2_{\MH^k}\nonumber\\
&\leq \lambda \sum_{|\alpha|=k-l}\|\bp^\alpha\pa_3^l U(t)\|_{L^2}^2+C(E_{k+1}) \|(U,G)(0)\|_{\MH^k}^2\nonumber\\
&\qquad+C(E_{k+1})\int_0^t\|U(s)\|_{\MH^k}^2+\|G(s)\|_{\MH^{k+1}}^2ds.
\end{align}
Integrating by parts   with respect to $x_1$ or $x_2$ for the highest order terms, using \eqref{5.33} and Lemma \ref{lemA.1}, the first term on RHS of \eqref{5.48} is bounded by
\begin{align}\label{5.50}
%&\int_0^t\int_{\R^2}  \bp^{\alpha}\pa_3^{l-1}\tilde{w}_3\times (\partial_t+\fu\cdot\nabla_x)\Big\{\sum_{i=1}^2\sum_{j=0}^{\frac{l}{2}-1} \{\pa_i \bp^{\alpha+2j} \pa_3^{l-1-2j}\tilde{w}_{i}+\bp^{\alpha+2j}\pa_3^{l-1-2j}\tilde{w}_{i}\} \nonumber\\
%&\quad +\bp^{\alpha}\partial_t^{l-1}\tilde{p}+\sum_{i=1}^2\sum_{|\beta|+j\leq l-2} \pa_i\bp^{\alpha+\beta}\pa_3^j (\tilde{p},\tilde{w})+\sum_{|\beta|+j\leq l-2} \bp^{\alpha+\beta}\pa_3^j (\tilde{p},\tilde{w},\tilde{\theta})+\bp^{\alpha}\pa^{l-1}G_0\Big\}dx_\sp ds\nonumber\\
&\sum_{i=1}^2\int_0^t\int_{\R^2}  \pa_i\bp^{\alpha}\pa_3^{l-1}\tilde{w}_3\times \sum_{j=0}^{\frac{l}{2}-1} \{\pa_i \bp^{\alpha+2j} \pa_3^{l-1-2j}\tilde{p}+\bp^{\alpha+2j} \pa_3^{l-1-2j}\tilde{p} \} \nonumber\\
&\quad+\int_0^t\int_{\R^2}  \bp^{\alpha}\pa_3^{l-1}\tilde{w}_3\cdot \bp^{\alpha}\partial_t^{l}\tilde{p} dx_\sp ds +C(E_{k+1})\int_0^t\|U(s)\|_{\MH^k}^2+\|G(s)\|_{\MH^{k+1}}^2ds\nonumber\\
& \cong \sum_{i=1}^2\int_0^t\int_{\R^2}  \pa_i\bp^{\alpha}\pa_3^{l-1}\tilde{w}_3\times \sum_{j=0}^{\frac{l}{2}-1} \pa_i \bp^{\alpha+2j} \pa_t^{l-2-2j}\tilde{p} dx_{\sp}ds\nonumber\\
&\quad+\int_0^t\int_{\R^2}  \bp^{\alpha}\pa_3^{l-1}\tilde{w}_3\cdot \bp^{\alpha}\partial_t^{l}\tilde{p} dx_{\sp} ds+C(E_{k+1})\int_0^t\|U(s)\|_{\MH^k}^2+\|G(s)\|_{\MH^{k+1}}^2ds\nonumber\\
&\cong \int_0^t\int_{\R^2}  \bp^{\alpha}\pa_3^{l-1}\tilde{w}_3\cdot \bp^{\alpha}\partial_t^{l}\tilde{p} dx_{\sp} ds+C(E_{k+1})\int_0^t\|U(s)\|_{\MH^k}^2+\|G(s)\|_{\MH^{k+1}}^2ds.
\end{align}
where we have used \eqref{5.47} for $\pa_3^{l-1-2j}\tilde{p}$.

If $\bp^\alpha=\partial_t^{\alpha_0}\partial_{x_1}^{\alpha_1}\partial_{x_2}^{\alpha_2}$ with $\alpha_1+\alpha_2\geq 1$, then by using Lemma \ref{lemA.1},  the first term on RHS of \eqref{5.50} is controlled by
\begin{equation*}%\label{5.51}
\int_0^t\int_{\R^2}  \bp^{\alpha}\pa_3^{l-1}\tilde{w}_3\cdot \bp^{\alpha}\partial_t^{l}\tilde{p} dx_{\sp} ds\leq \int_0^t\|U(s)\|_{\MH^k}^2ds.
\end{equation*}
The remaining case is  $\bp^\alpha\tilde{p} =\partial_t^{k-l} \tilde{p}$,  it follows from  \eqref{5.31} and Lemma \ref{lemA.1} that
\begin{align}\label{5.52}
&\int_0^t\int_{\R^2}  \pa_t^{k-l}\pa_3^{l-1}\tilde{w}_3\cdot \partial_t^{k}\tilde{p} dx_{\sp} ds\nonumber\\
&\cong C\int_0^t\|U(s)\|_{\MH^k}^2ds+\int_0^t\int_{\R^2}  \Big[\pa_t^{k-1}\tilde{p}+\sum_{i=1}^2\sum_{j\leq l-2} \pa_i\bp^{k-2-j}\pa_3^j (\tilde{p},\tilde{w})\nonumber\\
&\quad+\sum_{|\beta|\leq k-j-2,j\leq l-2} \bp^{\beta}\pa_3^j (\tilde{p},\tilde{w},\tilde{\theta})+\sum_{j=0}^{k-2}\pa^{j}G\Big]\cdot \partial_t^{k}\tilde{p} dx_{\sp} ds\nonumber\\
&\cong \int_0^t\|(U,G)(s)\|_{\MH^k}^2ds+\int_{\R^2} |\partial_t^{k-1}\tilde{p}|^2 dx_{\sp} \Big|_{0}^t+ \int_{\R^2}  \Big[\sum_{i=1}^2\sum_{j\leq l-2} \pa_i\bp^{k-2-j}\pa_3^j (\tilde{p},\tilde{w})\nonumber\\
&\qquad\qquad +\sum_{|\beta|\leq k-j-2,j\leq l-2} \bp^{\beta}\pa_3^j (\tilde{p},\tilde{w},\tilde{\theta})+\sum_{j=0}^{k-2}\pa^{j}G\Big]\cdot \partial_t^{k-1}\tilde{p} dx_{\sp} \Big|_0^t\nonumber\\
&\leq C(E_{k+1}) \Big\{\sum_{|\alpha|+j=k,j\leq l-1}\|\bp^\alpha\pa_3^j U(t)\|^2_{L^2}+\|(U,G)(0)\|_{\MH^k}^2\nonumber\\
&\qquad\qquad\qquad+\|G(t)\|^2_{\MH^k} +\int_0^t\|(U,G)(s)\|_{\MH^k}^2ds\Big\}.
\end{align}
Substituting \eqref{5.49}-\eqref{5.52} into \eqref{5.48}, for $l$ being even, we obtain
\begin{align}\label{5.53}
&\int_0^t\int_{\R^2} (p\, \bp^\alpha \pa_3^l \tilde{p}\cdot\bp^\alpha \pa_3^l \tilde{w}_3)(s,x_{\sp} ,0) dx_{\sp}  ds\nonumber\\
&\leq C(E_{k+1})\Big\{\|(U,G)(0)\|_{\MH^k}^2
+\int_0^t\|U(s)\|_{\MH^k}^2+\|G(s)\|_{\MH^{k+1}}^2ds\Big\}+\lambda \sum_{|\alpha|=k-l,}\|\bp^\alpha\pa_3^l U(t)\|^2_{L^2}.
\end{align}

For the case  $l$ being odd, by using \eqref{5.31} and \eqref{5.47}, one can also prove \eqref{5.53}, the proof is slight easier than even case,  and we omit the details here for simplicity of presentation. \vspace{1mm}

Combining \eqref{5.16}, \eqref{5.17} and \eqref{5.53}, and taking $\lambda$ small, we get
\begin{align*}%\label{5.54}
\sum_{|\alpha|=k-l}\|\bp^\alpha\pa_3^l U(t)\|_{L^2}^2
&\leq C(E_{k+1})\Big\{\|(U,G)(0)\|_{\MH^k}^2 +\int_0^t\|U(s)\|_{\MH^k}^2+\|G(s)\|_{\MH^{k+1}}^2ds\Big\}.
\end{align*}
This completes the induction argument. Therefore we can obtain
\begin{align}\nonumber
\|U(t)\|_{\MH^k}^2&\leq  C(E_{k+1})\Big\{\|(U,G)(0)\|_{\MH^k}^2
+\int_0^t\|U(s)\|_{\MH^k}^2+\|G(s)\|_{\MH^{k+1}}^2ds\Big\}.
\end{align}
which, together with Gronwall's inequality, yields that
\begin{equation}\label{5.48-1}
\|U(t)\|_{\MH^k}^2\leq C(E_{k+1})\left(\|(U,G)(0)\|_{\MH^k} +\int_0^t\|G(s)\|_{\MH^{k+1}}^2ds\right).
\end{equation}
Hence we conclude \eqref{5.5} by using \eqref{5.7-1}, \eqref{5.8} and \eqref{5.48-1}.  $\hfill\Box$

\

%%%%%%%%%%%%%%%%%%%%%%%%%%%%%%%%%%%%%%%%%%%%%%%%%%%%%%%%%%%%%%%%%%%%%%%%%%%%%%%%%%%%%%

\section{Existence of Solution for a Linear Parabolic System}\label{section4}
To construct the solution  of viscous boundary layer, we consider the following linear parabolic system of $(u,\theta)=(u_1,u_2,\theta)$
\begin{align}\label{5.57}
\begin{cases}
\dis \rho^0 \partial_t u_{i}+\rho^0 (\fu^0_{\sp} \cdot\nabla_{\sp} ) u_{i}+\rho^0\partial_3\fu_3^0\cdot y\partial_y u_{i}\\
\dis\qquad\qquad\qquad\quad+\rho^0 u\cdot \nabla_{\sp}  \fu^0_i
-\frac{\partial_i p^0}{3T^0} \theta=\mu(T^0) \partial_{yy} u_{i}+\mathfrak{f}_i,\ i=1,2,\\[3mm]
\dis\rho^0 \partial_t \theta+\rho^0 (\fu^0_{\sp} \cdot\nabla_{\sp} ) \theta+\rho^0\partial_3\fu_3^0\cdot y\partial_y\theta+\frac23\rho^0 \mbox{\rm div}\fu^0 \theta=\frac35\kappa(T^0) \partial_{yy}\theta+\mathfrak{g},
\end{cases}
\end{align}
where $(t,x_{\sp} ,y)\in[0,\tau]\times\R^2\times\R_+$, and $(\rho^0,\fu^0,T^0,\partial_3\fu_3^0, \mbox{\rm div}\fu^0,  \pa_1 p^0, \pa_2 p^0)$ are the corresponding values of Euler solution on the boundary $x_3=0$, which is independent of $y\in\R_+$. We impose the system \eqref{5.57} with non-homogenous Neumann boundary conditions, i.e.,
\begin{align}\label{5.58}
\begin{cases}
\pa_y u_i(t,x_{\sp} ,y)|_{y=0}=b_i(t,x_{\sp} ),\quad \pa_y \theta(t,x_{\sp} ,y)|_{y=0}=a(t,x_{\sp} ),\\[1mm]
\lim_{y\rightarrow\infty} (u,\theta)(t,x_{\sp} ,y)=0.
\end{cases}
\end{align}
%where $a$ and $b$ are some given functions.
We also impose \eqref{5.57} with initial data
\begin{align}\label{5.59}
u(t,x_{\sp} ,y)|_{t=0}=u_0(x_{\sp} ,y),\quad \theta(t,x_{\sp} ,0)|_{t=0}=\theta_0(x_{\sp} ,y).
\end{align}
The initial data $(u_0,\theta_0)$ should satisfies the corresponding compatibility condition.

Let $l\geq 0$, we define the  notations %$\|\cdot\|_{L^{2}_l}$  and $\|\cdot\|_{H^{k}_l}$:
\begin{equation}\label{5.59-1}
\|f\|_{L^2_l}^2=\iint(1+y)^l |f(x_{\sp} ,y)|^2 dx_{\sp}  dy,
% \quad \mbox{and}\quad
%\|f\|_{H^k_l}^2??=\sum_{j=0}^k \|\nabla^{j}f\|_{L^2_{l_j}}^2,
\end{equation}
and
\begin{equation}\label{5.59-2}
\bar{x}:=(x_{\sp} ,y),\quad \mbox{and}\quad \nabla_{\bar{x}}:=(\nabla_{\sp} ,\pa_y)\equiv (\pa_{x_1},\pa_{x_2},\pa_y).
\end{equation}

\begin{lemma}\label{lem4.1}
Let $l\geq 0$, $k\geq 3$, and the compatibility condition for the initial data \eqref{5.59} is satisfied.  Assume that
\begin{equation*}
\sup_{t\in[0,\tau]} \Big\{
\sum_{\beta+2\gamma\leq k+2}  \|\nabla^{\beta}_{\sp}  \pa_t^{\gamma}(a,b)(t)\|^2_{L^2(\R^2)}
+\sum_{j=0}^k\sum_{\beta+2\gamma=j}  \| \nabla_{\bar{x}}^{\beta} \pa_t^{\gamma} (\mathfrak{f},\ \mathfrak{g})(t)\|^2_{L^2_{l_j}}\Big\}<\infty.
\end{equation*}
with  $l_j:=l+2(k-j)$, $0\leq j\leq k$.
Then there exists a unique smooth solution $(u,\theta)$ of \eqref{5.57}-\eqref{5.59} over $t\in[0,\tau]$, and satisfies
\begin{align}\label{5.08}
&\sum_{j=0}^{k} \sum_{\beta+2\gamma=j} \sup_{t\in[0,\tau]}\left\{\|\pa_t^\gamma \nabla_{\bar{x}}^{\beta} (u,\theta)(t)\|^2_{L^2_{l_j}}+\int_0^t \|\pa_t^\gamma \nabla_{\bar{x}}^{\beta} \pa_y(u,\theta)\|_{L^2_{l_j}}^2ds\right\}\nonumber\\
&\leq C(\tau,E_{k+3}) \Bigg\{ \sum_{j=0}^{k} \sum_{\beta+2\gamma=j} \|\pa_t^\gamma \nabla_{\bar{x}}^{\beta} (u,\theta)(0)\|^2_{L^2_{l_j}}+ \sup_{t\in[0,\tau]}\Big[\sum_{\beta+2\gamma\leq k+2}  \|\nabla^{\beta}_{\sp}  \pa_t^{\gamma}(a,b)(t)\|^2_{L^2(\R^2)} \nonumber\\
&\qquad\qquad\qquad+\sum_{j=0}^k\sum_{\beta+2\gamma=j}
\| \nabla_{\bar{x}}^{\beta} \pa_t^{\gamma} (\mathfrak{f},\ \mathfrak{g})(t)\|^2_{L^2_{l_j}}\Big] \Bigg\},
\end{align}
where the notation $E_{k+3}$ is the one defined in Lemma \ref{lem3.1}.
\end{lemma}

\noindent{\bf Proof.} We define the following background functions
\begin{align}\label{5.09}
u_b:=y b(t,x_{\sp} ) \chi(y) \quad\mbox{and}\quad \theta_a:= y a(t,x_{\sp} ) \chi(y),
\end{align}
where $\chi$ is the smooth monotonic cut-off function defined in \eqref{5.7}. It is obvious to know that $u_b$ and $\theta_a$ are smooth function, and is compact with respect to $y$.

We define
\begin{align}\label{5.010}
\Psi= u-u_b \quad\mbox{and}\quad \Theta=\theta-\theta_a,
\end{align}
then \eqref{5.57} is reduced to
\begin{align}\label{5.63}
\begin{cases}
\dis \partial_t \Psi_{i}+(\fu^0_{\sp}\cdot\nabla_\sp) \Psi_{i}+\partial_3\fu_3^0\cdot y\partial_y \Psi_{i} +\Psi\cdot \nabla_{\sp} \fu^0_i
-\frac{\partial_i p^0}{3p^0} \Theta=\tilde{\mu}\partial_{yy} \Psi_{i}+\tilde{\mathfrak{f}}_i,\ i=1,2,\\[3mm]
\dis\partial_t \Theta+(\fu^0_{\s}\cdot\nabla_\sp) \Theta+\partial_3\fu_3^0\cdot y\partial_y\Theta+\frac2{3}\mbox{div}\fu^0\, \Theta=\tilde{\kappa} \partial_{yy}\Theta+\tilde{\mathfrak{g}},
\end{cases}
\end{align}
where $\dis\tilde{\mu}:=\frac{1}{\rho^0}\mu(T^0) $, $\dis\tilde{\kappa}:=\frac{3}{5\rho^0}\kappa(T^0) $, and
\begin{align}\label{5.011}
\begin{split}
\dis\tilde{\mathfrak{f}}_i&:=\frac{1}{\rho^0}\mathfrak{f}_i-\partial_t u_{b,i}-  (\fu^0_{\sp}\cdot\nabla_\sp) u_{b,i}-\partial_3\fu_3^0\cdot y\partial_y u_{b,i}-  u_b\cdot \nabla_\sp \fu^0_i
+\frac{\partial_i p^0}{3p^0} \theta_a+\tilde{\mu} \partial_{yy} u_{b,i},\\
\dis\tilde{\mathfrak{g}}&:=\frac{1}{\rho^0}\mathfrak{g}-\partial_t \theta_a-  (\fu^0_{\sp}\cdot\nabla_\sp) \theta_a-\partial_3\fu_3^0\cdot y\partial_y\theta_a-\frac2{3}\mbox{div}\fu^0\, \theta_a+\tilde{\kappa}\partial_{yy} \theta_a.
\end{split}
\end{align}
The boundary conditions \eqref{5.58}   becomes
\begin{align}\label{5.64}
\begin{cases}
\pa_y \Psi_i(t,x_{\sp} ,y)|_{y=0}=0,\quad \pa_y\Theta(t,x_{\sp} ,y)|_{y=0}=0,\  i=1,2,\\[1mm]
\lim_{y\rightarrow\infty} (\Psi,\Theta)(t,x_{\sp} ,y)=0.
\end{cases}
\end{align}
Noting the coefficient $\partial_3\fu_3^0\cdot y$ in \eqref{5.63} is singular as $y\rightarrow\infty$ and there are no horizontal viscous terms $\Delta_{\sp}  \Psi$ and $\Delta_{\sp}  \Theta$, we can not directly use the standard linear parabolic theory. To prove the existence of smooth solution to \eqref{5.63}-\eqref{5.64}, we divide the proof into several steps.

\vspace{2mm}

\noindent{\it Step 1. Approximate problem.}  We consider the following approximate problem
\begin{align}\label{5.66}
\begin{cases}
\dis \partial_t \Psi_{i}+(\fu^0_{\sp} \cdot\nabla_{\sp} ) \Psi_{i}+\partial_3\fu_3^0\cdot y\chi_\sigma(y)\partial_y \Psi_{i}\\
\dis\qquad\qquad\qquad\qquad+\Psi\cdot \nabla_{\sp}  \fu^0_i
-\frac{\partial_i p^0}{3p^0} \Theta=\tilde{\mu}\,\partial_{yy} \Psi_{i}+\x \Delta_{\sp} \Psi_i+\tilde{\mathfrak{f}}_i^{\sigma},\ i=1,2,\\[3mm]
\dis\partial_t \Theta+(\fu^0_{\sp} \cdot\nabla_{\sp} ) \Theta+ \partial_3\fu_3^0\cdot y\partial_y\Theta+\frac2{3} \mbox{\rm div}\fu^0 \Theta=\tilde{\kappa}\, \partial_{yy}\Theta+\x \Delta_{\sp}  \Theta+\tilde{\mathfrak{g}}^{\sigma},
\end{cases}
\end{align}
where $(t,x_{\sp} ,y)\in[0,\tau]\times\R^2\times [0,3/\s]$,  $\chi_{\sigma}(y)=\chi(\sigma y)$ with $0<\x, \sigma\ll 1$, and  $\tilde{\mathfrak{f}}^{\sigma}_i:=\chi_\s(y)\tilde{\mathfrak{f}}_i,\  \tilde{\mathfrak{g}}^{\sigma}:=\chi_\s(y)\tilde{\mathfrak{g}}$. Here $\chi(\cdot)$ is the one defined in \eqref{5.7}.  We impose \eqref{5.66} with the following boundary conditions
\begin{align}\label{5.67}
\begin{cases}
\dis\pa_y \Psi_i(t,x_{\sp} ,y)\big|_{y=0}=0,\quad \pa_y\Theta(t,x_{\sp} ,y)\big|_{y=0}=0,\\[2mm]
\dis\Psi(t,x_{\sp} ,y)\big|_{y=\frac3\s}=0,\quad  \Theta(t,x_{\sp} ,y)\big|_{y=\frac3\s}=0.
\end{cases}
\end{align}
We impose \eqref{5.66} with the following cut-off initial data
\begin{align*}%\label{5.67-1}
\Psi(0,x_{\sp} ,y)=(u_0-u_b)\chi_{\sigma}(y),\qquad \Theta(0,x_{\sp} ,y)=(\theta_0-\theta_a)\chi_{\sigma}(y),
\end{align*}
then the compatibility condition of initial data  at $y=\frac3{\sigma}$ is also satisfied due to the property of $\chi(s)$.
For the approximate problem \eqref{5.66}-\eqref{5.67}, now we can use the standard linear parabolic theory to obtain  the existence of  smooth solution in Sobolev space provided the initial data and $(\rho^0,\fu^0,T^0)$ are suitably smooth.  To prove the lemma, we need only to obtain some uniform estimates of $(\Psi,\Theta)$ with respect to $\sigma$ and $\xi$,  then take the limit $\sigma, \x\rightarrow0+$.

\vspace{2mm}

\noindent{\it Step 2. Uniform energy estimates.}  We use induction argument to prove the uniform estimates. Firstly we consider the zero-order derivatives estimation.  Multiplying $\eqref{5.66}_1$ by $(1+y)^{l_0}\Psi_i$ and  integrating the resultant equation over $[0,t]\times\R^2\times[0,\frac{3}{\sigma}]$, we have
\begin{align}\label{5.77}
& \iint \frac12(1+y)^{l_0}|\Psi_i(t,x_\sp,y)|^2 dx_{\sp}  dy+ \frac12 \int_0^t\iint  \pa_3\fu_3^0 y (1+y)^{l_0}\chi_{\sigma}(y) \pa_y(|\Psi_i|^2) dx_{\sp}  dyds\nonumber\\
&\leq \int_0^t\iint \tilde{\mu}\, (1+y)^{{l_0}}\partial_{yy} \Psi_{i} \Psi_i dx_{\sp}  dyds+\xi \int_0^t\iint (1+y)^{l_0}\Delta_{\sp} \Psi\cdot \Psi_i dx_{\sp}  dyds\nonumber\\
&\,+C\|(\rho,\fu,T)\|_{W^{1,\infty}}\int_0^t \|(\Psi, \Theta)\|_{L^2_{l_0}}^2ds+C\int_0^t\|\tilde{\mathfrak{f}}_i\|_{L^2_{l_0}}^2ds+ C  \|\Psi_i(0)\|_{L^2_{l_0}}^2.
\end{align}

It is easy to know that
\begin{align}\label{5.70}
\begin{cases}
\chi_{\sigma}(y)\equiv0,\  \forall\, y\geq \frac2\sigma,\\[1.5mm]
 |y\pa_y\chi_{\sigma}(y)|=|y\sigma \chi'(\sigma y)|\leq C, \ \forall\, y\in\R_+.
\end{cases}
\end{align}
For the second term on LHS of \eqref{5.77}, integrating by part w.r.t. $y$ and using \eqref{5.70}, we obtain
\begin{align}\label{5.78}
\left|\int_0^t\iint \pa_3\fu_3^0 y (1+y)^{l_0} \chi_{\sigma}(y) \pa_y(|\Psi_i|^2) dx_{\sp}  dyds\right|\leq C \|\fu\|_{W^{1,\infty}} \int_0^t\|\Psi_i(s)\|_{L^2_{l_0}}^2ds.
%&\leq C\int_0^{\frac{3}{\sigma}}\int_{\R^2} (1+y)^l |\Psi_i|^2dx_{\sp}  dy
\end{align}
For the viscous terms, integrating by parts to yield
\begin{align*}%\label{5.79}
\iint (1+y)^{l_0}\Delta_{\sp} \Psi_i\cdot \Psi_i dx_{\sp}  dy=- \|\nabla_{\sp} \Psi_i\|^2_{L^2_{l_0}},
\end{align*}
and
\begin{align}\label{5.80}
&\iint \tilde{\mu}\,(1+y)^{l_0} \partial_{yy} \Psi_{i} \Psi_i dx_{\sp}  dy\nonumber\\
&=-\iint \tilde{\mu}\, (1+y)^{l_0} |\partial_{y} \Psi_{i}|^2 dx_{\sp}  dy-\iint l_0\tilde{\mu}\, (1+y)^{l_0-1} \partial_{y} \Psi_{i} \Psi_i dx_{\sp}  dy
\nonumber\\
&\quad+\int_{\R^2} \tilde{\mu}\, (1+y)^{l_0}  \partial_{y} \Psi_{i} \Psi_i dx_{\sp} \Big|_{y=0}^{y=\f3\sigma}\nonumber\\
&\leq -\frac12\iint \tilde{\mu}\, (1+y)^{l_0} |\partial_{y} \Psi_{i}|^2 dx_{\sp}  dy+C\|\Psi_i\|^2_{L^2_{l_0}},
\end{align}
where we have used \eqref{5.67} in \eqref{5.80}.

Substituting \eqref{5.78}-\eqref{5.80} into \eqref{5.77}, we get that, for some positive constant $c_0>0$
\begin{align}\label{5.81}
& \|\Psi(t)\|_{L^2_{l_0}}^2+\int_0^t\|\pa_y\Psi(s)\|_{L^2_{l_0}}^2+\x \|\nabla_{\sp} \Psi(s)\|_{L^2_{l_0}}^2 ds\nonumber\\
&\leq C\Big(\|\Psi(0)\|_{L^2_{l_0}}^2+\int_0^t\|\tilde{\mathfrak{f}}\|_{L^2_{l_0}}^2ds\Big)+C(\|(\rho,\fu,T)\|_{W^{1,\infty}}) \int_0^t\|(\Psi,\Theta)(s) \|_{L^2_{l_0}}^2ds.
\end{align}
Similarly we can prove
\begin{align*}%\label{5.82}
& \|\Theta(t)\|_{L^2_{l_0}}^2+\int_0^t\|\pa_y\Theta(s)\|_{L^2_{l_0}}^2+\x \|\nabla_{\sp} \Theta(s)\|_{L^2_{l_0}}^2ds\nonumber\\
&\leq C\Big(\|\Theta(0)\|_{L^2_{l_0}}^2+\int_0^t\|\tilde{\mathfrak{g}}\|_{L^2_{l_0}}^2ds\Big)+C(\|(\rho,\fu,T)\|_{W^{1,\infty}}) \int_0^t\|\Theta(s)\|_{L^2_{l_0}}^2ds,
\end{align*}
which, together with \eqref{5.81}, yields
\begin{align}\label{5.83}
& \|(\Psi,\Theta)(t)\|_{L^2_{l_0}}^2+\int_0^t\| \pa_y(\Psi,\Theta)(s)\|_{L^2_{l_0}}^2+\x \|\nabla_{\sp} (\Psi,\Theta)(s)\|_{L^2_{l_0}}^2ds\nonumber\\
&\leq C\Big(\|(\Psi,\Theta)(0)\|_{L^2_{l_0}}^2+\int_0^t\|(\tilde{\mathfrak{f}},\tilde{\mathfrak{g}})(s)\|_{L^2_{l_0}}^2ds\Big)+C(\|(\rho,\fu,T)\|_{W^{1,\infty}}) \int_0^t\|(\Psi,\Theta)(s)\|_{L^2_{l_0}}^2ds.
\end{align}
Now applying the Gronwall's inequality to \eqref{5.83}, we have
\begin{align}\label{5.87}
&\|(\Psi,\Theta)(t)\|_{L^2_{l_0}}^2+\int_0^t\| \pa_y(\Psi,\Theta)(s)\|_{L^2_{l_0}}^2+\x \|\nabla_{\sp} (\Psi,\Theta)(s)\|_{L^2_{l_0}}^2ds\nonumber\\
&\leq C(t,\|(\rho,\fu,T)\|_{W^{1,\infty}})\Big\{\|(\Psi,\Theta)(0)\|_{L^2_{l_0}}^2 +\int_0^t\|(\tilde{\mathfrak{f}},\tilde{\mathfrak{g}})(s)\|_{L^2_{l_0}}^2ds\Big\}.
\end{align}
%\textcolor{blue}{(We must consider the weighted estimate because it is necessary when we work on derivatives estimates. )}

\vspace{2mm}

We  shall use induction arguments to close the uniform energy estimates. We assume, for $0\leq r\leq k-1$ ($r\geq1$), that
\begin{align}\label{5.88}
&\sum_{j=0}^{r-1}\sum_{2\alpha+|\beta|=j}  \Big\{\|\partial_t^{\alpha}\nabla_{\bar{x}}^{\beta}(\Psi,\Theta)(t)\|_{L^2_{l_j}}^2 +\int_0^t \|\partial_t^{\alpha}\nabla_{\bar{x}}^{\beta}\pa_y(\Psi,\Theta)(t)\|_{L^2_{l_j}}^2+
\xi \|\nabla_{\sp}\partial_t^{\alpha}\nabla_{\bar{x}}^{\beta}(\Psi,\Theta)\|^2_{L^2_{l_j}} ds\Big\}\nonumber\\
&\leq C(t,\|(\rho,\fu,T)\|_{W^{r,\infty}}) \Bigg\{\sum_{j=0}^{r-1}\sum_{2\alpha+|\beta|=j}  \Big[ \|\partial_t^{\alpha}\nabla_{\bar{x}}^{\beta}(\Psi,\Theta)(0)\|_{L^2_{l_j}}^2+\int_0^t \|\partial_t^{\alpha}\nabla_{\bar{x}}^{\beta}(\tilde{\mathfrak{f}},\tilde{\mathfrak{g}})(s)\|_{L^2_{l_j}}^2ds\Big]\nonumber\\
&\qquad+\sum_{j=0}^{r-3}\sum_{\beta+2\gamma=j}\|\pa_t^{\gamma}\nabla_{\bar{x}}^{\beta}  (\tilde{\mathfrak{f}},\tilde{\mathfrak{g}})(0)\|^2_{L^2_{l_j}}\Bigg\}.
\end{align}
Here we point out that one order of time derivative $\partial_t$ is equal to two orders of space derivatives.

\smallskip

Now we consider the $r$-order derivative estimates. Let  $\pa_{\sp}^{\beta}=\pa_{x_1}^{\b_1}\pa_{x_2}^{\beta_2}$. Applying $\pa_t^{\a}\pa_{\sp}^{\b} $ to $\eqref{5.66}_1$, we have
\begin{align}\label{5.91}
&\pa_t \pa_t^{\a}\pa_{\sp}^{\b}  \Psi_i+\fu_{\sp} ^0\cdot \nabla_{\sp} \pa_t^{\a}\pa_{\sp}^{\b} \Psi_i+\pa_3\fu_3^0\cdot y\chi_{\sigma}(y) \pa_y \pa_t^{\a}\pa_{\sp}^{\b} \Psi_i\nonumber\\
&=\pa_t^{\a}\pa_{\sp}^{\b} \big(\tilde{\mu} \pa_y^2 \Psi_i \big)+\xi \Delta_{\sp}  \pa_t^{\a}\pa_{\sp}^{\b}  \Psi_i-[\pa_t^{\a}\pa_{\sp}^{\b}, \, \fu_{\sp} ^0\cdot \nabla_{\sp}] \Psi_i
-y\chi_{\sigma}(y) [\pa_t^{\a}\pa_{\sp}^{\b} , \pa_3\fu_3^0  \pa_y] \Psi_i\nonumber\\
&\qquad-\pa_t^{\a}\pa_{\sp}^{\b}  \Big\{\Psi\cdot \nabla_{\sp} \fu_i^0-\frac{\pa_i p^0}{3p^0} \Theta\Big\}+\pa_t^{\a}\pa_{\sp}^{\b}   \tilde{\mathfrak{f}}_i^\sigma.
\end{align}

Let $2\alpha+|\beta|=r$.  Multiplying \eqref{5.91} by  $(1+y)^{l_r} \pa_t^{\a}\pa_{\sp}^{\b}  \Psi_i$, and integrating the resultant equation over $[0,t]\times\R^2\times[0,\frac{3}{\sigma}]$, we have
\begin{align}\label{5.92}
&\frac12 \|\pa_t^{\a}\pa_{\sp}^{\b} \Psi_i(t)\|^2_{L^2_{l_r}}+\xi \int_0^t  \|\nabla_{\sp}\pa_t^{\a}\pa_{\sp}^{\b} \Psi_i(s)\|^2_{L^2_{l_r}} ds\nonumber\\
&\leq \frac12\|\pa_t^{\a}\pa_{\sp}^{\b} \Psi_i(0)\|^2_{L^2_{l_r}}
+\int_0^t\iint (1+y)^{l_r} \pa_t^{\a}\pa_{\sp}^{\b}  \Psi_i \cdot \pa_t^{\a}\pa_{\sp}^{\b} \big(\tilde{\mu} \pa_y^2 \Psi_i \big) dx_{\sp}  dy ds
\nonumber\\
&\quad
+\Big|\int_0^t\iint  y\chi_{\sigma}(y) (1+y)^{l_r} [\pa_t^{\a}\pa_{\sp}^{\b} , \pa_3\fu_3^0  \pa_y] \Psi_i\cdot \pa_t^{\a}\pa_{\sp}^{\b}  \Psi_i  dx_{\sp}  dyds\Big|
\nonumber\\
&\quad+C(t,\|(\rho,\fu,T)\|_{W^{r+1,\infty}})\sum_{2\tilde\alpha+|\tilde\beta|\leq r}  \int_0^t \|\pa_t^{\tilde\a}\nabla_{\bar{x}}^{\tilde\b} (\Psi,\Theta)(s)\|^2_{L^2_{l_r}}ds\nonumber\\
&\quad+C\sum_{2\tilde\alpha+|\tilde\beta|\leq r}\int_0^t \| \pa_t^{\tilde\a}\pa_{\sp}^{\tilde\b}   \tilde{\mathfrak{f}}_i(s)\|^2_{L^2_{l_r}}ds.
\end{align}

Using \eqref{5.67}, we have that
\begin{align*}
\pa_t^{\a}\pa_{\sp}^{\b}  \Psi_i \Big|_{y=\frac{3}{\sigma}}=0\quad\mbox{and} \quad \pa_t^{\a}\pa_{\sp}^{\b}  \partial_y\Psi_i \Big|_{y=0}=0,
\end{align*}
which, together with  integrating by parts, yields  that
\begin{align}\label{5.97}
&\int_0^t\iint (1+y)^{l_r} \pa_t^{\a}\pa_{\sp}^{\b}  \Psi_i \cdot \pa_t^{\a}\pa_{\sp}^{\b} \big(\tilde{\mu} \pa_y^2 \Psi_i \big) dx_{\sp}  dy ds\nonumber\\
&=-\int_0^t\iint (1+y)^{l_r} \pa_t^{\a}\pa_{\sp}^{\b} \pa_y \Psi_i \cdot \pa_t^{\a}\pa_{\sp}^{\b} \big(\tilde{\mu} \pa_y \Psi_i \big) dx_{\sp}  dy ds\nonumber\\
&\quad-l_r\int_0^t\iint (1+y)^{l_r-1} \pa_t^{\a}\pa_{\sp}^{\b}  \Psi_i \cdot \pa_t^{\a}\pa_{\sp}^{\b} \big(\tilde{\mu} \pa_y \Psi_i \big) dx_{\sp}  dy ds\nonumber\\
&\leq -\frac12\int_0^t\iint \tilde{\mu}\, (1+y)^{l_r} |\pa_t^{\a}\pa_{\sp}^{\b}\pa_y\Psi_i|^2 dx_{\sp}  dy ds\nonumber\\
&\quad+C(t,\|(\rho,\fu,T)\|_{W^{r+1,\infty}}) \sum_{j=0}^{r} \sum_{2\tilde\alpha+|\tilde\beta|=j}  \int_0^t \|\pa_t^{\tilde\a}\nabla_{\bar{x}}^{\tilde\b} \Psi(s)\|^2_{L^2_{l_j}}ds.
\end{align}

The third term on RHS of \eqref{5.92} is bounded by
\begin{align}\label{5.96}
&\Big|\int_0^t\iint  y\chi_{\sigma}(y) (1+y)^{l_r} [\pa_t^{\a}\pa_{\sp}^{\b} , \pa_3\fu_3^0  \pa_y] \Psi_i\cdot \pa_t^{\a}\pa_{\sp}^{\b}  \Psi_i  dx_{\sp}  dyds\Big|\nonumber\\
&\leq
C(t,\|(\rho,\fu,T)\|_{W^{r+1,\infty}})\sum_{2\tilde\alpha+|\tilde\beta|\leq r}  \int_0^t \|\pa_t^{\tilde\a}\nabla_{\bar{x}}^{\tilde\b} \Psi(s)\|^2_{L^2_{l_r}}ds\nonumber\\
&\quad+\sum_{j=0}^{r-1}  \int_0^t \sum_{2\tilde\alpha+|\tilde\beta|=j} \|\pa_y \pa_t^{\tilde\a}\nabla_{\bar{x}}^{\tilde\b} \Psi(s)\|^2_{L^2_{l_j}}ds,
\end{align}
where we have used the fact $l_r+2=l_{r-1}\leq l_j$ for $0\leq j\leq r-1$.

Combining \eqref{5.92}-\eqref{5.96}, and taking $\lambda>0$ suitably small,  we obtain
\begin{align}\label{5.98}
&\sum_{2\alpha+|\beta|=r} \Bigg\{\|\pa_t^{\a}\pa_{\sp}^{\b} \Psi(t)\|^2_{L^2_{l_r}}+ \int_0^t \|\pa_y\pa_t^{\a}\pa_{\sp}^{\b} \Psi(s)\|^2_{L^2_{l_r}}+ \xi\|\nabla_{\sp}\pa_t^{\a}\pa_{\sp}^{\b} \Psi(s)\|^2_{L^2_{l_r}} ds\Bigg\} \nonumber\\
&\leq C \Bigg\{\sum_{2\alpha+|\beta|=r}\|\pa_t^{\a}\pa_{\sp}^{\b} \Psi(0)\|^2_{L^2_{l_r}}
+\sum_{j=0}^{r} \sum_{2\tilde\alpha+|\tilde\beta|=j} \int_0^t \| \pa_t^{\tilde\a}\pa_{\sp}^{\tilde\b}   \tilde{\mathfrak{f}}(s)\|^2_{L^2_{l_j}}ds\Bigg\}\nonumber\\
&\quad+C(t,\|(\rho,\fu,T)\|_{W^{r+1,\infty}}) \sum_{j=0}^{r} \sum_{2\tilde\alpha+|\tilde\beta|=j} \int_0^t \|\pa_t^{\tilde\a}\nabla_{\bar{x}}^{\tilde\b} (\Psi,\Theta)(s)\|^2_{L^2_{l_j}}ds\nonumber\\
&\quad +C\sum_{j=0}^{r-1} \sum_{2\tilde\alpha+|\tilde\beta|=j} \int_0^t  \|\pa_y \pa_t^{\tilde\a}\nabla_{\bar{x}}^{\tilde\b} \Psi(s)\|^2_{L^2_{l_j}}ds.
\end{align}

\vspace{2mm}

For the normal derivative estimate, applying $\pa_y$ to $\eqref{5.91}$, we have
\begin{align}\label{5.99}
&\pa_t \pa_t^{\a}\pa_{\sp}^{\b} \pa_y \Psi_i+\fu_{\sp} ^0\cdot \nabla_{\sp} \pa_t^{\a}\pa_{\sp}^{\b} \pa_y\Psi_i+\pa_3\fu_3^0\cdot y\chi_{\sigma}(y)   \pa_t^{\a}\pa_{\sp}^{\b} \pa_y^2\Psi_i\nonumber\\
&=\pa_t^{\a}\pa_{\sp}^{\b} \big(\tilde{\mu} \pa_y^3 \Psi_i \big)+\xi \Delta_{\sp}  \pa_t^{\a}\pa_{\sp}^{\b} \pa_y \Psi_i-y\chi_{\sigma}(y) [\pa_t^{\a}\pa_{\sp}^{\b} , \pa_3\fu_3^0  \pa_y^2]\Psi_i\nonumber\\
&\quad-\pa_3\fu_3^0\cdot \pa_y\big(y\chi_{\sigma}(y)\big)  \pa_t^{\a}\pa_{\sp}^{\b} \pa_y\Psi_i-[\pa_t^{\a}\pa_{\sp}^{\b}, \, \fu_{\sp} ^0\cdot \nabla_{\sp}] \pa_y\Psi_i
\nonumber\\
&\quad-\pa_y\big(y\chi_{\sigma}(y)\big) [\pa_t^{\a}\pa_{\sp}^{\b} , \pa_3\fu_3^0  \pa_y] \Psi_i-\pa_t^{\a}\pa_{\sp}^{\b} \pa_y \Big\{\Psi\cdot \nabla_{\sp} \fu_i^0-\frac{\pa_i p^0}{3p^0} \Theta\Big\}+\pa_t^{\a}\pa_{\sp}^{\b}  \pa_y \tilde{\mathfrak{f}}_i^\sigma.
\end{align}
Multiplying \eqref{5.99} by  $(1+y)^{l_r} \pa_t^{\a}\pa_{\sp}^{\b} \pa_y\Psi_i$, and integrating the resultant equation over $[0,t]\times\R^2\times [0,\frac{3}{\sigma}]$, we have
\begin{align}\label{5.100}
&\sum_{2\alpha+|\beta|=r-1} \Bigg\{\frac12\|\pa_t^{\a}\pa_{\sp}^{\b} \pa_y\Psi(t)\|^2_{L^2_{l_r}}+ \int_0^t %\|\pa_y\pa_t^{\a}\pa_{\sp}^{\b} \pa_y \Psi(s)\|^2_{L^2_{l_r}}+
\xi\|\nabla_{\sp}\pa_t^{\a}\pa_{\sp}^{\b} \pa_y\Psi(s)\|^2_{L^2_{l_r}} ds\Bigg\} \nonumber\\
&\leq \sum_{2\alpha+|\beta|=r-1}  \int_0^t\iint (1+y)^{l_r} \pa_t^{\a}\pa_{\sp}^{\b}  \pa_y\Psi_i \cdot \pa_t^{\a}\pa_{\sp}^{\b} \big(\tilde{\mu} \pa_y^3\Psi_i \big) dx_{\sp}  dy ds
\nonumber\\
&\quad+\Bigg\{\sum_{2\alpha+|\beta|=r-1} \frac12\|\pa_t^{\a}\pa_{\sp}^{\b} \pa_y\Psi(0)\|^2_{L^2_{l_r}}
+C\sum_{j=0}^{r-1} \sum_{2\tilde\alpha+|\tilde\beta|=j} \int_0^t \| \pa_t^{\tilde\a}\pa_{\sp}^{\tilde\b}   \pa_y\tilde{\mathfrak{f}}(s)\|^2_{L^2_{l_j}}ds\Bigg\}\nonumber\\
&\quad+C(t,\|(\rho,\fu,T)\|_{W^{r,\infty}}) \sum_{j=0}^{r} \sum_{2\tilde\alpha+|\tilde\beta|=j} \int_0^t \|\pa_t^{\tilde\a}\nabla_{\bar{x}}^{\tilde\b} (\Psi,\Theta)(s)\|^2_{L^2_{l_j}}ds\nonumber\\
&\quad +C\sum_{j=0}^{r-1} \sum_{2\tilde\alpha+|\tilde\beta|=j} \int_0^t  \|\pa_y \pa_t^{\tilde\a}\nabla_{\bar{x}}^{\tilde\b} \Psi(s)\|^2_{L^2_{l_j}}ds.
\end{align}

Using \eqref{5.67} and \eqref{5.91}, we have
\begin{align}\nonumber
\begin{cases}
\pa_t^{\a}\pa_{\sp}^{\b}  \pa_y\Psi_i \Big|_{y=0}=0,\\[3.5mm]
\pa_t^{\a}\pa_{\sp}^{\b} \big(\tilde{\mu} \pa_y^2\Psi_i \big) \Big|_{y=\f3\sigma} =-\xi \Delta_{\sp} \pa_t^{\a}\pa_{\sp}^{\b}\Psi_i \Big|_{y=\f3\sigma} =0,
\end{cases}
\end{align}
which, together with integrating by parts, yields that
\begin{align}\label{5.102}
&\int_0^t\iint (1+y)^{l_r} \pa_t^{\a}\pa_{\sp}^{\b}  \pa_y\Psi_i \cdot \pa_t^{\a}\pa_{\sp}^{\b} \big(\tilde{\mu} \pa_y^3\Psi_i \big) dx_{\sp}  dy ds
\nonumber\\
&=-\int_0^t\iint (1+y)^{l_r} \pa_t^{\a}\pa_{\sp}^{\b}  \pa_y^2\Psi_i \cdot \pa_t^{\a}\pa_{\sp}^{\b} \big(\tilde{\mu} \pa_y^2\Psi_i \big) dx_{\sp}  dy ds
\nonumber\\
&\quad-l_r \int_0^t\iint (1+y)^{l_r-1} \pa_t^{\a}\pa_{\sp}^{\b}  \pa_y\Psi_i \cdot \pa_t^{\a}\pa_{\sp}^{\b} \big(\tilde{\mu} \pa_y^2\Psi_i \big) dx_{\sp}  dy ds
\nonumber\\
&\leq -\frac12\int_0^t\iint \tilde{\mu}\, (1+y)^{l_r} |\pa_t^{\a}\pa_{\sp}^{\b}  \pa_y^2\Psi_i|^2 dx_{\sp}  dy ds\nonumber\\
&\quad+C(t,\|(\rho,\fu,T)\|_{W^{r+1,\infty}}) \sum_{j=0}^{r} \sum_{2\tilde\alpha+|\tilde\beta|=j}  \int_0^t \|\pa_t^{\tilde\a}\nabla_{\bar{x}}^{\tilde\b} \Psi(s)\|^2_{L^2_{l_j}}ds.
\end{align}
Combining \eqref{5.100}-\eqref{5.102}, we obtain
\begin{align}\label{5.103}
&\sum_{2\alpha+|\beta|=r-1} \Bigg\{\frac12\|\pa_t^{\a}\pa_{\sp}^{\b} \pa_y\Psi(t)\|^2_{L^2_{l_r}}+ \int_0^t \|\pa_y^2\pa_t^{\a}\pa_{\sp}^{\b}  \Psi(s)\|^2_{L^2_{l_r}}+
\xi\|\nabla_{\sp}\pa_t^{\a}\pa_{\sp}^{\b} \pa_y\Psi(s)\|^2_{L^2_{l_r}} ds\Bigg\} \nonumber\\
&\leq C\Bigg\{\sum_{2\alpha+|\beta|=r-1} \frac12\|\pa_t^{\a}\pa_{\sp}^{\b} \pa_y\Psi(0)\|^2_{L^2_{l_r}}
+C\sum_{j=0}^{r-1} \sum_{2\tilde\alpha+|\tilde\beta|=j} \int_0^t \| \pa_t^{\tilde\a}\pa_{\sp}^{\tilde\b}   \pa_y\tilde{\mathfrak{f}}(s)\|^2_{L^2_{l_j}}ds\Bigg\}\nonumber\\
&\quad+C(t,\|(\rho,\fu,T)\|_{W^{r,\infty}}) \sum_{j=0}^{r} \sum_{2\tilde\alpha+|\tilde\beta|=j} \int_0^t \|\pa_t^{\tilde\a}\nabla_{\bar{x}}^{\tilde\b} (\Psi,\Theta)(s)\|^2_{L^2_{l_j}}ds\nonumber\\
&\quad +C\sum_{j=0}^{r-1} \sum_{2\tilde\alpha+|\tilde\beta|=j} \int_0^t  \|\pa_y \pa_t^{\tilde\a}\nabla_{\bar{x}}^{\tilde\b} \Psi(s)\|^2_{L^2_{l_j}}ds.
\end{align}

For the higher order normal derivatives estimate, we shall use the  equation directly but not the energy method. In fact, applying $\pa_t^{\a}\pa_{\sp}^{\b}\pa_y^n$ with $n\geq0$ to $\eqref{5.66}_1$, then we have
\begin{align*}
\tilde{\mu} \pa_t^{\a}\pa_{\sp}^{\b}\pa_y^{n+2} \Psi_i&=\pa_t \pa_t^{\a}\pa_{\sp}^{\b}  \pa_y^{n} \Psi_i-\xi \Delta_{\sp}  \pa_t^{\a}\pa_{\sp}^{\b}  \pa_y^{n}\Psi_i
-[ \pa_t^{\a}\pa_{\sp}^{\b} ,\tilde{\mu}] \pa_{y}^{n+2}\Psi_i\nonumber\\
&+\pa_t^{\a}\pa_{\sp}^{\b}\Big\{\pa_3\fu_3^0\cdot y\chi_{\sigma}(y) \pa_y^{n+1} \Psi_i\Big\}
+\pa_t^{\a}\pa_{\sp}^{\b}\Big\{ \pa_3\fu_3^0\,  [\pa_y^n, \, y\chi_{\sigma}(y)]\pa_y \Psi_i\Big\}\nonumber\\
&+\pa_t^{\a}\pa_{\sp}^{\b}  \Big\{ \fu_{\sp} ^0\cdot \nabla_{\sp} \pa_y^{n}\Psi_i+\pa_y^{n}\Psi\cdot \nabla_{\sp} \fu_i^0-\frac{\pa_i p^0}{3p^0} \pa_y^{n}\Theta\Big\}-\pa_t^{\a}\pa_{\sp}^{\b}\pa_y^{n}\tilde{\mathfrak{f}}_i^\sigma,
\end{align*}
which yields that
\begin{align}\label{5.104}
&\sum_{2\alpha+|\beta|=m} | \pa_t^{\a}\pa_{\sp}^{\b}\pa_y^{n+2} \Psi|\nonumber\\
&\leq C\sum_{2\tilde\alpha+|\tilde\beta|=m+2} | \pa_t^{\tilde\a}\pa_{\sp}^{\tilde\b}\pa_y^{n} \Psi|
+C(\|(\rho,\fu,T)\|_{W^{r,\infty}}) \Big\{\sum_{2\tilde\alpha+|\tilde\beta|\leq m+n+1}   |\pa_t^{\tilde\a}\nabla_{\bar{x}}^{\tilde\b} (\Psi,\Theta)|\nonumber\\
&\qquad +\sum_{2\tilde\alpha+|\tilde\beta|\leq m+n} y |\pa_y \pa_t^{\tilde\a}\nabla_{\bar{x}}^{\tilde\b} \Psi|\Big\}+\sum_{2\tilde\alpha+|\tilde\beta|\leq m+n} |\pa_t^{\tilde\a}\nabla_{\bar{x}}^{\tilde\b}\tilde{\mathfrak{f}}|.
\end{align}

For any fixed $0\leq n\leq r-2$, using \eqref{5.104}, a direct calculation shows that
\begin{align}\label{5.105}
&\sum_{2\alpha+|\beta|=r-n-2} \Bigg\{ \| \pa_y^{n+2}\pa_t^{\a}\pa_{\sp}^{\b} \Psi(t)\|^2_{L^2_{l_r}}+\int_0^t \|\pa_y^{n+3} \pa_t^{\a}\pa_{\sp}^{\b} \Psi(s)\|^2_{L^2_{l_r}}ds\Bigg\}\nonumber\\
&\leq C\sum_{2\tilde\alpha+|\tilde\beta|=r-n} \Bigg\{ \|\pa_y^{n}\pa_t^{\tilde\a}\pa_{\sp}^{\tilde\b} \Psi(t)\|^2_{L^2_{l_r}}+\int_0^t\|\pa_y^{n+1}\pa_t^{\tilde\a}\pa_{\sp}^{\tilde\b} \Psi(s)\|^2_{L^2_{l_r}}ds\Bigg\}\nonumber\\
&\, +C(\|(\rho,\fu,T)\|_{W^{r,\infty}}) \sum_{j=0}^{r-1}\sum_{2\tilde\alpha+|\tilde\beta|=j}   \Bigg\{\|\pa_t^{\tilde\a}\nabla_{\bar{x}}^{\tilde\b} (\Psi,\Theta)(t)\|^2_{L^2_{l_j}}+\int_0^t \|\pa_y\pa_t^{\tilde\a}\nabla_{\bar{x}}^{\tilde\b}\Psi(s)\|^2_{L^2_{l_j}} ds\Bigg\}\nonumber\\
&\,\,
+C \sum_{j=0}^{r-2}\sum_{2\tilde\alpha+|\tilde\beta|=j} \|\pa_t^{\tilde\a}\nabla_{\bar{x}}^{\tilde\b}\tilde{\mathfrak{f}}(t)\|^2_{L^2_{l_j}}
+ C\sum_{j=0}^{r-1}\sum_{2\tilde\alpha+|\tilde\beta|=j}\int_0^t\|\pa_t^{\tilde\a}\nabla_{\bar{x}}^{\tilde\b}\tilde{\mathfrak{f}}(s)\|^2_{L^2_{l_j}} ds.
\end{align}
Noting \eqref{5.98} and  \eqref{5.103}, by using induction arguments on $0\leq n\leq r-2$ in \eqref{5.105}, we can obtain
\begin{align}\label{5.106}
&\sum_{2\alpha+|\beta|=r} \Bigg\{ \| \pa_t^{\a}\nabla_{\bar{x}}^{\b} \Psi(t)\|^2_{L^2_{l_r}}+\int_0^t \|\pa_y \pa_t^{\a}\nabla_{\bar{x}}^{\b} \Psi(s)\|^2_{L^2_{l_r}}
+\xi\|\nabla_{\sp} \pa_t^{\a}\nabla_{\bar{x}}^{\b} \Psi(s)\|^2_{L^2_{l_r}}ds\Bigg\}\nonumber\\
&\leq
C(\|(\rho,\fu,T)\|_{W^{r,\infty}}) \sum_{j=0}^{r-1}\sum_{2\tilde\alpha+|\tilde\beta|=j}   \Bigg\{\|\pa_t^{\tilde\a}\nabla_{\bar{x}}^{\tilde\b} (\Psi,\Theta)(t)\|^2_{L^2_{l_j}}+\int_0^t \|\pa_y\pa_t^{\tilde\a}\nabla_{\bar{x}}^{\tilde\b}\Psi(s)\|^2_{L^2_{l_j}} ds\Bigg\}\nonumber\\
&\quad+C(t,\|(\rho,\fu,T)\|_{W^{r+1,\infty}}) \sum_{j=0}^{r} \sum_{2\tilde\alpha+|\tilde\beta|=j} \int_0^t \|\pa_t^{\tilde\a}\nabla_{\bar{x}}^{\tilde\b} (\Psi,\Theta)(s)\|^2_{L^2_{l_j}}ds\nonumber\\
&\quad
+C \sum_{j=0}^{r-2}\sum_{2\tilde\alpha+|\tilde\beta|=j} \|\pa_t^{\tilde\a}\nabla_{\bar{x}}^{\tilde\b}\tilde{\mathfrak{f}}(t)\|^2_{L^2_{l_j}}
+ C\sum_{j=0}^{r}\sum_{2\tilde\alpha+|\tilde\beta|=j}\int_0^t\|\pa_t^{\tilde\a}\nabla_{\bar{x}}^{\tilde\b}\tilde{\mathfrak{f}}(s)\|^2_{L^2_{l_j}} ds\nonumber\\
&\quad+C \sum_{2\alpha+|\beta|=r}\|\pa_t^{\a}\nabla_{\bar{x}}^{\b} \Psi(0)\|^2_{L^2_{l_r}}\nonumber\\
&\leq C(t,\|(\rho,\fu,T)\|_{W^{r+1,\infty}}) \Bigg\{\sum_{j=0}^{r} \sum_{2\tilde\alpha+|\tilde\beta|=j} \Big[\|\pa_t^{\a}\nabla_{\bar{x}}^{\b} \Psi(0)\|^2_{L^2_{l_j}}+\int_0^t\|\pa_t^{\tilde\a}\nabla_{\bar{x}}^{\tilde\b}(\tilde{\mathfrak{f}},\tilde{\mathfrak{g}})(s)\|^2_{L^2_{l_j}} ds\Big]\nonumber\\
&\,\,
+\int_0^t \sum_{2\alpha+|\beta|=r}\|\pa_t^{\a}\nabla_{\bar{x}}^{\b} (\Psi,\Theta)(s)\|^2_{L^2_{l_j}}ds
+C \sum_{j=0}^{r-2}\sum_{2\tilde\alpha+|\tilde\beta|=j}  \|\pa_t^{\tilde\a}\nabla_{\bar{x}}^{\tilde\b}(\tilde{\mathfrak{f}},\tilde{\mathfrak{g}})(0)\|^2_{L^2_{l_j}}\Bigg\},
\end{align}
where we have used \eqref{5.88} in the second inequality.

\smallskip

By similar argument as \eqref{5.106}, we can also get
\begin{align}\label{5.107}
&\sum_{2\alpha+|\beta|=r} \Bigg\{ \| \pa_t^{\a}\nabla_{\bar{x}}^{\b} \Theta(t)\|^2_{L^2_{l_r}}+\int_0^t \|\pa_y \pa_t^{\a}\nabla_{\bar{x}}^{\b} \Theta(s)\|^2_{L^2_{l_r}}+\xi\|\nabla_{\sp} \pa_t^{\a}\nabla_{\bar{x}}^{\b} \Theta(s)\|^2_{L^2_{l_r}}ds\Bigg\}\nonumber\\
&\leq C(t,\|(\rho,\fu,T)\|_{W^{r+1,\infty}}) \Bigg\{\sum_{j=0}^{r} \sum_{2\tilde\alpha+|\tilde\beta|=j} \Big[\|\pa_t^{\a}\nabla_{\bar{x}}^{\b} \Theta(0)\|^2_{L^2_{l_j}}+\int_0^t\|\pa_t^{\tilde\a}\nabla_{\bar{x}}^{\tilde\b}(\tilde{\mathfrak{f}},\tilde{\mathfrak{g}})(s)\|^2_{L^2_{l_j}} ds\Big]\nonumber\\
&\,\,
+\int_0^t \sum_{2\alpha+|\beta|=r}\|\pa_t^{\a}\nabla_{\bar{x}}^{\b} \Theta(s)\|^2_{L^2_{l_j}}ds
+C \sum_{j=0}^{r-2}\sum_{2\tilde\alpha+|\tilde\beta|=j} \|\pa_t^{\tilde\a}\nabla_{\bar{x}}^{\tilde\b}(\tilde{\mathfrak{f}},\tilde{\mathfrak{g}})(0)\|^2_{L^2_{l_j}}\Bigg\}
\end{align}
Now combining \eqref{5.106} and \eqref{5.107}, then using the Gronwall's inequality, we obtain
\begin{align}
&\sum_{2\alpha+|\beta|=r} \Bigg\{ \| \pa_t^{\a}\nabla_{\bar{x}}^{\b} (\Psi,\Theta)(t)\|^2_{L^2_{l_r}}+\int_0^t \|\pa_y \pa_t^{\a}\nabla_{\bar{x}}^{\b} (\Psi,\Theta)(s)\|^2_{L^2_{l_r}}+\xi\|\nabla_{\sp} \pa_t^{\a}\nabla_{\bar{x}}^{\b} (\Psi,\Theta)(s)\|^2_{L^2_{l_r}}ds\Bigg\}\nonumber\\
&\leq C(t,\|(\rho,\fu,T)\|_{W^{r+1,\infty}}) \Bigg\{\sum_{j=0}^{r} \sum_{2\tilde\alpha+|\tilde\beta|=j} \Big[\|\pa_t^{\a}\nabla_{\bar{x}}^{\b} (\Psi,\Theta)(0)\|^2_{L^2_{l_j}}+\int_0^t\|\pa_t^{\tilde\a}\nabla_{\bar{x}}^{\tilde\b}(\tilde{\mathfrak{f}},\tilde{\mathfrak{g}})(s)\|^2_{L^2_{l_j}} ds\Big]\nonumber\\
&\quad
+C \sum_{j=0}^{r-2}\sum_{2\tilde\alpha+|\tilde\beta|=j} \|\pa_t^{\tilde\a}\nabla_{\bar{x}}^{\tilde\b}(\tilde{\mathfrak{f}},\tilde{\mathfrak{g}})(0)\|^2_{L^2_{l_j}}\Bigg\}.\nonumber
\end{align}
Hence by the induction arguments, we have
\begin{align}\label{5.142}
&\sum_{j=0}^{k}\sum_{2\alpha+|\beta|=j} \Bigg\{ \| \pa_t^{\a}\nabla_{\bar{x}}^{\b} (\Psi,\Theta)(t)\|^2_{L^2_{l_j}}+\int_0^t \|\pa_y \pa_t^{\a}\nabla_{\bar{x}}^{\b} (\Psi,\Theta)(s)\|^2_{L^2_{l_j}}+\xi\|\nabla_{\sp} \pa_t^{\a}\nabla_{\bar{x}}^{\b} (\Psi,\Theta)(s)\|^2_{L^2_{l_j}}ds\Bigg\}\nonumber\\
&\leq C(t,\|(\rho,\fu,T)\|_{W^{k+1,\infty}}) \Bigg\{\sum_{j=0}^{k} \sum_{2\tilde\alpha+|\tilde\beta|=j} \Big[\|\pa_t^{\a}\nabla_{\bar{x}}^{\b} (\Psi,\Theta)(0)\|^2_{L^2_{l_j}}+\int_0^t\|\pa_t^{\tilde\a}\nabla_{\bar{x}}^{\tilde\b}(\tilde{\mathfrak{f}},\tilde{\mathfrak{g}})(s)\|^2_{L^2_{l_j}} ds\Big]\nonumber\\
&\quad
+C \sum_{j=0}^{k-2}\sum_{2\tilde\alpha+|\tilde\beta|=j} \|\pa_t^{\tilde\a}\nabla_{\bar{x}}^{\tilde\b}(\tilde{\mathfrak{f}},\tilde{\mathfrak{g}})(0)\|^2_{L^2_{l_j}}\Bigg\}.
\end{align}

\

{\it Step 3. Taking limits $\sigma, \xi\rightarrow0+$.} Based on the uniform estimates \eqref{5.142}, we can  firstly take the limit $\sigma\rightarrow0+$, and then $\xi\rightarrow0+$.  Then, by using \eqref{5.09}-\eqref{5.011} and \eqref{5.142},  Lemma \ref{lem4.1} is proved, the details are omitted for simplicity of presentation.  $\hfill\Box$

\begin{remark}
In the proof of Lemma \ref{lem4.1}, it is direct to check that we can improve the polynomial decay of $y$ to exponential decay if $\pa_3\fu_3^0<0$. The key point is that we can obtain some  good term from \eqref{5.78}.
\end{remark}

%\vspace{1.5mm}

%%%%%%%%%%%%%%%%%%%%%%%%%%%%%%%%%%%%%%%%%%%%%%%%%%%%%%%%%%%%%%%%%%%%%%%%%%%%%%%%%%

\section{Construction on the Solutions of Expansions}\label{section5}

We define the velocity  weight functions
\begin{align}\label{ewf}
\tw_{\k_i}(v)=w_{\k_i}(v) \mu^{-\fa},
\quad \fw_{\bar{\k}_i}(v)=w_{\bar{\k}_i}(v) \mu_0^{-\fa} \,
\mbox{and}\
\fw_{\hat{\k}_i}(v)=w_{\hat{\k}_i}(v) \mu_0^{-\fa},
\end{align}
for constants $\k_i, \bar{\k}_i, \hat{\k}_i\geq 0,\, 1\leq i\leq N$ and $0\leq \fa <\frac{1}{2}$. It is noted that the weight function $\tw_{\k_i}$ depends on $(t,x)$, while $\fw_{\bar{\k}_i}$ and $\fw_{\hat{\k}_i}$ depend on $(t,x_{\sp} )$.  For later use, we define the notations
\begin{align}\nonumber
\quad \hat{x}=(x_{\sp} ,\eta)\in \R_+^3, \quad \nabla_{\hat{x}}:=(\nabla_{\sp} ,\pa_{\eta}),
\end{align}
and recall the notations $\bar{x}=(x_{\sp} ,y) \in \R_+^3$ and $\nabla_{\bar{x}}=(\nabla_{\sp} , \pa_y)$ in \eqref{5.59-2}, and the weighted $L^2_{l}$-norm in \eqref{5.59-1}.

%We construct the solutions of interior expansion, viscous and Knudsen boundary layers.
\begin{proposition} \label{prop5.1}
Let $0\leq \fa<\frac12$ in \eqref{ewf}. Let  $s_0, s_i, \bar{s}_i, \hat{s}_i\in \mathbb{N}_+$, $\k_i, \bar{\k}_{i}, \hat{\k}_i \in \R_+$ for $1\leq i\leq N$; and define  $l_j^i:=\bar{l}_i+2(\bar{s}_i-j)$ for $1\leq i \leq N, \  0\leq j\leq \bar{s}_i$. For these parameters, we assume the restrictions  \eqref{7.52}-\eqref{7.54} hold. We impose the IBVP \eqref{2.2}, \eqref{2.51}, \eqref{2.50} with the initial data $(\rho_{i}, u_i, \theta_i)(0)$; and impose IBVP \eqref{2.19-2}-\eqref{2.20}, \eqref{2.55}-\eqref{2.56} with initial data $(\bar{u}_{i,\sp},\bar{\theta}_{i})(0)$ such that
\begin{align}\label{7.1-0}
\sum_{i=0}^{N}\Big\{\sum_{\gamma+\beta\leq s_{i}}\|\pa_t^\gamma\nabla_x^{\beta}(\rho_{i}, u_{i},  \theta_{i})(0)\|_{L^2_{x}}+\sum_{j=0}^{\bar{s}_{i}} \sum_{j=2\gamma+\beta}   \|\pa_t^{\gamma} \nabla_{\bar{x}}^{\beta} (\bar{u}_{i,\sp},\bar{\theta}_{i})(0)\|^2_{L^2_{l_j^{i}}}
\Big\}<\infty.
\end{align}
And we also assume that the compatibility conditions are satisfied at the boundary.
Then there exist solutions $F_i=\sqrt{\mu}f_i, \, \bar{F}_i=\sqrt{\mu_0}\bar{f}_i, \, \hat{F}_i=\sqrt{\mu_0}\hat{f}_i$ to \eqref{1.7-1}, \eqref{1.14}, \eqref{1.19-1} over the time interval  $t\in[0,\tau]$, respectively. Moreover, we have  the following uniform estimates
\begin{align}\label{7.2-0}
\begin{split}
&\sup_{t\in[0,\tau]} \sum_{i=1}^N\Bigg\{\sum_{\gamma+\beta\leq s_i} \|\tw_{\k_i}\pa_t^\gamma \nabla_x^{\beta} f_{i}(t)\|_{L^2_xL^\infty_v}
+\sum_{j=0}^{\bar{s}_i}\sum_{j=2\gamma+\beta} \|\fw_{\bar{\k}_i}\pa_t^\gamma \nabla_{\bar{x}}^{\beta} \bar{f}_{i}(t)\|_{L^2_{l^i_j}L^\infty_v}\\
&\qquad\qquad\qquad+\sum_{\gamma+\beta\leq \hat{s}_i} \| e^{\zeta_i\cdot \eta}\fw_{\hat{\k}_i}\pa_t^\gamma \nabla_{\sp} ^{\beta} \hat{f}_{i}(t)\|_{L^\infty_{\hat{x},v}\cap L^2_{x_{\sp} }L^\infty_{\eta,v}}\Bigg\}\\
&\leq C\Bigg( \tau, \|(\varphi_0,\Phi_0, \vartheta_0)\|_{H^{s_0}}+\sum_{i=0}^{N}\sum_{\gamma+\beta\leq s_{i}}\|\pa_t^\gamma\nabla_x^{\beta}(\rho_{i}, u_{i},  \theta_{i})(0)\|_{L^2_{x}}\\
&\qquad\qquad\qquad+\sum_{i=0}^{N}\sum_{j=0}^{\bar{s}_{i}} \sum_{j=2\gamma+\beta}   \|\pa_t^{\gamma} \nabla_{\bar{x}}^{\beta} (\bar{u}_{i,\sp},\bar{\theta}_{i})(0)\|^2_{L^2_{l_j^{i}}}
\Bigg),
\end{split}
\end{align}
where the positive constants $\zeta_i>0 \, (i=1,\cdots,N)$ satisfying $\zeta_{i+1}=\frac12\zeta_i$ and $\zeta_1=1$.
\end{proposition}

\begin{remark}
Since the Knudsen boundary layer $\hat{f}_i$ is indeed a stationary problem with $(t,x_{\sp} )$ as parameters, hence there is no necessary to give initial data for $\hat{f}_i$.
Here the functions  $f_i, \bar{f}_i$ are smooth, however  $\pa_t^{\gamma}\nabla_{\sp} ^{\beta}\hat{f}_{i}$  is only continuous away from the grazing set $[0,\tau]\times \gamma_0$. %$\{(x_{\sp} ,0, v)\, | \, x_{\sp} \in\R^2, \, v_{\sp} \in\R^2, v_3\neq0 \}$.
\end{remark}

\noindent{\bf Proof.} Since the proof is very complicate, we divide the proof into several steps.\vspace{1.5mm}

\noindent{\it Step 1. Construction of solutions $f_1, \bar{f}_1$ and $\hat{f}_1$.}  \vspace{1.5mm}

{\it Step 1.1. Construction of  solution  $f_1$.} Noting $f_1\in \mathcal{N}$,  we need only to construct the macroscopic part $(\rho_1,u_1,\theta_1)$. Hence we consider \eqref{2.2} with $k=0$, and impose it with the boundary condition \eqref{2.51}. Then by using Lemma \ref{lem3.1}, we establish the existence of smooth solution of \eqref{2.2} (with $k=0$), \eqref{2.51} with the following estimate
\begin{align*}
\sup_{t\in[0,\tau]} \sum_{\gamma+\beta\leq s_1}\|\pa_{t}^\gamma \nabla_x^\beta (\rho_1,u_1,\theta_1)(t)\|_{L^{2}(\R_+^3)}&\leq
C(\tau, E_{2+s_1}) \sum_{\gamma+\beta\leq s_1}\|\pa_{t}^\gamma \nabla_x^\beta (\rho_1,u_1,\theta_1)(0)\|_{L^{2}(\R_+^3)},
\end{align*}
with $s_1\gg1$ such that $s_0\geq 2+s_1>0$, and $E_k$ is the one defined in Lemma \ref{lem3.1}. Therefore we have proved the existence of smooth solution $f_1$ over $[0,\tau]$ with
\begin{align*}%\label{7.1}
\sup_{t\in[0,\tau]} \|\tw_{\k_1}f_1(t)\|_{L^{2}_{x}L^\infty_v}&\leq
C(\tau, E_{2+s_1}) \sum_{\gamma+\beta\leq s_1}\|\pa_{t}^\gamma \nabla_x^\beta (\rho_1,u_1,\theta_1)(0)\|_{L^{2}(\R_+^3)},
\end{align*}
for any $\k_1>0$.
 \vspace{1.5mm}

{\it Step 1.2.  Construction of  solution $\bar{f}_1$.} Noting  \eqref{Boussinesq}, \eqref{2.12} and \eqref{2.14}, we need only to calculate $(\bar{u}_{1,\sp}, \bar{\theta}_1)$. Taking $k=2$ in \eqref{2.55}-\eqref{2.56}, then using \eqref{2.14} and the facts
\begin{equation*}
\bar{J}_0=0, \quad (\hat{A}_2,\hat{B}_2,\hat{C}_2)=(0,0,0),
\end{equation*}
we have
\begin{align}\label{7.2}
\begin{cases}
\dis\partial_{y} \bar{u}_{1,i}(t,x_{\sp} ,0)
=\frac{1}{\mu(T^0)} \big\langle T^0 \mathcal{A}_{3i}^0, \  (\mathbf{I-P_0})f_2\big\rangle(t,x_{\sp} ,0), \  i=1,2,\\[3mm]
\dis\partial_y\bar{\theta}_{1}(t,x_{\sp} ,0)
=\frac{1}{\kappa(T^0)} \big\langle 2(T^0)^{\frac32}\mathcal{B}_3^0,\   (\mathbf{I-P_0})f_2 \big\rangle(t,x_{\sp} ,0).
\end{cases}
\end{align}

It follows from \eqref{2.1} that
\begin{equation*}%\label{7.3}
\dis \{\mathbf{I-P}\}f_2=\mathbf{L}^{-1}\left(-\f{\{\partial_t+v\cdot\nabla_{x}\}\mu}{\sqrt{\mu}}+\frac{1}{\sqrt{\mu}}Q(\sqrt{\mu}f_1,\sqrt{\mu}f_1)\right).
\end{equation*}
Noting \eqref{0.8}, it is direct to  have that
\begin{align}\label{7.4}
\{\mathbf{I-P}\} \left(-\f{\{\partial_t+v\cdot\nabla_{x}\}\mu}{\sqrt{\mu}}\right)
%&=\{\mathbf{I-P}\} \left(-\f{(v-\fu)\cdot\nabla_{x}\mu}{\sqrt{\mu}}\right)
=\sum_{j,l=1}^3\partial_j \fu_l \mathcal{A}_{jl}+\sum_{j=1}^3 \frac{\partial_j T}{\sqrt{T}} \mathcal{B}_j.
\end{align}
Since $f_1\in\mathcal{N}$, by similar arguments as in \eqref{0.1}-\eqref{0.2}, it holds that
\begin{align}
&\mathbf{L}^{-1}\left(\frac{1}{\sqrt{\mu}}Q(\sqrt{\mu}f_1,\sqrt{\mu}f_1)\right)\nonumber\\
&=\sum_{j,l=1}^3\frac{1}{2T}u_{1,l}\, u_{1,j} \mathcal{A}_{lj}+\frac{\theta_1}{3T^{\frac32}} u_1\cdot \mathcal{B}+\frac{\theta_1^2}{72T^2} \{\mathbf{I-P}\}\left\{ (\frac{|v-\fu|^2}{T}-5)^2\sqrt{\mu} \right\},\nonumber
\end{align}
which, together  with  \eqref{7.4}, yields that
\begin{align}\label{7.6}
\{\mathbf{I-P}\}f_2&=-\mathbf{L}^{-1}\left\{\sum_{j,l=1}^3\partial_j \fu_l \mathcal{A}_{jl}+\sum_{j=1}^3 \frac{\partial_j T}{T} \mathcal{B}_j \right\}+\sum_{j,l=1}^3\frac{1}{2T}u_{1,l}\, u_{1,j} \mathcal{A}_{lj}\nonumber\\
&\quad+\frac{\theta_1}{3T^{\frac32}} u_1\cdot \mathcal{B}+\frac{\theta_1^2}{72T^2} \{\mathbf{I-P}\}\left\{ (\frac{|v-\fu|^2}{T}-5)^2\sqrt{\mu} \right\}.
\end{align}
Since $\mathbf{L}^{-1}$ preserves the decay property of $v$, it is direct to check that
\begin{equation}\nonumber
|\{\mathbf{I-P}\}f_2(t,x,v)|\leq C(|\nabla(\fu,T)|+|(u_1,\theta_1)|^2)(t,x) \, (1+|v|)^{4} \sqrt{\mu}.
\end{equation}

Noting \eqref{2.51} and $\pa_1\fu_3(t,x_{\sp} ,0)=\pa_2\fu_3(t,x_{\sp} ,0)=0$,  and using \eqref{2.15-2},  \eqref{7.6},  a direct calculation shows that
\begin{align}\label{7.7}
\big\langle T^0 \mathcal{A}_{3i}^0, \,  (\mathbf{I-P})f_2 \big\rangle(t,x_{\sp} ,0)
=-\pa_3\fu_i^0 \big\langle T^0 \mathcal{A}_{3i}^0, \,  \mathbf{L}_0^{-1}\mathcal{A}_{3i}^0\big\rangle=-\mu(T^0)\, \pa_3\fu_i^0,\quad \mbox{for}\  i=1,2,
\end{align}
and
\begin{align}\label{7.8}
\langle 2(T^0)^{\frac32}\mathcal{B}_3^0,\   (\mathbf{I-P})f_2\rangle(t,x_{\sp} ,0)
=-2T^0\partial_3T^0 \langle \mathcal{B}_3^0,\,  \mathbf{L}_0^{-1} \mathcal{B}_3^0 \rangle
=-3\kappa(T^0)\, \partial_3T^0.
\end{align}
Substituting \eqref{7.7} and \eqref{7.8} into \eqref{7.2}, we get the exact expression of  boundary condition for viscous boundary layer $(\bar{u}_{1,i},\bar{\theta}_{1})$, i.e.,
\begin{equation}\label{7.9}
\partial_{y} \bar{u}_{1,i}(t,x_{\sp} ,0)=-\partial_3\fu_i^0, \quad \partial_{y} \bar{\theta}_{1}(t,x_{\sp} ,0)=-3\pa_3 T^0.
\end{equation}

By using Lemma \ref{lem4.1}, we can obtain the existence of smooth solution of \eqref{2.15}-\eqref{2.15-1} and \eqref{7.9} over $[0,\tau]\times\R_+^3$ satisfying
\begin{align}\label{7.10}
&\sum_{j=0}^{\bar{s}_1} \sum_{\beta+2\gamma=j}\left\{\|\pa_t^\gamma \nabla_{\bar{x}}^{\beta} (\bar{u}_{1,\sp}, \bar{\theta}_1)(t)\|^2_{L^2_{l_{j}^{1}}}+\int_0^\tau \|\pa_t^\gamma \nabla_{\bar{x}}^{\beta} \pa_y(\bar{u}_{1,\sp}, \bar{\theta}_1)(t)\|_{L^2_{l_{j}^{1}}}^2dt\right\}\nonumber\\
&\leq C(\tau,E_{3+\bar{s}_1}) \Bigg\{\sum_{j=0}^{\bar{s}_1} \sum_{\beta+2\gamma=j}  \|\pa_t^\gamma \nabla_{\bar{x}}^{\beta} (\bar{u}_{1,\sp}, \bar{\theta}_1)(0)\|^2_{L^2_{l_{j}^{1}}}\nonumber\\
&\qquad+\sup_{t\in[0,\tau]}\sum_{\beta+2\gamma\leq \bar{s}_1+2}  \|\pa^{\beta}_{\sp}  \pa_t^{\gamma}(\pa_3\fu^0_1,\, \pa_3\fu^0_2,\, \pa_3T^0)(t)\|^2_{L^2(\R^2)} \Bigg\}\nonumber\\
&\leq C\Big(\tau,E_{4+\bar{s}_1}, \sum_{j=0}^{s_1} \sum_{\beta+2\gamma=j}  \|\pa_t^\gamma \nabla_{\bar{x}}^{\beta} (\bar{u}_{1,\sp}, \bar{\theta}_1)(0)\|^2_{L^2_{l_{j}^{1}}}\Big),
\end{align}
where  $l_{j}^{1}=\bar{l}_1+2(\bar{s}_1-j)$ with $\bar{l}_1\gg1$ and $s_0\geq 4+\bar{s}_1$. Combining \eqref{7.10} with \eqref{2.14}, we get
\begin{align*}
\sum_{j=0}^{\bar{s}_1} \sum_{\beta+2\gamma=j}\|\fw_{\bar{\k}_1}\partial_t^{\gamma} \nabla_{\bar{x}}^{\beta} \bar{f}_{1}(t)\|_{L^2_{l^1_j}L^\infty_v}\leq C\Big(\tau,E_{4+\bar{s}_1}, \sum_{j=0}^{\bar{s}_1} \sum_{\beta+2\gamma=j}  \|\pa_t^\gamma \nabla_{\bar{x}}^{\beta} (\bar{u}_{1,\sp},\bar{\theta}_1)(0)\|^2_{L^2_{l_{j}^{1}}}\Big),
\end{align*}
for any  $\bar{\k}_1>0$.

\

{\it Step 1.3.  Construction of solution $\hat{f}_1$.}  From \eqref{2.43-2}, we know $\hat{f}_{1,1}\equiv0$, we need only to consider $\hat{f}_{1,2}$. Noting \eqref{2.45}(with $k=1$), using \eqref{2.51} and $\eqref{2.14}_1$, one has the boundary condition for $\hat{f}_{1,2}$
\begin{equation}\nonumber
\hat{f}_{1,2}(t,x_{\sp} ,0,v_{\sp} ,v_3)|_{v_3>0}=\hat{f}_{1,2}(t,x_{\sp} ,0,v_{\sp} ,-v_3),
\end{equation}
which, together with Lemma \ref{lem2.6}, yields the existence of $\hat{f}_{1,2}$ with $\hat{f}_{1,2}^\v\equiv0$. Therefore we have proved the existence of $\hat{f}_{1}$ with
\begin{equation}\label{7.11}
\hat{f}_{1}\equiv0.
\end{equation}
That is the $\v$-th order Knudsen boundary layer does not appear. This is reasonable since the Knudsen boundary layer is used to mend the boundary condition at higher orders.

\

{\it Step 2. Construction of solutions  $f_{k}, \bar{f}_k$ and $\hat{f}_k$.} We shall use induction argument. Suppose we have already proved the existence of $f_{i}, \bar{f}_{i}$ and $\hat{f}_{i}$ for $1\leq i\leq k$ such that
\begin{align}\label{7.11-1}
D_k+\bar{D}_k+\hat{D}_k&\leq C\bigg(\tau, E_{s_0}, \sum_{i=0}^{k}\sum_{j=0}^{\bar{s}_{i}} \sum_{j=2\gamma+\beta}   \|\pa_t^{\gamma} \nabla_{\bar{x}}^{\beta} (\bar{u}_{i,\sp},\bar{\theta}_{i})(0)\|^2_{L^2_{l_j^{i}}}\nonumber\\
&\qquad\qquad+\sum_{i=0}^{k}\sum_{\gamma+\beta\leq s_{i}}\|\pa_t^\gamma\nabla^{\beta}(\rho_{i}, u_{i},  \theta_{i})(0)\|_{L^2_{x}}\bigg)
\end{align}
where
\begin{align*}%\label{7.11-2}
\begin{split}
&D_k:=\sup_{t\in[0,\tau]} \Bigg\{\sum_{i=1}^k\sum_{\gamma+\beta\leq s_i} \|\tw_{\k_i}\pa_t^\gamma \nabla_x^{\beta} f_{i}(t)\|_{L^2_xL^\infty_v}\Bigg\}<\infty,\\
&\bar{D}_k:=\sup_{t\in[0,\tau]} \Bigg\{\sum_{i=1}^k\sum_{j=0}^{\bar{s}_i}\sum_{j=2\gamma+\beta} \|\fw_{\bar{\k}_i}\pa_t^\gamma \nabla_{\bar{x}}^{\beta} \bar{f}_{i}(t)\|_{L^2_{l^i_j}L^\infty_v}\Bigg\}<\infty,\\
&\hat{D}_k:=\sup_{t\in[0,\tau]} \Bigg\{\sum_{i=1}^k\sum_{\gamma+\beta\leq \hat{s}_i} \| e^{\zeta_i\cdot \eta}\fw_{\hat{\k}_i}\pa_t^\gamma \nabla_{\sp} ^{\beta} \hat{f}_{i}(t)\|_{L^\infty_{\hat{x},v}\cap L^2_{x_{\sp} }L^\infty_{\eta,v}}\Bigg\}<\infty,
\end{split}
\end{align*}
for some  $s_{i}>\bar{s}_{i}>\hat{s}_{i}\geq s_{i+1}>\bar{s}_{i+1}>\hat{s}_{i+1}\geq 1$, $\k_{i}\gg \bar{\k}_{i}\gg\hat{\k}_{i}\gg \k_{i+1}\gg \bar{\k}_{i+1}\gg\hat{\k}_{i+1}\gg1$ with $1\leq i\leq k-1$, and  $l_j^i=\bar{l}_i+2(\bar{s}_i-j)$, $\bar{l}_i\gg1$ with $1\leq i\leq k$ and $0\leq j\leq \bar{s}_i$. In the following, we consider the existence of $f_{k+1}, \bar{f}_{k+1}$ and $\hat{f}_{k+1}$.

\smallskip

{\it Step 2.1. Construction of solution $f_{k+1}$.} Let $r\geq 3$, Since $\mathbf{L}^{-1}$ preserves the decay property of $v$, by using  \eqref{2.1} and Sobolev inequality, we have   that
\begin{align}\label{7.12}
\sum_{\gamma+|\beta|\leq r}\|\tw_{\k_{k+1}}\pa_t^\gamma \nabla_x^{\beta} \{\mathbf{I-P}\} f_{k+1}(t)\|_{L^2_xL^\infty_v}
\leq C\Big(E_{r+1},\sum_{j=1}^k\sum_{\gamma+\beta\leq r+1} \|\tw_{\k_j}\pa_t^\gamma \nabla_x^{\beta} f_{j}(t)\|_{L^2_xL^\infty_v}\Big).
\end{align}

To obtain the solution of $f_{k+1}$, we still need to obtain the estimate for macroscopic part. For the source terms on RHS of \eqref{2.2}, it follows from \eqref{2.3}, \eqref{7.12} that
\begin{align}\label{7.14}
\sum_{\gamma+|\beta|\leq r+1}\|\pa_t^\gamma \nabla_x^{\beta} (\mathfrak{f}_{k},\mathfrak{g}_{k})(t)\|_{L^2_x} \leq C\Big(E_{r+3},\sum_{j=1}^k\sum_{\gamma+\beta\leq r+3} \|\tw_{\k_j}\pa_t^\gamma \nabla_x^{\beta} f_{j}(t)\|_{L^2_xL^\infty_v}\Big).
\end{align}

Before applying Lemma \ref{lem3.1}, we  need to estimate the boundary condition. Noting \eqref{2.50} (with $k$ replaced by $k+1$), we have that %It is direct to check that the boundary value  of \eqref{2.2}  given by \eqref{2.50}  is already a know function.
\begin{align}\label{7.15}
\sum_{\gamma+\beta\leq r+2}\|\partial_t^\gamma \pa_{\sp} ^{\beta} u_{k+1,3}(t,\cdot,0)\|_{L^2}
&\leq C(E_{r+4}) \Big\{ \sum_{\gamma+\beta\leq r+3}\int_0^\infty\|\partial_t^\gamma \pa_{\sp} ^{\beta} (\bar{\rho}_k,\bar{u}_k)(t,\cdot, y)\|_{L^2(\R^2)} dy\nonumber\\
&+ \sum_{\gamma+\beta\leq r+2}\|\partial_t^\gamma \pa_{\sp} ^{\beta} (\hat{A}_{k+1},\hat{C}_{k+1})(t,\cdot,0)\|_{L^2(\R^2)}\Big\}.
\end{align}
Noting $l_{j}^{k}>1$, a direct calculation shows that
\begin{align}\label{7.16}
 \sum_{\gamma+\beta\leq r+3} \int_0^\infty\|\partial_t^\gamma \pa_{\sp} ^{\beta} (\bar{\rho}_k,\bar{u}_k)(t,\cdot, y)\|_{L^2(\R^2)} dy
 \leq C\sum_{j=0}^{2(r+3)}\sum_{\beta+2\gamma=j}  \|\partial_t^\gamma \pa_{\sp} ^{\beta} \bar{f}_k(t)\|_{L^2_{l_j^{k}}L^\infty_v}.
\end{align}
By using \eqref{2.40}, \eqref{2.37} and \eqref{2.35}, one obtains
\begin{align}\label{7.17}
&\sum_{\gamma+\beta\leq r+2}\|\partial_t^\gamma \pa_{\sp} ^{\beta} (\hat{A}_{k+1},\hat{B}_{k+1}, \hat{C}_{k+1})(t,\cdot,0)\|_{L^2(\R^2)}\nonumber\\
&\leq C(E_{r+4}) \sum_{\gamma+\beta\leq r+2} \int_0^\infty \|\partial_t^\gamma \pa_{\sp} ^{\beta} (\hat{a}_{k+1},\hat{b}_{k+1}, \hat{c}_{k+1})(t,\cdot,\eta)\|_{L^2_{x_{\sp} }} d\eta\nonumber\\
&\leq C(E_{r+4}) \sum_{\gamma+\beta\leq r+2}  \|e^{\f12 \zeta_{k-1} \cdot \eta}\partial_t^\gamma \pa_{\sp} ^{\beta} (\hat{a}_{k+1},\hat{b}_{k+1}, \hat{c}_{k+1})(t)\|_{L^2_{x_{\sp} }L^\infty_{\eta}} \nonumber\\
%&\leq C(E_{r+5})  \sum_{\gamma+\beta\leq r+3} \| e^{\frac14\zeta_{k-1}\cdot \eta} \, w_{\hat{\alpha}_{k-1}} \partial_t^\gamma \pa_{\sp} ^{\beta} \hat{f}^\v_{k-1}(t)\|_{L^2_{\hat{x}}L^\infty_v}\nonumber\\
&\leq C(E_{r+5})  \sum_{\gamma+\beta\leq r+3} \| e^{\frac12\zeta_{k-1}\cdot \eta} \, \fw_{\hat{\k}_{k-1}} \partial_t^\gamma \pa_{\sp} ^{\beta} \hat{f}_{k-1}(t)\|_{L^2_{x_{\sp} }L^\infty_{\eta,v}}.
\end{align}
Here we emphasize that although only time and tangential derivatives are available for the knudsen boundary layer $\hat{f}_i$, but the trace $(\hat{A}_{k+1},\hat{B}_{k+1}, \hat{C}_{k+1})(t,\cdot,0)$  (equivalent to $\hat{f}_{k+1,1}(t,x_{\sp},0,v)$ ) is indeed well-defined from \eqref{2.40}.

Substituting \eqref{7.16}-\eqref{7.17} into \eqref{7.15}, it holds that
\begin{align}\label{7.18}
\sum_{\gamma+\beta\leq r+2}\|\partial_t^\gamma \pa_{\sp} ^{\beta} u_{k+1,3}(t,\cdot,0)\|_{L^2}
&\leq C(E_{r+5}) \Big\{ \sum_{j=0}^{2(r+3)}\sum_{\beta+2\gamma=j}  \|\partial_t^\gamma \pa_{\sp} ^{\beta} \bar{f}_k(t)\|_{L^2_{l_j^{k}}L^\infty_v} \nonumber\\
&\hspace{-10mm}+ \sum_{\gamma+\beta\leq r+3} \| e^{\frac12\zeta_{k-1}\cdot \eta}  \fw_{\hat{\k}_{k-1}}  \partial_t^\gamma \pa_{\sp} ^{\beta} \hat{f}_{k-1}(t)\|_{L^2_{x_{\sp} }L^\infty_{\eta,v}}\Big\}.
\end{align}
Noting \eqref{5.4-1}, the above time and tangential derivatives are enough when using Lemma \ref{lem3.1}.

\smallskip

Now applying Lemma \ref{lem3.1}, using  \eqref{7.14} and  \eqref{7.18}, we obtain
\begin{align}\label{7.19}
&\sum_{\gamma+\beta\leq r}\|\pa_t^\gamma\nabla_x^{\beta}(\rho_{k+1}, u_{k+1},\theta_{k+1})(t)\|_{L^2_{x}}\nonumber\\
&\leq C\bigg(\tau, E_{r+5}, \sum_{\gamma+\beta\leq r}\|\pa_t^\gamma\nabla_x^{\beta}(\rho_{k+1}, u_{k+1},\theta_{k+1})(0)\|_{L^2_{x}}, \sup_{t\in[0,\tau]}\Big[\sum_{j=0}^{2(r+3)}\sum_{\beta+2\gamma=j}  \|\partial_t^\gamma \pa_{\sp} ^{\beta} \bar{f}_k(t)\|_{L^2_{l_j^{k}}L^\infty_v} \nonumber\\
&\qquad+ \sum_{\gamma+\beta\leq r+3} \| e^{\frac12\zeta_{k-1}\cdot \eta}  \fw_{\hat{\k}_{k-1}}  \partial_t^\gamma \pa_{\sp} ^{\beta} \hat{f}_{k-1}(t)\|_{L^2_{x_{\sp} }L^\infty_{\eta,v}}+\sum_{j=1}^k\sum_{\gamma+\beta\leq r+3} \|\tw_{\k_j}\pa_t^\gamma  \nabla_x^{\beta} f_{j}(t)\|_{L^2_{x}L^\infty_v}\Big] \bigg).
\end{align}
Now taking $r=s_{k+1}$ so that
\begin{equation}\nonumber
s_0\geq 5+s_{k+1},\quad s_{k}\geq s_{k+1}+3,\quad \hat{s}_{k-1}\geq s_{k+1}+3\quad \mbox{and}\   \bar{s}_{k}\geq 2(s_{k+1}+3).
\end{equation}
then, combining \eqref{7.19} with \eqref{7.12}, we  get
\begin{align}\label{7.20}
&\sum_{\gamma+\beta\leq s_{k+1}}\|\tw_{\k_{k+1}}\pa_t^\gamma\nabla_x^{\beta}f_{k+1}(t)\|_{L^2_xL^\infty_v} \nonumber\\
&\leq C\Big(\tau, E_{s_0}, D_k+\bar{D}_{k}+\hat{D}_{k-1}, \sum_{\gamma+\beta\leq s_{k+1}}\|\pa_t^\gamma \nabla_x^{\beta}(\rho_{k+1}, u_{k+1},  \theta_{k+1})(0)\|_{L^2_{x}}\Big).
\end{align}
Since we have obtained $f_i$ with $1\leq i\leq k+1$, by using \eqref{2.1}, it holds that
\begin{align}\label{7.12-1}
&\sum_{\gamma+|\beta|\leq s_{k+1}-1}\|\tw_{\k_{k+2}}\pa_t^\gamma \nabla_x^{\beta} \{\mathbf{I-P}\} f_{k+2}(t)\|_{L^2_xL^\infty_v}\nonumber\\
&\leq C\Big(E_{s_{k+1}},\sum_{i=1}^{k+1}\sum_{\gamma+\beta\leq s_{k+1}} \|\tw_{\k_i}\pa_t^\gamma \nabla_x^{\beta} f_{i}(t)\|_{L^2_xL^\infty_v}\Big)\nonumber\\
&\leq C\Big(\tau, E_{s_0}, D_k+\bar{D}_{k}+\hat{D}_{k-1}, \sum_{\gamma+\beta\leq s_{k+1}}\|\pa_t^\gamma\nabla_x^{\beta}(\rho_{k+1}, u_{k+1},  \theta_{k+1})(0)\|_{L^2_{x}}\Big),
\end{align}
which will be used when consider the trace of $(\pa_y\bar{u}_{k+1,\sp},\ \pa_y\bar{\theta}_{k+1} )(t,x_{\sp} ,0)$ in the following.

\

{\it Step 2.2. Construction of solution for $\bar{f}_{k+1}$.}  Noting $\mathbf{L}_0^{-1}$ preserves the decay property of $v$, then it follows from \eqref{2.21} that
\begin{align}\label{7.21}
&\sum_{2\gamma+\beta=j} \|\fw_{\bar{\k}_{k+1}} \pa_t^{\gamma} \nabla_{\bar{x}}^{\beta} \{\mathbf{I-P_0}\}\bar{f}_{k+1}(t)\|_{L^2_{l}L^\infty_v}\leq  C\bigg(E_{j+3}, \sum_{2\gamma+\beta\leq j+3} \|\tw_{\k_1}\pa_t^{\gamma}\nabla_x^{\beta}f_1(t)\|_{L^2_xL^\infty_v}\nonumber\\
&\qquad\qquad+\sum_{2\gamma+\beta\leq j} \|\fw_{\bar{\k}_1}\pa_t^{\gamma}\nabla_{\bar{x}}^{\beta}\bar{f}_1(t)\|_{L^2_{\bar{x}}L^\infty_v}+\sum_{2\gamma+\beta\leq j+1} \|\fw_{\bar{\k}_k}\pa_t^{\gamma}\nabla_{\bar{x}}^{\beta}\bar{f}_{k}(t)\|_{L^2_{2+l}L^\infty_v}
\bigg)\nonumber\\
&\qquad\qquad\quad +\sum_{2\gamma+\beta=j} \|\fw_{\bar{\k}_{k+1}} \pa_t^{\gamma}\nabla_{\bar{x}}^{\beta} \bar{J}_{k-1}(t)\|_{L^2_{l}L^\infty_v}.
\end{align}
Noting \eqref{2.26}, a direct calculation shows that
\begin{align}\label{7.22}
&\sum_{2\gamma+\beta=j} \|\fw_{\bar{\k}_{k+1}} \pa_t^{\gamma}\nabla_{\bar{x}}^{\beta} \bar{J}_{k-1}(t)\|_{L^2_{l}L^\infty_v}
\nonumber\\
&\leq C\Big(E_{j+\fb+2}, \sum_{i=0}^{k-1}\sum_{2\gamma+\beta\leq j+2} \|\fw_{\bar{\k}_i}\pa_t^{\gamma}\nabla_{\bar{x}}^{\beta} \bar{f}_{i}(t)\|_{L^2_{l+2\fb}L^\infty_v}
+ \sum_{i=0}^{k}\sum_{2\gamma+\beta\leq j+2+\fb} \|\tw_{\k_i}\pa_t^{\gamma}\nabla_x^{\beta} f_{i}(t)\|_{L^2_xL^\infty_v} \Big)\nonumber\\
&\quad+C(E_{j+\fb+2}) \sum_{2\gamma+\beta\leq j+1} \|\fw_{\bar{\k}_{k}}\pa_t^{\gamma}\nabla_{\bar{x}}^{\beta} \{\mathbf{I-P_0}\}\bar{f}_{k}(t)\|_{L^2_{l+2}L^\infty_v}.
\end{align}
Substituting \eqref{7.22} into \eqref{7.21}, and noting $l_j^{k+1}\ll l_j^i$ for $1\leq i\leq k$, we can obtain the estimate for microscopic part $\{\mathbf{I-P_0}\}\bar{f}_{k+1}$
\begin{align}\label{7.23}
\sum_{2\gamma+\beta=j} \|\tw_{\bar{\k}_{k+1}} \pa_t^{\gamma} \nabla_{\bar{x}}^{\beta} \{\mathbf{I-P_0}\}\bar{f}_{k+1}(t)\|_{L^2_{l}L^\infty_v}&\leq  C\Big(E_{j+\fb+2}, \,  \sum_{i=0}^{k}\sum_{2\gamma+\beta\leq j+\fb+2} \|\tw_{\k_i}\pa_t^{\gamma}\nabla_x^{\beta} f_{i}(t)\|_{L^2_xL^\infty_v}\nonumber\\
&+\sum_{i=1}^{k}\sum_{2\gamma+\beta\leq j+2} \|\fw_{\bar{\k}_i}\pa_t^{\gamma}\nabla_{\bar{x}}^{\beta}\bar{f}_{i}(t)\|_{L^2_{l+2\fb}L^\infty_v}
\Big).
\end{align}

On the other hand, substituting \eqref{7.23} (with $k$ replaced by $k-1$) into \eqref{7.22}, one can obtain a better estimate for $\bar{J}_{k-1}$
\begin{align}\label{7.24}
\sum_{2\gamma+\beta=j} \|\fw_{\bar{\k}_{k+1}} \pa_t^{\gamma}\nabla_{\bar{x}}^{\beta} \bar{J}_{k-1}(t)\|_{L^2_{l}L^\infty_v}
&\leq C\Big(E_{j+\fb+3}, \, \sum_{i=0}^{k-1}\sum_{2\gamma+\beta\leq j+3} \|\fw_{\bar{\k}_i}\pa_t^{\gamma}\nabla_{\bar{x}}^{\beta} \bar{f}_{i}(t)\|_{L^2_{l+2\fb}L^\infty_v} \nonumber\\
&\quad+ \sum_{i=0}^{k}\sum_{2\gamma+\beta\leq j+3+\fb} \|\tw_{\k_i}\pa_t^{\gamma}\nabla_x^{\beta} f_{i}(t)\|_{L^2_xL^\infty_v} \Big).
\end{align}
For $\bar{W}_{k-1}, \bar{H}_{k-1}$, it follows from \eqref{2.24-1}-\eqref{2.24-2} and \eqref{7.23}  that
\begin{align}\label{7.25}
&\sum_{2\gamma+\beta=j} \|\pa_t^{\gamma}\nabla_{\bar{x}}^{\beta}( \bar{W}_{k-1}, \bar{H}_{k-1})(t)\|_{L^2_l} \leq C(E_{j+3})  \sum_{2\gamma+\beta\leq j+1} \|\pa_t^{\gamma}\nabla_{\bar{x}}^{\beta} \{\mathbf{I-P_0}\}\bar{f}_{k}(t)\|_{L^2_lL^\infty_v}\nonumber\\
&\leq C\Big(E_{j+\fb+3}, \, \sum_{i=0}^{k-1}\sum_{2\gamma+\beta\leq j+3+\fb} \|\tw_{\k_i}\pa_t^{\gamma}\nabla_x^{\beta} f_{i}(t)\|_{L^2_xL^\infty_v}+\sum_{i=1}^{k-1}\sum_{2\gamma+\beta\leq j+3} \|\fw_{\bar{\k}_i}\pa_t^{\gamma}\nabla_{\bar{x}}^{\beta}\bar{f}_{i}(t)\|_{L^2_{l+2\fb}L^\infty_v}
\Big).
\end{align}

For $\bar{u}_{k+1,3}(t,x_{\sp} ,y)$, it follows from \eqref{2.24} and \eqref{2.11} that
\begin{align}\label{7.26}
\bar{u}_{k+1,3}(t,x_{\sp} ,y)=-\int_y^\infty\frac{1}{\rho^0}\Big\{\partial_t\bar{\rho}_{k}+\mbox{\rm div}_{\sp} (\rho^0\bar{u}_{k,\sp}+\bar{\rho}_{k} \fu^0_{\sp} )\Big\}(t,x_{\sp} ,z)dz,
\end{align}
which yields that
\begin{align}\label{7.27}
&\sum_{2\gamma+\beta=j}\|\pa_t^\gamma\nabla_{\sp} ^\beta\bar{u}_{k+1,3}(t)\|_{L^2_l}\nonumber\\
&\leq C(E_{j+3}) \int_0^\infty (1+y)^{l} dy \int_y^{\infty} \sum_{2\gamma+\beta\leq j+2} \|\pa_t^\gamma\nabla_{\sp} ^\beta(\bar{\rho}_{k}, \bar{u}_{k})(t,\cdot,z)\|_{L^2(\R^2)} dz\nonumber\\
&\leq C(E_{j+3}) \int_0^\infty (1+y)^{l} \Big(\int_y^\infty (1+z)^{-2l-4}dz\Big)^{\frac12} dy\sum_{2\gamma+\beta\leq j+2} \|\pa_t^\gamma\nabla_{\sp} ^\beta(\bar{\rho}_{k}, \bar{u}_{k})(t)\|_{L^2_{2l+4}} \nonumber\\
&\leq C(E_{j+3}) \sum_{2\gamma+\beta\leq j+2} \|\pa_t^\gamma\nabla_{\sp} ^\beta(\bar{\rho}_{k}, \bar{u}_{k})(t)\|_{L^2_{2l+4}}\nonumber\\
&\leq C(E_{j+3}) \sum_{2\gamma+\beta\leq j+2} \|\pa_t^\gamma\nabla_{\sp} ^\beta \bar{f}_{k}(t)\|_{L^2_{2l+4}L^\infty_v},\quad \mbox{for} \quad l\geq0.
\end{align}
On the other hand, we assume that  $\nabla_{\bar{x}}^{\beta}$  contains at least one $\pa_y$, then it follows from \eqref{2.24} that
\begin{align}\nonumber
\sum_{2\gamma+\beta=j}\|\pa_t^\gamma\nabla_{\bar{x}}^\beta\bar{u}^\v_{k+1,3}(t)\|_{L^2_l}
&\leq \sum_{2\gamma+\beta=j}\|\pa_t^\gamma\bar{\nabla}^{\beta-1} \pa_y\bar{u}^\v_{k+1,3}(t)\|_{L^2_l}\nonumber\\
&\leq C(E_{j+3}) \sum_{2\gamma+\beta\leq j+1} \|\pa_t^\gamma\nabla_{\bar{x}}^\beta \bar{f}^\v_{k}(t)\|_{L^2_{l}L^\infty_v},\quad \mbox{for} \quad l\geq0,\nonumber
\end{align}
which, together with \eqref{7.27}, yields
\begin{align}\label{7.28}
\sum_{2\gamma+\beta=j}\|\pa_t^\gamma\nabla_{\bar{x}}^\beta\bar{u}_{k+1,3}(t)\|_{L^2_l}
\leq C(E_{j+3}) \sum_{2\gamma+\beta\leq j+2} \|\pa_t^\gamma\nabla_{\bar{x}}^\beta \bar{f}_{k}(t)\|_{L^2_{2l+4}L^\infty_v},\quad \mbox{for} \quad l\geq0.
\end{align}

By using \eqref{2.25} and similar arguments as in \eqref{7.26}-\eqref{7.28}, one can obtain the estimate for $\bar{p}^\v_{k+1}$
\begin{align}\label{7.30}
&\sum_{2\gamma+\beta=j}\|\pa_t^\gamma\nabla_{\bar{x}}^\beta\bar{p}_{k+1,3}(t)\|_{L^2_l}
\leq C(E_{j+3}) \Big(\sum_{2\gamma+\beta\leq j+2} \|\pa_t^\gamma\nabla_{\bar{x}}^\beta \bar{f}_{k}(t)\|_{L^2_{2l+6}L^\infty_v}\nonumber\\
&\qquad\qquad+\sum_{2\gamma+\beta\leq j+1} \|\pa_t^\gamma\nabla_{\bar{x}}^\beta \bar{J}_{k-1}(t)\|_{L^2_{2l+4}L^\infty_v}+\sum_{2\gamma+\beta\leq j} \|\pa_t^\gamma\nabla_{\bar{x}}^\beta \bar{W}_{k-1}(t)\|_{L^2_{2l+4}}\Big) \nonumber\\
&\leq C\Big(E_{j+4+\fb}, \, \sum_{i=1}^k\sum_{2\gamma+\beta\leq j+4} \|\fw_{\bar{\k}_i}\pa_t^\gamma\nabla_{\bar{x}}^\beta \bar{f}_{i}(t)\|_{L^2_{2l+4+2\fb}L^\infty_v} \nonumber\\
&\qquad\qquad\qquad\qquad\qquad+\sum_{i=1}^k\sum_{2\gamma+\beta\leq j+4+\fb} \|\tw_{\k_i}\pa_t^\gamma\nabla_x^\beta f_{i}(t)\|_{L^2_xL^\infty_v}\Big),
\end{align}
where we have used \eqref{7.24}-\eqref{7.25} in the last inequality.

\

By similar arguments as  \eqref{7.27}, we can have the following trace estimate
\begin{align}\label{7.29}
\sum_{2\gamma+\beta=j}\|\pa_t^\gamma\nabla_{\sp} ^\beta\bar{u}_{k+1,3}(t,\cdot,0)\|_{L^2(\R^2)}
&\leq C(E_{j+3}) \int_0^{\infty} \sum_{2\gamma+\beta\leq j+2} \|\pa_t^\gamma\nabla_{\sp} ^\beta(\bar{\rho}_{k}, \bar{u}_{k})(t,\cdot,z)\|_{L^2(\R^2)} dz\nonumber\\
&\leq C(E_{j+3}) \sum_{2\gamma+\beta\leq j+2} \|\pa_t^\gamma\nabla_{\sp} ^\beta(\bar{\rho}_{k}, \bar{u}_{k})(t)\|_{L^2_l} \nonumber\\
&\leq  C(E_{j+3}) \sum_{2\gamma+\beta\leq j+2} \|\pa_t^\gamma\nabla_{\sp} ^\beta \bar{f}_{k}(t)\|_{L^2_{l}L^\infty_v},\quad \mbox{for} \quad l>1.
\end{align}
Taking $s_0\geq \bar{s}_{k+1}+\fb+6$, $s_k>\bar{s}_k>\hat{s}_k>s_{k+1}\geq \bar{s}_{k+1}+7+\fb$, and $l^i_j\geq 2\fb$ for $i\leq k, j\leq \bar{s}_{k}$,  then it follows from  \eqref{2.55}-\eqref{2.56} that
\begin{align}\label{7.33}
&\sum_{2\gamma+\beta\leq \bar{s}_{k+1}+2}\|\pa_t^\gamma\nabla_{\sp} ^\beta(\pa_y\bar{u}_{k+1,\sp},\ \pa_y\bar{\theta}_{k+1} )(t,\cdot,0)\|_{L^2(\R^2)}\nonumber\\
&\leq C\bigg(  E_{\bar{s}_{k+1}+4},  \sum_{2\gamma+\beta\leq \bar{s}_{k+1}+4} \big[\|\tw_{\k_1}\pa_t^\gamma\nabla_x^\beta f_1(t)\|_{L^2}+\|\fw_{\bar{\k}_{1}}\pa_t^\gamma \nabla_{\bar{x}}^\beta \bar{f}_1(t)\|_{L^2}\big]\nonumber\\
&\quad+\sum_{2\gamma+\beta\leq \bar{s}_{k+1}+2}\Big[ \|\pa_t^\gamma\nabla_{\sp} ^\beta\bar{u}_{k+1,3}(t,\cdot,0)\|_{L^2_{x_{\sp} }}+\|\pa_t^\gamma\nabla_{\sp} ^\beta\{\mathbf{I-P}\}f_{k+2}(t,\cdot,0,\cdot)\|_{L^2_{x_{\sp} }L^\infty_v}\nonumber\\
&\quad+ \|\pa_t^\gamma\nabla_{\sp} ^\beta\bar{J}_{k}(t,\cdot,0,\cdot)\|_{L^2_{x_{\sp} }L^\infty_v}+\|\pa_t^\gamma\nabla_{\sp} ^\beta(\hat{A}_{k+2},\hat{B}_{k+2},\hat{C}_{k+2})(t,\cdot,0)\|_{L^2_{x_{\sp} }}\Big]\bigg)\nonumber\\
&\leq C\Big(\tau, E_{s_0},\, D_k+\bar{D}_k+\hat{D}_k, \sum_{\gamma+\beta\leq s_{k+1}}\|\pa_t^\gamma\nabla_x^{\beta}(\rho_{k+1}, u_{k+1},  \theta_{k+1})(0)\|_{L^2_{x}}\Big)<\infty,
\end{align}
where we have used \eqref{7.29},\eqref{7.20}, \eqref{7.12-1}, \eqref{7.24},  \eqref{7.17} (with $k-1$ replaced by $k$ in \eqref{7.24},  \eqref{7.17}) and the trace theorem.

Similarly, for the source terms of \eqref{2.19-2}-\eqref{2.20}, by using \eqref{7.24}-\eqref{7.25} and \eqref{7.28}-\eqref{7.30}, a direct calculation shows that
\begin{align}\label{7.34}
\sum_{j=0}^{\bar{s}_{k+1}}\sum_{\beta+2\gamma=j}\|(\pa_y \pa_t^{\gamma}\nabla_{\bar{x}}^{\beta} (\bar{\mathfrak{f}}_{k},\bar{\mathfrak{g}}_{k}), \pa_t^{\gamma}\nabla_{\bar{x}}^{\beta}  (\bar{\mathfrak{f}}_{k},\bar{\mathfrak{g}}_{k}))(t)\|^2_{L^2_{l_j^{k+1}}}
\leq C(\tau, E_{s_0}, D_k+\bar{D}_k)
\end{align}
where we have taken
\begin{equation}\label{7.34-2}
s_0\geq \bar{s}_{k+1}+\fb+6,\quad s_{k+1}\geq \bar{s}_{k+1}+8+\fb,\quad l_{j}^i\geq 2 l_{j}^{k+1}+18+2\fb,\quad \mbox{for} \ 1\leq i\leq k,
\end{equation}
with $l_{j}^{k+1}=\bar{l}_{k+1}+2(\bar{s}_{k+1}-j)$ and  $\bar{l}_{k+1}\gg1$.

\

Using Lemma \ref{lem4.1} (the time and tangential derivatives estimate \eqref{7.33} for the boundary condition  are used enough when using Lemma \ref{lem4.1}) and \eqref{7.33}-\eqref{7.34}, and noting \eqref{7.34-2}, one can obtain
\begin{align}\label{7.35}
&\sum_{j=0}^{\bar{s}_{k+1}} \sum_{j=2\gamma+\beta} \Bigg( \|\pa_t^{\gamma} \nabla_{\bar{x}}^{\beta} (\bar{u}_{k+1,\sp},\bar{\theta}_{k+1})(t)\|^2_{L^2_{l_j^{k+1}}}+\int_{0}^\tau  \|\pa_t^{\gamma} \nabla_{\bar{x}}^{\beta} \pa_y(\bar{u}_{k+1,\sp},\bar{\theta}_{k+1})(t)\|^2_{L^2_{l_j^{k+1}}} dt\Bigg)\nonumber\\
&\leq C\bigg(\tau, E_{s_0}, D_k+\bar{D}_k+\hat{D}_k,  \sum_{j=0}^{\bar{s}_{k+1}} \sum_{j=2\gamma+\beta}   \|\pa_t^{\gamma} \nabla_{\bar{x}}^{\beta} (\bar{u}_{k+1,\sp},\bar{\theta}_{k+1})(0)\|^2_{L^2_{l_j^{k+1}}}\nonumber\\
&\qquad+\sum_{\gamma+\beta\leq s_{k+1}}\|\pa_t^\gamma\nabla_x^{\beta}(\rho_{k+1}, u_{k+1},  \theta_{k+1})(0)\|_{L^2_{x}}\bigg).
\end{align}
Finally, combining \eqref{7.23}, \eqref{7.28}, \eqref{7.30} and \eqref{7.35}, and noting \eqref{7.34-2}, we have
\begin{align}\label{7.39}
&\sum_{j=0}^{\bar{s}_{k+1}} \sum_{j=2\gamma+\beta} \|\fw_{\hat{\k}_{k+1}}\pa_t^{\gamma} \nabla_{\bar{x}}^{\beta} \bar{f}_{k+1}(t)\|^2_{L^2_{l_j^{k+1}}L^\infty_v} \nonumber\\
&\leq C\bigg(\tau, E_{s_0}, D_k+\bar{D}_k+\hat{D}_k,  \sum_{j=0}^{\bar{s}_{k+1}} \sum_{j=2\gamma+\beta}   \|\pa_t^{\gamma} \nabla_{\bar{x}}^{\beta} (\bar{u}_{k+1,\sp},\bar{\theta}_{k+1})(0)\|^2_{L^2_{l_j^{k+1}}}\nonumber\\
&\qquad+\sum_{\gamma+\beta\leq s_{k+1}}\|\pa_t^\gamma\nabla_x^{\beta}(\rho_{k+1}, u_{k+1},  \theta_{k+1})(0)\|_{L^2_{x}}\bigg)<\infty.
\end{align}

\

{\it Step 2.3.  Construction of solution $\hat{f}_{k+1}$.} Let $0\leq\zeta\leq \zeta_{k-1}$, by using \eqref{2.40}, \eqref{2.37} and \eqref{2.35}, we have
\begin{align}\label{7.40}
&\sum_{\gamma+\beta\leq \hat{s}_{k+1}}|e^{\zeta\cdot\eta}\partial_t^\gamma \pa_{\sp} ^{\beta} (\hat{A}_{k+1},\hat{B}_{k+1}, \hat{C}_{k+1})(t,x_{\sp} ,\eta)|\nonumber\\
&\leq C(E_{2+\hat{s}_{k+1}}) \sum_{\gamma+\beta\leq \hat{s}_{k+1}} e^{\zeta\cdot\eta}\int_{\eta}^{\infty}e^{-\zeta_{k-1}\cdot z} dz \, \|e^{\zeta_{k-1}\cdot z}\partial_t^\gamma \pa_{\sp} ^{\beta} (\hat{a}_{k+1},\hat{b}_{k+1}, \hat{c}_{k+1})(t,x_{\sp} ,\cdot)\|_{L^\infty_{z}}\nonumber\\
&\leq C(E_{2+\hat{s}_{k+1}}) \sum_{\gamma+\beta\leq\hat{s}_{k+1}} e^{-(\zeta_{k-1}-\zeta)\cdot \eta}\, \|e^{\zeta_{k-1}\cdot z}\partial_t^\gamma \pa_{\sp} ^{\beta} (\hat{a}_{k+1},\hat{b}_{k+1}, \hat{c}_{k+1})(t,x_{\sp} ,\cdot)\|_{L^\infty_{z}}\nonumber\\
&\leq C(E_{3+\hat{s}_{k+1}})  \sum_{\gamma+\beta\leq 1+\hat{s}_{k+1}} \| e^{\zeta_{k-1}\cdot \eta} \, \fw_{\hat{\k}_{k-1}} \partial_t^\gamma \pa_{\sp} ^{\beta} \hat{f}_{k-1}(t,x_{\sp} ,\cdot,\cdot)\|_{L^\infty_{\eta,v}},
\end{align}
and
\begin{align}\label{7.40-1}
&\sum_{\gamma+\beta\leq \hat{s}_{k+1}}|e^{\zeta\cdot\eta}\partial_t^\gamma \pa_{\sp} ^{\beta} \pa_{\eta} (\hat{A}_{k+1},\hat{B}_{k+1}, \hat{C}_{k+1})(t,x_{\sp} ,\eta)|\nonumber\\
&\leq C(E_{2+\hat{s}_{k+1}}) \sum_{\gamma+\beta\leq \hat{s}_{k+1}}  \|e^{\zeta\cdot\eta}\partial_t^\gamma \pa_{\sp} ^{\beta} (\hat{a}_{k+1},\hat{b}_{k+1}, \hat{c}_{k+1})(t,x_{\sp} ,\cdot)\|_{L^\infty_{z}}\nonumber\\
&\leq C(E_{3+\hat{s}_{k+1}})  \sum_{\gamma+\beta\leq 1+\hat{s}_{k+1}} \| e^{\zeta_{k-1}\cdot \eta} \, \fw_{\hat{\k}_{k-1}} \partial_t^\gamma \pa_{\sp} ^{\beta} \hat{f}_{k-1}(t,x_{\sp} ,\cdot,\cdot)\|_{L^\infty_{\eta,v}}.
\end{align}
Then, combining \eqref{7.40}, \eqref{7.40-1} with \eqref{2.38}, we obtain  the existence of  $\hat{f}_{k+1,1}(t)$ with
\begin{align}\label{7.43}
&\sum_{i=0,1}\sum_{\gamma+\beta\leq \hat{s}_{k+1}} \|e^{\frac32\zeta_{k+1}\cdot \eta} \fw_{\hat{\k}_{k+1,1}}\partial_t^\gamma \pa_{\sp} ^{\beta} \pa_{\eta}^{i}\hat{f}_{k+1,1}(t)\|_{L^\infty_{\hat{x},v}\cap L^2_{x_{\sp} }L^\infty_{\eta,v}} \nonumber\\
&\leq C(E_{3+s_{k+1}})  \sum_{\gamma+\beta\leq 1+\hat{s}_{k+1}} \| e^{\zeta_{k-1}\cdot \eta} \, \fw_{\hat{\k}_{k-1}} \partial_t^\gamma \pa_{\sp} ^{\beta} \hat{f}_{k-1}(t)\|_{L^\infty_{\hat{x},v}\cap L^2_{x_{\sp} }L^\infty_{\eta,v}},
\end{align}
where we have used $\hat{\k}_{k+1,1}\ll \hat{\k}_{k-1}$,  $1+\hat{s}_{k+1}\leq \hat{s}_{k-1}$, and  $\zeta=\frac32 \zeta_{k+1}$ such that  $0<\frac32\zeta_{k+1}\leq \zeta_{k-1}$. Moreover, from \eqref{2.40} and \eqref{2.38},  it direct to know that $\hat{f}_{k+1,1}$ is a continuous function over  $(t,x_{\sp} ,\eta,v)\in [0,\tau]\times\R^2\times \R_+\times \R^3$.

Using \eqref{7.43},  \eqref{2.38} and  the  trace theorem, one can obtain
\begin{align}\label{7.45}
&\sum_{\gamma+\beta\leq \hat{s}_{k+1}} \|e^{\frac32\zeta_{k+1}\cdot \eta}\fw_{\hat{\k}_{k+1,1}}\partial_t^\gamma \pa_{\sp} ^{\beta}\hat{f}_{k+1,1}(t,\cdot,0,\cdot)\|_{L^\infty_{x_{\sp} ,v}\cap L^2_{x_{\sp} }L^\infty_{v}} \nonumber\\
&\leq C(E_{3+s_{k+1}})  \sum_{\gamma+\beta\leq 1+\hat{s}_{k+1}} \| e^{\zeta_{k-1}\cdot \eta} \, \fw_{\hat{\k}_{k-1}} \partial_t^\gamma \pa_{\sp} ^{\beta} \hat{f}_{k-1}(t)\|_{L^\infty_{\hat{x},v}\cap L^2_{x_{\sp} }L^\infty_{\eta,v}}.
\end{align}

We still need to construct $\hat{f}_{k+1,2}$.  Firstly, it follows from \eqref{2.46},  \eqref{7.20}, \eqref{7.39}, the trace theorem and \eqref{7.45} that
\begin{align*}%\label{7.46}
&\sum_{\gamma+\beta\leq \hat{s}_{k+1}} \|\fw_{\hat{\k}_{k+1,1}}\pa_t^{\gamma}\nabla_{\sp} ^{\beta} \hat{g}_{k+1}(t)\|_{L^\infty_{x_{\sp} ,v}\cap L^2_{x_{\sp} }L^\infty_v}\nonumber\\
&\leq \sum_{\gamma+\beta\leq \hat{s}_{k+1}} \|\fw_{\hat{\k}_{k+1,1}}\pa_t^{\gamma}\nabla_{\sp} ^{\beta} (f_{k+1}, \bar{f}_{k+1}, \hat{f}_{k+1,1})(t,\cdot,0,\cdot)\|_{L^\infty_{x_{\sp} ,v}\cap L^2_{x_{\sp} }L^\infty_v}\nonumber\\
&\leq C\bigg(\tau, E_{s_0}, D_k+\bar{D}_k+\hat{D}_k,  \sum_{j=0}^{\bar{s}_{k+1}} \sum_{j=2\gamma+\beta}   \|\pa_t^{\gamma} \nabla_{\bar{x}}^{\beta} (\bar{u}_{k+1,\sp},\bar{\theta}_{k+1})(0)\|^2_{L^2_{l_j^{k+1}}}\nonumber\\
&\qquad+\sum_{\gamma+\beta\leq s_{k+1}}\|\pa_t^\gamma \nabla_x^{\beta}(\rho_{k+1}, u_{k+1},  \theta_{k+1})(0)\|_{L^2_{x}}\bigg),
\end{align*}
provided $2+\hat{s}_{k+1}\leq \bar{s}_{k+1}$.

On the other hand, using  \eqref{2.34} and Sobolev inequality, a direct calculation shows that
\begin{align*}%\label{7.47}
&\sum_{\gamma+\beta\leq \hat{s}_{k+1}} \|e^{\frac34 \zeta_{k}\cdot\eta} \fw_{\hat{\k}_{k+1,1}} \pa_t^\gamma \pa_{\sp} ^{\beta} \hat{S}_{k+1,2}(t)\|_{L^\infty_{\bar{x},v}\cap L^2_{x_{\sp} }L^\infty_{\eta,v}}\nonumber\\
&\leq C\Big(E_{\fb+2+\hat{s}_{k+1}},\sum_{i=1}^k\sum_{\gamma+\beta\leq \hat{s}_{k+1}+\fb} \big[\|\tw_{\k_i} \pa_t^\gamma \pa_{\sp} ^{\beta} f_i(t)\|_{L^\infty_{x,v}}
+\|\fw_{\hat{\k}_i} \pa_t^\gamma \pa_{\sp} ^{\beta} \bar{f}_{i}(t)\|_{L^\infty_{\bar{x},v}}\big]\nonumber\\
&\qquad\qquad+\sum_{i=1}^k\sum_{\gamma+\beta\leq 1+\hat{s}_{k+1}}\|\eta^{\fb}e^{\frac34 \zeta_{k}\cdot\eta} \fw_{\hat{\k}_i} \pa_t^\gamma \pa_{\sp} ^{\beta} \hat{f}_{i}(t)\|_{L^\infty_{\bar{x},v}\cap L^2_{x_{\sp} }L^\infty_{\eta,v}}\Big)\nonumber\\
&\leq C\bigg(\tau, E_{s_0}, D_k+\bar{D}_k+\hat{D}_k\Big),
\end{align*}
provided $s_0\geq m+2+\hat{s}_{k+1}$, $\hat{s}_k\geq 2+\hat{s}_{k+1}$ and $\hat{s}_{k+1}+\fb+2\leq \frac12 \bar{s}_k$. For $\hat{S}_{k+1,1}$, it is direct to have
\begin{align}\label{7.48}
&\sum_{\gamma+\beta\leq \hat{s}_{k+1}} \|e^{\frac34 \zeta_{k}\cdot\eta} \fw_{\hat{\k}_{k+1,1}} \pa_t^\gamma \pa_{\sp} ^{\beta} \hat{S}_{k+1,1}(t)\|_{L^\infty_{\bar{x},v}\cap L^2_{x_{\sp} }L^\infty_{\eta,v}}\nonumber\\
&\leq C(E_{2+\hat{s}_{k+1}}) \sum_{\gamma+\beta\leq \hat{s}_{k+1}}  \|e^{\frac34 \zeta_{k}\cdot\eta}\partial_t^\gamma \pa_{\sp} ^{\beta} (\hat{a}_{k+1},\hat{b}_{k+1}, \hat{c}_{k+1})(t,x_{\sp} ,\cdot)\|_{L^\infty_{\bar{x}}\cap L^2_{x_{\sp} }L^\infty_{\eta}}\nonumber\\
&\leq C(E_{s_0}, \hat{D}_{k-1}).
\end{align}

Let $0<\zeta_{k+1}\leq \frac12\zeta_k$ and $1\ll\hat{\k}_{k+1}\ll \hat{\k}_{k+1,1}\ll \bar{\k}_{k+1}$. Then, by using Lemma \ref{lem2.6},  \eqref{2.53-2}, and \eqref{7.43}-\eqref{7.48}, one establish the existence of solution $\hat{f}_{k+1,2}(t)$ over $t\in[0,\tau]$ with
 \begin{align}
& \sum_{\beta+\gamma\leq \hat{s}_{k+1}} \Big\{\|\fw_{\hat{\k}_{k+1}}e^{\zeta_{k+1}\eta}\partial_t^{\gamma}\nabla^{\beta}_{\sp} \hat{f}_{k+1,2}(t)\|_{L^\infty_{x_{\sp} ,\eta,v}\cap L^2_{x_{\sp} }L^\infty_{\eta,v}}
\nonumber\\
&\qquad\qquad\qquad+\|\fw_{\hat{\k}_{k+1}}\partial_t^{\gamma}\nabla^{\beta}_{\sp} \hat{f}_{k+1,2}(t,\cdot,0,\cdot)\|_{L^\infty_{x_{\sp} ,v}\cap L^2_{x_{\sp} }L^\infty_{v}}\Big\}\nonumber\\
&\leq  C\bigg(\tau, E_{s_0}, D_k+\bar{D}_k+\hat{D}_k,  \sum_{j=0}^{\bar{s}_{k+1}} \sum_{j=2\gamma+\beta}   \|\pa_t^{\gamma} \nabla_{\bar{x}}^{\beta} (\bar{u}_{k+1,\sp},\bar{\theta}_{k+1})(0)\|^2_{L^2_{l_j^{k+1}}}\nonumber\\
&\qquad+\sum_{\gamma+\beta\leq s_{k+1}}\|\pa_t^\gamma\nabla_x^{\beta}(\rho_{k+1}, u_{k+1},  \theta_{k+1})(0)\|_{L^2_{x}}\bigg),\nonumber
\end{align}
which, together with \eqref{7.43} and \eqref{7.45}, yields the existence of solution $\hat{f}_{k+1}$ satisfying
 \begin{align*}%\label{7.49}
& \sum_{\beta+\gamma\leq \hat{s}_{k+1}} \Big\{\|\fw_{\hat{\k}_{k+1}}e^{\zeta_{k+1}\eta}\partial_t^{\gamma}\nabla^{\beta}_{\sp} \hat{f}_{k+1}(t)\|_{L^\infty_{x_{\sp} ,\eta,v}\cap L^2_{x_{\sp} }L^\infty_{\eta,v}}
\nonumber\\
&\qquad\qquad\qquad+\|\fw_{\hat{\k}_{k+1}}\partial_t^{\gamma}\nabla^{\beta}_{\sp} \hat{f}_{k+1}(t,\cdot,0,\cdot)\|_{L^\infty_{x_{\sp} ,v}\cap L^2_{x_{\sp} }L^\infty_{v}}\Big\}\nonumber\\
&\leq  C\bigg(\tau, E_{s_0}, D_k+\bar{D}_k+\hat{D}_k,  \sum_{j=0}^{\bar{s}_{k+1}} \sum_{j=2\gamma+\beta}   \|\pa_t^{\gamma} \nabla_{\bar{x}}^{\beta} (\bar{u}_{k+1,\sp},\bar{\theta}_{k+1})(0)\|^2_{L^2_{l_j^{k+1}}}\nonumber\\
&\qquad+\sum_{\gamma+\beta\leq s_{k+1}}\|\pa_t^\gamma\nabla_x^{\beta}(\rho_{k+1}, u_{k+1},  \theta_{k+1})(0)\|_{L^2_{x}}\bigg).
\end{align*}

{\it Step 3. } Combining all above estimates and the induction assumption \eqref{7.11-1}, we have proved the existence of solutions $f_i, \bar{f}_i, \hat{f}_i, \, i=1,\cdots, N$ with
\begin{align*}
D_N+\bar{D}_N+\hat{D}_N&\leq C\bigg(\tau, E_{s_0}, \sum_{i=0}^{N}\sum_{j=0}^{\bar{s}_{i}} \sum_{j=2\gamma+\beta}   \|\pa_t^{\gamma} \nabla_{\bar{x}}^{\beta} (\bar{u}_{i,\sp},\bar{\theta}_{i})(0)\|^2_{L^2_{l_j^{i}}}\nonumber\\
&\qquad\qquad+\sum_{i=0}^{N}\sum_{\gamma+\beta\leq s_{i}}\|\pa_t^\gamma \nabla_x^{\beta}(\rho_{i}, u_{i},  \theta_{i})(0)\|_{L^2_{x}}\bigg),
\end{align*}
where we have chosen $s_i,\ \bar{s}_i,\ \hat{s}_i$ such that
\begin{align}\label{7.52}
\begin{split}
&s_0\geq s_1+\fb+6,\quad s_1=\bar{s}_1=\hat{s}_1\gg1;\\
& s_1>s_i>\bar{s}_i>\hat{s}_i\geq s_{i+1}>\bar{s}_{i+1}>\hat{s}_{i+1}\geq \cdots\gg1,\quad\mbox{for}\ i=2,\cdots,  N-1;\\
&s_{i+1}\leq \min\{\hat{s}_i,\, \frac12\bar{s}_{i}-3\},\, \bar{s}_{i+1}\leq s_{i+1}-8-\fb,\, \hat{s}_{i+1}\leq \frac12 \bar{s}_{i+1}-2-\fb,\mbox{for}\ i=1,\cdots,  N-1;
\end{split}
\end{align}
and taken $l_j^i=\bar{l}_j+2(\bar{s}_i-j)$ with $0\leq j\leq \bar{s}_i$ so that
\begin{equation}\label{7.53}
l^{N}_{j}\gg 2\fb \quad \mbox{and} \quad  l_j^i\geq 2 l_j^{i+1} +18+2\fb,\,  \mbox{for}\, 1\leq i\leq  N-1.
\end{equation}
Here we can taken $s_1=\bar{s}_1=\hat{s}_1$ because $f_1, \bar{f}_1$ and $\hat{f}_1$ depend only on the Euler solution, and do not depend on each other.  We also point out that $f_i, \bar{f}_i$ are smooth, but $\hat{f}_{i}$ is only continuous away from the grazing set $\{(x_{\sp} ,0, v)\, | \, x_{\sp} \in\R^2,  \, v_{\sp} \in\R^2, v_3\neq0 \}$. For the velocity weight functions, we demand
\begin{equation}\label{7.54}
\k_i\gg \bar{\k}_i\gg  \hat{\k}_i\gg \k_{i+1}\gg \bar{\k}_{i+1}\gg  \hat{\k}_{i+1}\gg1
\end{equation}
 for $1\leq i\leq N-1$, and we do not describe the precise relations between $\k_i, \bar{\k}_i$ and $\hat{\k}_i$ because the functions  $F_i, \bar{F}_i$ and $\hat{F}_i$ indeed decay exponentially with respect to particle velocity $v$. Therefore this completes the proof. $\hfill\Box$

\begin{remark}
To establish  the interior expansion $F_k$, viscous boundary layer $\bar{F}_k$ and Knudsen boundary layer $\hat{F}_k$, one should deal with the boundary interplay very carefully. In fact,  due to the  boundary effects,  one can only obtain the uniform estimates of  time and tangential derivatives for the Knudsen boundary layer $\hat{F}_k$. Fortunately, such  time and tangential derivatives estimates of Knudsen layer are enough to control the boundary interplay, see \eqref{7.15}-\eqref{7.17} and \eqref{7.33} for details.
\end{remark}

%%%%%%%%%%%%%%%%%%%%%%%%%%%%%%%%%%%%%%%%%%%%%%%%%%%%%%%%%%%%%%%%%%%%%%%%%%%%%%%%%%

\section{Hilbert Expansion: Proof of Theorem \ref{theorem} }\label{section6}

In this section, with the uniform estimates in Proposition \ref{prop5.1},  we shall use $L^2$-$L^\infty$ method to estimate the remainder term $F_R^\v$ in \eqref{1.22}  over half-space. Firstly,  from the formulation of boundary condition in section \ref{sec2.4}, it is easy to know that $F^\v_R$ satisfies the specular reflection boundary conditions, i.e.,
\begin{equation}\label{3.0}
F^\v_R(t,x,v)|_{\gamma_-}=F_R^\v(t,x_{\sp} ,0,v_{\sp} ,-v_3).
\end{equation}

\subsection{ $L^2$-energy Estimate} Recalling the definition of $f^\v_R$ in \eqref{1.20},  we  rewrite the equation terms of $f^\v_R$ as
\begin{align}\label{3.1}
\dis \partial_t f^\v_R+v\cdot\nabla_x f^\v_R+\frac{1}{\v^2}\mathbf{L}f^\v_R
&=-\frac{\{\partial+v\cdot\nabla_x\}\sqrt{\mu}}{\sqrt{\mu}}f^\v_R+\v^{3}\frac1{\sqrt{\mu}}Q(\sqrt{\mu}f^\v_R,\sqrt{\mu}f^\v_R)\nonumber\\
&\hspace{-6mm}+\sum_{i=1}^N\v^{i-2}\frac1{\sqrt{\mu}}\Big\{Q(F_i+\bar{F}_i+\hat{F}_i,\sqrt{\mu}f^\v_R)
+Q(\sqrt{\mu}f^\v_R, F_i+\bar{F}_i+\hat{F}_i)\Big\}\nonumber\\
&+ \frac1{\sqrt{\mu}}R^\v+\frac1{\sqrt{\mu}}\bar{R}^\v+\frac1{\sqrt{\mu}}\hat{R}^\v,
\end{align}
where $R^\v,\bar{R}^\v,\hat{R}^\v$ are defined in \eqref{1.9-3}, \eqref{1.9-4} and \eqref{1.9-5}, respectively. From \eqref{3.0}, we know that $f^\v_R$ satisfies specular reflection boundary conditions
\begin{equation}\label{3.1-1}
f^\v_R(t,x_1,x_2,0,v_1,v_2,v_3)|_{v_3>0}=f^\v_R(t,x_1,x_2,0,v_1,v_2,-v_3).
\end{equation}

\begin{lemma}\label{lem5.1}
Let $0<\frac{1}{2\alpha}(1-\alpha)<\fa<\frac12$,  $\k\geq 7$,  $N\geq 6$ and $\fb\geq 5$. Let $\tau>0$ be the life span of compressible Euler solution obtained in Lemma \ref{lem2.1-1}, then there exists a suitably small constant $\v_0>0$  such that for all $\v\in (0,\v_0)$,  it holds that
\begin{align}\label{3.3}
&\quad \frac{d}{dt}\|f^\v_R(t)\|_{L^2}^{2}+\frac{c_{0}}{2\v^2}\|\{\mathbf{I-P}\}f^\v_R(t)\|_{\nu}^{2}\leq C \Big\{1+ \v^{8}\|h_R^\v(t)\|_{L^\infty}^2 \Big\}
\cdot(\|f_R^\v(t)\|_{L^2}^{2}+1),
\end{align}
for  $t\in[0,\tau]$.
\end{lemma}

\noindent{\bf Proof.}  Multiplying \eqref{3.1} by $f^\v_R$ and integrating the resultant equation over $\mathbb{R}^3_+\times\mathbb{R}^3$, one obtains that
\begin{align}\label{3.5}
&\frac12 \frac{d}{dt} \|f^{\v}_R\|_{L^2}^2+\frac{c_0}{\v^2} \|\{\mathbf{I-P}\}f^{\v}_R\|_{\nu}^2
-\frac12\int_{\partial\mathbb{R}^3_+}\int_{\mathbb{R}^3} v_3 |f^\v_R(t,x_1,x_2,0,v)|^2dx_1dx_2 dv \nonumber\\
&= -\int_{\mathbb{R}^3_+}\intr \frac{\{\partial_t+v\cdot\nabla_x\}\sqrt{\mu}}{\sqrt{\mu}} |f^\v_R|^2 dv dx + \v^3\int_{\mathbb{R}^3_+}\intr \frac1{\sqrt{\mu}}Q(\sqrt{\mu}f^\v_R,\sqrt{\mu}f^\v_R) f^{\v}_R dvdx \nonumber\\
&\quad +\int_{\mathbb{R}^3_+}\intr\sum_{i=1}^N \v^{i-2}\frac1{\sqrt{\mu}}\Big\{Q(F_i+\bar{F}_i+\hat{F}_i,\sqrt{\mu}f^\v_R)
+Q(\sqrt{\mu}f^\v_R, F_i+\bar{F}_i+\hat{F}_i)\Big\} f^\v_Rdvdx \nonumber\\
&\quad+ \int_{\mathbb{R}^3_+}\intr \left\{\frac1{\sqrt{\mu}}R^\v+\frac1{\sqrt{\mu}}\bar{R}^\v+\frac1{\sqrt{\mu}}\hat{R}^\v\right\} f^\v_R dvdx.
\end{align}

Using the boundary condition \eqref{3.1-1}, it is direct to have
\begin{equation*}%\label{3.6}
\int_{\partial\mathbb{R}^3_+}\int_{\mathbb{R}^3} v_3 |f^\v_R(t,x_1,x_2,0,v)|^2dx_1dx_2 dv=0.
\end{equation*}

For any $\lambda>0$, by using similar arguments as in \cite{Guo Jang Jiang} and taking $\kappa\geq 7$, we obtain
\begin{align*}
\int_{\mathbb{R}^3_+}\intr \frac{\{\partial_t+(v\cdot\nabla_x)\}\sqrt{\mu}}{\sqrt{\mu}} |f^\v_R|^2 dv dx &\leq C \int_{\mathbb{R}^3_+}\intr |(\nabla_x\rho,\nabla_x\fu,\nabla_x T)| (1+|v|)^3 |f^\v_R|^2 dv dx \nonumber\\
&\leq C\left\{\int_{\mathbb{R}^3_+}\int_{|v|\geq \frac{\lambda}{\v}}+  \int_{\mathbb{R}^3_+}\int_{|v|\leq \frac{\lambda}{\v}}\right\} (\cdots) dvdx\nonumber\\
&\leq C\frac{\lambda}{\v^2} \|\{\mathbf{I-P}\}f^{\v}_R\|_{\nu}^2+C_\lambda (1+\v^4 \|h^\v_R\|_{L^\infty})  \|f^\v_R\|_{L^2},
\end{align*}
and
\begin{align*}
&\v^3 \int_{\mathbb{R}^3_+}\intr \frac1{\sqrt{\mu}}Q(\sqrt{\mu}f^\v_R,\sqrt{\mu}f^\v_R) f^{\v}_R dvdx=\v^3 \int_{\mathbb{R}^3_+}\intr \frac1{\sqrt{\mu}}Q(\sqrt{\mu}f^\v_R,\sqrt{\mu}f^\v_R) \{\mathbf{I-P}\}f^{\v}_R dvdx\nonumber\\
&\leq \v^3 \|\{\mathbf{I-P}\}f^{\v}_R \|_{\nu} \|h^\v_R\|_{L^\infty} \|f^\v_R\|_{L^2}
\leq \frac{\lambda}{\v^2} \|\{\mathbf{I-P}\}f^{\v}_R \|_{\nu}^2+C_{\lambda} \v^8 \|h^\v_R\|_{L^\infty}^2 \|f^\v_R\|_{L^2}^2.
\end{align*}

\

From \eqref{7.52}, it is noted that
\begin{equation*}
s_N>\bar{s}_N\geq 2\fb+4+\hat{s}_N,\quad \hat{s}_N\geq 1,
\end{equation*}
which, together with \eqref{7.53}, \eqref{7.54},  \eqref{7.2-0} and Sobolev imbedding theorem, yields  that, for $1\leq i\leq N$ and $t\in[0,\tau]$,
\begin{align}\label{3.10}
\begin{split}
&\sum_{k=0}^{2\fb+2}\Big\{\left\|\tw_{\k_i}(v)   \nabla^k_{t,x}f_i(t)\right\|_{L^2_{x,v}}+\left\|\tw_{\k_i}\nabla^k_{t,x} f_i(t)\right\|_{L^\infty_{x,v}}\Big\}\leq C_R(\tau),\\
&\sum_{k=0}^{\fb+2}\Bigg\{\left\|\fw_{\bar{\k}_i} (1+y)^{\fb+9} \nabla^k_{t,\bar{x}}\bar{f}_i(t)\right\|_{L^2_{\bar{x},v}}+\left\|\fw_{\bar{\k}_i} (1+y)^{\fb+9} \nabla^k_{t,\bar{x}}\bar{f}_i(t)\right\|_{L^\infty_{\bar{x},v}} \Bigg\}\leq C_R(\tau), \\
&\sum_{k=0,1}\Bigg\{\left\|\fw_{\bar{\k}_i} e^{\frac{1}{2^N}\cdot\eta} \nabla_{t,x_{\sp} }\hat{f}_i(t)\right\|_{L^2_{\hat{x},v}}+ \left\|\fw_{\bar{\k}_i} e^{\frac{1}{2^N}\cdot\eta} \nabla_{t,x_{\sp} }\hat{f}_i(t)\right\|_{L^\infty_{\hat{x},v}}\Bigg\}\leq C_R(\tau),
\end{split}
\end{align}
where we have denoted
\begin{align}
C_R(\tau):&=C\Bigg( \tau, \|(\varphi_0,\Phi_0, \vartheta_0)\|_{H^{s_0}}+\sum_{i=0}^{N}\sum_{\gamma+\beta\leq s_{i}}\|\pa_t^\gamma\nabla_x^{\beta}(\rho_{i}, u_{i},  \theta_{i})(0)\|_{L^2_{x}}\nonumber\\
&\qquad\qquad\qquad+\sum_{i=0}^{N}\sum_{j=0}^{\bar{s}_{i}} \sum_{j=2\gamma+\beta}   \|\pa_t^{\gamma} \nabla_{\bar{x}}^{\beta} (\bar{u}_{i,\sp},\bar{\theta}_{i})(0)\|^2_{L^2_{l_j^{i}}}
\Bigg).\nonumber
\end{align}

\

Noting \eqref{relation of mu and muM}, we have, for $1\leq i\leq N$,  that
\begin{align}\label{3.9}
\begin{split}
\left|w_{\k}(v) \frac{\sqrt{\mu_0}}{\sqrt{\mu}} \bar{f}_i(t,x_{\sp} ,y,v)\right|
&\leq C |w_{\k}(v) \mu_0^{-\fa} \bar{f}_i(t,x_{\sp} ,y,v)|\cdot \frac{\mu_0^{\frac12+\fa}}{\mu^{\frac12}}\\
&\leq C |w_{\k}(v) \mu_0^{-\fa} \bar{f}_i(t,x_{\sp} ,y,v)|\cdot (\mu_{M})^{(\frac12+\fa)\alpha-\frac12},\\[2mm]
\left|w_{\k}(v) \frac{\sqrt{\mu_0}}{\sqrt{\mu}} \hat{f}_i(t,x_{\sp} ,\eta,v)\right|
&\leq C |w_{\k}(v) \mu_0^{-\fa} \hat{f}_i(t,x_{\sp} ,y,v)|\cdot (\mu_{M})^{(\frac12+\fa)\alpha-\frac12}.
\end{split}
\end{align}
Taking $0<\frac{1}{2\alpha}(1-\alpha)<\fa<\frac12$, then it holds that
$ (\frac12+\fa)\alpha-\frac12>0$
which, together with  \eqref{ewf}, \eqref{3.10} and  \eqref{3.9}, yields that  the third term on RHS of \eqref{3.5} is bounded by
\begin{align}
&\|\{\mathbf{I-P}\}f^{\v}_R \|_{\nu} \|f^\v_R\|_{\nu} \cdot \sum_{i=1}^N \v^{i-2}  \left\{\|w_{\k} f_i\|_{L^\infty_{x,v}}+\|w_{\k} \frac{\sqrt{\mu_0}}{\sqrt{\mu}} \bar{f}_i\|_{L^\infty_{x,v}}+\|w_{\k} \frac{\sqrt{\mu_0}}{\sqrt{\mu}} \hat{f}_i\|_{L^\infty_{x,v}}\right\}\nonumber\\
&\leq C\frac1{\v}\|\{\mathbf{I-P}\}f^{\v}_R \|_{\nu} \|f^\v_R\|_{\nu}
\leq \frac{\lambda}{\v^2} \|\{\mathbf{I-P}\}f^{\v}_R \|_{\nu}^2+C_\lambda \|f^{\v}_R \|_{\nu}^2\nonumber\\
&\leq ( \lambda+C_{\lambda} \v^2) \frac{1}{\v^2} \|\{\mathbf{I-P}\}f^{\v}_R \|_{\nu}^2+C_\lambda \|f^{\v}_R \|_{L^2}^2.\nonumber
\end{align}

From \eqref{1.9-3} and $\eqref{3.10}_1$, a direct calculation shows that
\begin{align}\label{3.12}
\left(\int_{\mathbb{R}^3_+}\intr |\frac1{\sqrt{\mu}}R^\v|^2 dvdx\right)^{\frac12}\leq C \v^{N-6}.
\end{align}
It follows from \eqref{1.9-4}, \eqref{1.9-5} and \eqref{3.10} that
\begin{align}\label{3.14}
\begin{split}
\left(\int_{\mathbb{R}^3_+}\intr |\frac1{\sqrt{\mu}}\bar{R}^\v|^2 dvdx\right)^{\frac12}&\leq C (\v^{N-5.5}+\v^{\fb-4.5}),\\
\left(\int_{\mathbb{R}^3_+}\intr |\frac1{\sqrt{\mu}}\hat{R}^\v|^2 dvdx\right)^{\frac12}&\leq C (\v^{N-5}+\v^{\fb-3}).
\end{split}
\end{align}
Combining \eqref{3.12}-\eqref{3.14} and Cauchy inequality, it holds that
\begin{equation*}
\left| \int_{\mathbb{R}^3_+}\intr \left\{\frac1{\sqrt{\mu}}R^\v+\frac1{\sqrt{\mu}}\bar{R}^\v+\frac1{\sqrt{\mu}}\hat{R}^\v\right\} f^\v_R dvdx\right|\leq C(\v^{N-6} +\v^{\fb-5}) \|f^\v_R\|_{L^2}.
\end{equation*}
Hence \eqref{3.3} follows from above estimates. This completes the proof of Lemma \ref{lem5.1}. $\hfill\Box$

\subsection{Weighted $L^\infty$-estimate}

Given $(t,x,v),$ let $[X(s),V(s)]$ %=[X(s;t,x,v),V(s;t,x,v)]=[x+(s-t)v,v]$
be the backward bi-characteristics of the Boltzmann equation, which is determined by %\eqref{1.9}
\begin{equation*}%\label{ode}
\begin{cases}
\displaystyle\frac{dX(s)}{ds}=V(s), \quad \frac{dV(s)}{ds}=0,\\
\displaystyle[X(t),V(t)]=[x,v].
\end{cases}
\end{equation*}
The solution is then given by
\begin{equation*}%\label{def.bic}
[X(s),V(s)]=[X(s;t,x,v),V(s;t,x,v)]=[x-(t-s)v,v].
\end{equation*}

For each $(x,v)$ with $x\in \overline{\mathbb{R}}^3_+$ and $v_3\neq 0,$ we define its {\it backward exit time} $t_{\mathbf{b}}(x,v)\geq 0$ to be the last moment at which the
back-time straight line $[X(s;0,x,v),V(s;0,x,v)]$ remains in $\overline{\mathbb{R}}^3_+:$
\begin{equation*}%\label{exit}
t_{\mathbf{b}}(x,v)=\sup \{\tau \geq 0:x-\tau v\in\mathbb{R}^3_+\}.
\end{equation*}
We therefore have $x-t_{\mathbf{b}}v\in \partial \mathbb{R}^3_+ $ and $ x_3-t_{\mathbf{b}}v_3=0$.  We also define
\begin{equation*}%\label{xb}
x_{\mathbf{b}}(x,v)=x(t_{\mathbf{b}})=x-t_{\mathbf{b}}v\in \partial \mathbb{R}^3_+.
\end{equation*}
%We always have
Note that $v\cdot \vec{n}(x_{\mathbf{b}})=v\cdot \vec{n}(x_{\mathbf{b}}(x,v)) < 0$ always holds true.

For half space problem the back-time trajectory is very simple, and the particle hit  the boundary at most one time. More precisely, for the case $v_3<0$,  the back-time cycle is a straight line and does not hit the boundary; on the other hand, for $v_3>0$, the back-time cycle will hit the boundary for one time. Now  let $x\in \overline{\mathbb{R}}^3_+$, $(x,v)\notin \gamma _{0}\cup \g_{-}$ and
$
(t_{0},x_{0},v_{0})=(t,x,v)$,  the back-time cycle is defined as
\begin{equation*}%\label{4.4}
\left\{\begin{aligned}
X_{cl}(s;t,x,v)&=\mathbf{1}_{[t_1,t_0)}(s)\{x-v(t-s)\}+\mathbf{1}_{(-\infty,t_1)}(s)\{x-R_{x_{\mathbf{b}}}v(t-s)\},\\[1.5mm]
V_{cl}(s;t,x,v)&=\mathbf{1}_{[t_1,t_0)}(s)v+\mathbf{1}_{(-\infty,t_1)}(s)R_{x_{\mathbf{b}}}v,
\end{aligned}\right.
\end{equation*}
with
\begin{equation*}%\label{4.5}
(t_1,x_{\mathbf{b}})=(t-t_{\mathbf{b}}(x,v),x_{\mathbf{b}}(x,v)).
\end{equation*}
The explicit formula is
\begin{align}\label{4.6}
\begin{split}
&t_{\mathbf{b}}(x,v)=\frac{x_3}{v_3},\quad  \mbox{for} \ v_3>0 \quad \mbox{and} \quad t_{\mathbf{b}}(x,v)=\infty,\quad \mbox{for} \ v_3<0;\\
&V_{cl}(s)
=\begin{cases}
(v_1,v_2,v_3),\quad \mbox{if} \  s\in [t_1,t]\\
(v_1,v_2,-v_3), \quad \mbox{if} \  s\in (-\infty,t_1),
\end{cases}\\
&X_{cl}(s)=x-v(t-s),\quad \mbox{if}\  s\in [t_1,t],\\
&X_{cl}(s)
=\begin{cases}
x_1-v_1(t-s),\\
x_2-v_2(t-s),\\
-x_3+v_3(t-s),
\end{cases}
\mbox{if}\  s\in (-\infty,t_1).
\end{split}
\end{align}

\

As in \cite{Guo Jang Jiang-1,Guo Jang Jiang}, we denote
\begin{equation}
L_{M}g=-\frac{1}{\smum}\Big\{Q(\mu,\sqrt{\mu_M} g)+Q(\sqrt{\mu_M}  g,\mu)\Big\}=\nu(\mu)g-Kg,\nonumber
\end{equation}
where the frequency $\nu(\mu)$ has been defined in \eqref{1.19} and  $Kg=K_{2}g-K_{1}g$ with
\begin{align}
K_{1}g=&\intr\ints B(\theta)|u-v|^{\gamma}\sqrt{\mu_{M}(u)}\frac{\mu(v)}{\sqrt{\mu_{M}(v)}}g(u)dud\omega,
\nonumber\\
K_{2}g=&\intr\ints B(\theta)|u-v|^{\gamma}\mu(u')\frac{\sqrt{\mu_{M}(v')}}{\sqrt{\mu_{M}(v)}}g(v')dud\omega
\nonumber\\
&+\intr\ints B(\theta)|u-v|^{\gamma}\mu(v')\frac{\sqrt{\mu_{M}(u')}}{\sqrt{\mu_{M}(v)}}g(u')dud\omega.\nonumber
\end{align}

\begin{lemma}[\cite{Guo Jang Jiang}]\label{lem5.2-1}
It holds that $\displaystyle Kg(v)=\intr l(v,v')g(v') dv'$ where the kernel $l(v,v')$ satisfies
\begin{equation}\label{4.8}
|l(v,v')|\leq C\frac{\exp{\{-c|v-v'|^{2}\}}}{|v-v'|}.
\end{equation}
\end{lemma}

Letting $\displaystyle K_{w}g\equiv w_{\k} K(\frac{g}{w_{\k}})$, it follows from \eqref{1.9-2} and \eqref{def of h} that
\begin{align}\label{4.9}
\dis &\partial_t h^\v_R+v\cdot\nabla_x h^\v_R+\frac{\nu(\mu)}{\v^2} h^\v_R-\frac1{\v^2}K_w h^\v_{R} \nonumber\\
&=\sum_{i=1}^{N}\v^{i-2}\frac{w_{\k}(v)}{\sqrt{\mu_M(v)}}\Big\{Q(F_i+\bar{F}_i+\hat{F}_i,\frac{\sqrt{\mu_M} h^\v_R}{w_{\k}})
+Q(\frac{\sqrt{\mu_M} h^\v_R}{w_{\k}}, F_i+\bar{F}_i+\hat{F}_i)\Big\}\nonumber\\
&\quad+\v^{3}\frac{w_{\k}}{\sqrt{\mu_M}}Q\Big(\frac{\sqrt{\mu_M} h^\v_R}{w_{\k}},\frac{\sqrt{\mu_M} h^\v_R}{w_{\k}}\Big)+ \frac{w_\k}{\sqrt{\mu_M}}\big[R^\v+\bar{R}^\v+\hat{R}^\v\big].
\end{align}

\begin{lemma}\label{lem5.2}
For $t\in[0,\tau]$, it holds that
\begin{equation*}%\label{4.10}
\sup_{0\leq s\leq t}\|\v^{3}h^\v_R(s)\|_{L^\infty}
\leq C(t) \{\|\v^{3}h^\v_R(0)\|_{L^\infty}+\v^{N-1}+\v^{\fb}\}+\sup_{0\leq s\leq t}\|f^\v_R(s)\|_{L^2}.
\end{equation*}
\end{lemma}

\noindent{\bf Proof.}	For any $(t,x,v)$, integrating \eqref{4.9}  along the backward trajectory, one has that
\begin{align}\label{4.11}
&h^\v_R(t,x,v)
=\exp{\left\{-\frac{1}{\v^2}\int_{0}^{t}\nu(\xi)d\xi\right\}} h^\v_R(0,X_{cl}(0),V_{cl}(0))\nonumber\\
&\qquad +\frac{1}{\v^2}\int_{0}^{t}\exp{\left\{-\frac{1}{\v^2}\int_{s}^{t}\nu(\xi)d\xi\right\}} (K_{w}h^\v_R)(s,X_{cl}(s),V_{cl}(s))ds\nonumber\\
 &\qquad +\v^{3}\int_{0}^{t}\exp{\left\{-\frac{1}{\v^2}\int_{s}^{t}\nu(\xi)d\xi\right\}}\left(\frac{w_{\k}}{\smum}Q\left(\frac{h^\v_R\smum}{w_\k},\frac{h^\v_R\smum}{w_\k}\right)\right)(s,X_{cl}(s),V_{cl}(s))ds\nonumber\\
&\qquad +\int_{0}^{t}\exp{\left\{-\frac{1}{\v^2}\int_{s}^{t}\nu(\x)d\x\right\}} \cdot \left\{\sum_{i=1}^{N}\v^{i-2}\frac{w_\k}{\smum}Q\left(F_{i}+\bar{F}_i+\hat{F}_i,\frac{h^\v_R\sqrt{\mu_M}}{w_\k}\right)\right.\nonumber\\
&\qquad\qquad\qquad\qquad+\left.\sum_{i=1}^{N}\v^{i-2}\frac{w_\k}{\smum}Q\left(\frac{h^\v_R\sqrt{\mu_M}}{w_\k},F_{i}+\bar{F}_i+\hat{F}_i\right)\right\}(s,X_{cl}(s),V_{cl}(s))ds\nonumber\\
&\qquad +\int_{0}^{t}\exp{\left\{-\frac{1}{\v^2}\int_{s}^{t}\nu(\xi)d\xi\right\}} \left( \frac{w_\k}{\sqrt{\mu_M}}[R^\v+\bar{R}^\v+\hat{R}^\v]\right)(s,X_{cl}(s),V_{cl}(s))ds,
\end{align}
where we have denoted
\begin{equation}\nonumber
\nu(\xi)=\nu(\mu)(\xi,V_{cl}(\xi),X_{cl}(\xi)),
\end{equation}
for simplicity of presentation.

It follows from \eqref{4.6}  that $|V_{cl}(s)|\equiv|v|$. Then a direct calculation shows that
\begin{equation*}%\label{4.13}
\nu(\mu)\sim\nu_{M}(v):=\intr\ints B(v-u,\t)\mu_M(u)d\o du \cong (1+|v|),
\end{equation*}
and
\begin{equation*}%\label{4.14}
\int_{0}^{t}\exp{\left\{-\frac{1}{\v^2}\int_{s}^{t}\nu(\xi)d\xi\right\}}\nu(\mu)ds \lesssim \int_{0}^{t}\exp{\left\{-\frac{\nu_{M}(v)(t-s)}{C\v^2}\right\}}\nu_{M}(v)ds\leq O(\v^2).
\end{equation*}

For the first term on RHS of  \eqref{4.11}, it is easy to know that
\begin{equation*}%\label{4.15}
\left|\exp{\left\{-\frac{1}{\v^2}\int_{0}^{t}\nu(\xi)d\xi\right\}}h^\v_R(0,X_{cl}(0),V_{cl}(0))\right|\leq C\exp{\Big(- \frac{\nu_M(v)t}{C\v^2}\Big)}\|h^\v_R(0)\|_{L^\infty}.
\end{equation*}

\

It is direct to know that
\begin{align}
\left|\frac{w_\k(v)}{\smum}Q\left(\frac{h^\v_R\sqrt{\mu_M}}{w_\k},\frac{h^\v_R\sqrt{\mu_M}}{w_\k}\right)(s)\right| &\leq C\nu_{M}(v) \|h^\v_R(s)\|_{L^\infty}^{2}\leq  C\nu(s)\|h^\v_R(s)\|_{L^\infty}^{2},\nonumber
\end{align}
then the third term on RHS of \eqref{4.11} is bounded by
\begin{equation*}%\label{4.17}
C\v^{3}\int_{0}^{t}\exp{\left\{-\frac{1}{\v^2}\int_{s}^{t}\nu(\xi)d\xi\right\}}\nu(s)\|h^\v_R(s)\|_{L^\infty}^{2}ds\leq C\v^{5}\sup_{0\leq s\leq t}\|h^\v_R(s)\|_{L^\infty}^{2}.
\end{equation*}

Noting \eqref{relation of mu and muM}, \eqref{ewf}  and  \eqref{3.10}, a direct calculation shows that
\begin{align}
&\left|\sum_{i=1}^{N}\v^{i-2}\frac{w_\k}{\smum}\left\{Q\left(F_{i}+\bar{F}_i+\hat{F}_i,\frac{h^\v_R\sqrt{\mu_M}}{w_\k}\right)+Q\left(\frac{h^\v_R\sqrt{\mu_M}}{w_\k},F_{i}+\bar{F}_i+\hat{F}_i\right)\right\}(s)\right|\nonumber\\
&\leq C\nu_{M}(v)\|h^\v_R(s)\|_{L^\infty}\cdot\sum_{i=1}^{N}\v^{i-2}\left\|w_\k \, \Big[\frac{\sqrt{\mu}}{\sqrt{\mu_M}}f_i(s)+\frac{\sqrt{\mu_0}}{\sqrt{\mu_M}}\bar{f}_i(s)
+\frac{\sqrt{\mu_0}}{\sqrt{\mu_M}} \hat{f}_i(s) \Big]\right\|_{L^\infty}\nonumber\\
&\leq C\nu_{M}(v)\|h^\v_R(s)\|_{L^\infty}\cdot\sum_{i=1}^{N}\v^{i-2} \Big[\|\tw_{k_i}  f_i(s)\|_{L^\infty}+\|\fw_{\bar{k}_i}\bar{f}_i(s)\|_{L^\infty}
+\|\fw_{\hat{\k}_i} \hat{f}_i(s)\|_{L^\infty} \Big]\nonumber\\
&\leq C_R(s)\frac{\nu_{M}(v)}{\v}\|h^\v_R(s)\|_{L^\infty},\nonumber
\end{align}
where we have used  $\frac{1-\alpha}{2\alpha}<\fa<\frac{1}{2}$. Then  the fourth term on RHS of  \eqref{4.11} is bounded by
\begin{align*}%\label{4.18}
C_R(t)\frac1{\v}\int_{0}^{t}\exp{\left\{-\frac{1}{\v^2}\int_{s}^{t}\nu(\xi)d\xi\right\}}\nu_{M}(v)\|h^\v_R(s)\|_{L^\infty}ds\leq C_R(t) \v \sup_{0\leq s\leq t}\|h^\v_R(s)\|_{L^\infty}.
\end{align*}

Similarly, it follows from \eqref{1.9-3}-\eqref{1.9-5}, \eqref{relation of mu and muM} and  \eqref{3.10} that
\begin{align*}
\left|\left( \frac{w_\k}{\sqrt{\mu_M}}[R^\v+\bar{R}^\v+\hat{R}^\v]\right)(s)\right|\leq C_R(s) (\v^{N-6}+\v^{\fb-5}),
\end{align*}
which yields that   the last term on RHS of  \eqref{4.11} is bounded by $C_R(t) [\v^{N-4}+\v^{\fb-3}]$.

\

Let $l_w(v,v')$ be the corresponding kernel associated with $K_w$. Recalling \eqref{4.8} we have
\begin{align}\label{4.19}
|l_w(v,v')|\leq C\frac{w_\k(v')\exp{\{-c|v-v'|^{2}\}}}{w_\k(v)|v-v'|} \leq C\frac{\exp{\{-\frac{3}{4}c|v-v'|^{2}\}}}{|v-v'|}.
\end{align}
Now we can bound  the second term on RHS of \eqref{4.11} by
\begin{equation}\label{4.20}
\frac{1}{\v^2}\int_{0}^{t}\exp{\left\{-\frac{1}{\v^2}\int_{s}^{t}\nu(\xi)d\xi\right\}}\intr|l_w(V_{cl}(s),v')
h^\v_R(s,X_{cl}(s),v')|dv'ds.
\end{equation}
We denote $V'_{cl}(s_1)=V_{cl}(s_1;s,X_{cl}(s),v')$ and  $X'_{cl}(s_1)=X_{cl}(s_1;s,X_{cl}(s),v')$.
Using \eqref{4.11} again to \eqref{4.20}, then \eqref{4.20} is bounded by
\begin{align}\label{4.21}
&\frac{1}{\v^{4}}\int_{0}^{t}\exp{\left\{-\frac{1}{\v^2}\int_{s}^{t}\nu(\xi)d\xi \right\}}\intr\intr|l_w(V_{cl}(s),v')l_w(V'_{cl}(s_1),v'')|\nonumber\\
& \qquad\quad\times \int_{0}^{s}\exp{\left\{-\frac{1}{\v^2}\int_{s_{1}}^{s}\nu(v')(\xi)d\xi\right\}}|h^\v_R(s_{1},X'_{cl}(s_1),v'')|dv'dv''ds_{1}ds\nonumber\\
&+C \|h^\v_R(0)\|_{L^\infty}+ C_R(t)\Big\{\v\sup_{0\leq s\leq t}\|h^\v_R(s)\|_{L^\infty}+\v^{N-4}+\v^{\fb-3}\Big\}+C(t) \v^{5}\sup_{0\leq s\leq t}\|h^\v_R(s)\|_{L^\infty}^{2},
\end{align}
where we have denoted $\nu(v')(s)=\nu(\mu)(s,X_{cl}'(s), V_{cl}'(s))$ for simplicity of presentation. And we also used the following fact
\begin{align}\label{4.22}
\intr|l_w(v,v')|dv'\leq C(1+|v|)^{-1},
\end{align}
which is followed from \eqref{4.19}.

We now concentrate on the first term in \eqref{4.21}.  As in \cite{Guo Jang Jiang}, we divide it into the following several cases.
	
\noindent{\it Case 1.}  For $|v|\geq \fN$, by using  \eqref{4.22}, one deduces the following bound:
\begin{align*}%\label{4.25}
& \frac{C}{\v^4}\sup_{0\leq s\leq t}\|h^\v_R(s)\|_{L^\infty}\int_{0}^{t}\exp{\left\{-\frac{\nu_{M}(v)(t-s)}{C\v^2}\right\}}\intr|l_w(V_{cl}(s),v')|\nonumber\\
&\quad \times \int_{0}^{s}\exp{\left\{-\frac{\nu_{M}(v')(s-s_{1})}{C\v^2}\right\}}\intr|l_w(V_{cl}'(s_1),v'')|dv''ds_{1}dv'ds\nonumber\\
&\leq \frac{C}{\fN}\sup_{0\leq s\leq t}\|h^\v_R(s)\|_{L^\infty}.
\end{align*}
	
\noindent{\it Case 2.} For either $|v|\leq \fN, |v'|\geq 2\fN$ or $|v'|\leq 2\fN, |v''|\geq 3\fN$,  noting $|V_{cl}(s)|=|v|$ and $|V'_{cl}(s_1)|=|v'|$ we get either $|V_{cl}(s)-v'|\geq \fN$ or $|V'_{cl}(s_1)-v''|\geq \fN$, then either one of the following is valid for some small positive constant $0<c_1\leq \frac{c}{32}$ (where $c>0$ is the one in Lemma \ref{lem5.2-1}):
\begin{equation}\nonumber%\label{exppart}
\begin{split}
\dis |l_w(V_{cl}(s),v')|&\leq e^{-c_1\fN^{2}}|l_w(V_{cl}(s), v') \exp{\left(  c_1|V_{cl}(s)-v'|^{2}\right)}|,\\[1mm]
\dis |l_w(V'_{cl}(s),v'')|&\leq e^{-c_1\fN^{2}}|l_w(V'_{cl}(s_1),v'')\exp{\left(  c_1|V'_{cl}(s_1)-v'|^{2}\right)}|,
\end{split}
\end{equation}
which, together with \eqref{4.19}, yields that
\begin{equation}\label{4.26}
\begin{split}
\dis \intr|l_w(v,v')e^{c_1|v-v'|^{2}}|dv'\leq \frac{C}{1+|v|} \quad \mbox{and}\
\dis \intr|l_w(v',v'')e^{c_1|v'-v''|^{2}}|dv''\leq C\frac{C}{1+|v'|}.
\end{split}
\end{equation}
Hence, for  the case of $|v-v'|\geq \fN$ or $|v'-v''|\geq \fN$, it follows from \eqref{4.26} that
\begin{align*}%\label{4.27}
& \int_{0}^{t}\int_{0}^{s}\left\{\int\int_{|v|\leq \fN,|v'|\geq 2\fN}+\int\int_{|v'|\leq 2\fN,|v''|\geq 3\fN}\right\}(\cdots) dv''dv' ds_1ds\nonumber\\
& \quad \leq \frac{C}{\v^{4}}e^{-c_1\fN^2}\sup_{0\leq s\leq t}\|h^\v_R(s)\|_{L^\infty}\int_{0}^{t}\int_{0}^{s}\int|l_w(v,v')|\exp{\left\{-\frac{\nu_{M}(v)(t-s)}{C\v^2}\right\}}\nonumber\\
& \quad\quad \times\exp{\left\{-\frac{\nu_{M}(v')(s-s_{1})}{C\v^2}\right\}} dv'ds_{1}ds\nonumber\\
&\quad \leq C e^{-c_1\fN^2}\sup_{0\leq s\leq t}\|h^\v_R(s)\|_{L^\infty}.
\end{align*}
	
\noindent{\it Case 3a.} $|v|\leq \fN, |v'|\leq 2\fN, |v''|\leq 3\fN$. We note $\nu_{M}(v)\geq \nu_0>0$ where $\nu_0$ is a positive constant independent of $v$. Furthermore, we assume that $s-s_{1}\leq \lambda\v^2$ for some small $\lambda>0$ determined later. Then the corresponding part of the first term in \eqref{4.21} is bounded by
\begin{align*}%\label{4.28}
&\frac{C}{\v^{4}}\int_{0}^{t}\int_{s-\lambda\v^2}^{s}\exp{\left\{-\frac{\nu_0 (t-s)}{\v^2}\right\}}\exp{\left\{-\frac{\nu_0 (s-s_{1})}{\v^2}\right\}}\|h^\v_R(s_{1})\|_{L^\infty}ds_{1}ds\nonumber\\
&\leq C \sup_{0\leq s\leq t}\{\|h^\v_R(s)\|_{L^\infty}\}\cdot \frac{1}{\v^4}\int_{0}^{t}\exp{\left\{-\frac{\nu_0 (t-s)}{\v^2}\right\}}ds\cdot \int_{s-\lambda\v^2}^{s} ds_{1}\nonumber\\
&\leq C \lambda \sup_{0\leq s\leq t}\{\|h^\v_R(s)\|_{L^\infty}\}.
	\end{align*}
	
\noindent{\it Case 3b.} $|v|\leq \fN, |v'|\leq 2\fN, |v''|\leq 3\fN$ and $s-s_{1}\geq \lambda\v^2$. This is the last remaining case. We can bound the corresponding part of the first term in \eqref{4.21} by
\begin{align}\label{4.29}
&\frac{C}{\v^{4}}\int_{0}^{t}\int_{D}\int_{0}^{s-\lambda\v^2}\exp{\left\{-\frac{\nu_{0}(t-s)}{C\v^2}\right\}}\exp{\left\{-\frac{\nu_{0}(s-s_{1})}{C\v^2}\right\}}\nonumber\\[1mm]
&\quad \times |l_w(V_{cl}(s),v')l_w(V'_{cl}(s_1),v'') \cdot  h^\v_R(s_{1},X'_{cl}(s_1),v'')|ds_{1}dv'dv''ds,
\end{align}
where $D=\{|v'|\leq 2\fN,|v''|\leq 3\fN\}$. It follows from \eqref{4.19} that
\begin{align}\nonumber
\intr |l_w(v,v')|^2 dv'\leq C,
\end{align}
which, together with Cauchy inequality, yields that  \eqref{4.29} is bounded by
\begin{align}\label{4.33} &\frac{C}{\v^{4}}\Bigg\{\int_{0}^{t}\int_{D}\int_{0}^{s-\lambda\v^2}\exp{\left\{-\frac{\nu_0(t-s_{1})}{C\v^2}\right\}} |h^\v_R(s_{1},X'_{cl}(s_1),v'')|^2 dv''dv'ds_1ds\Bigg\}^{\frac12}\nonumber\\
&\times \Bigg\{\int_{0}^{t}\int_{D}\int_{0}^{s-\lambda\v^2}\exp{\left\{-\frac{\nu_0(t-s_{1})}{C\v^2}\right\}} |l_w(V_{cl}(s),v')l_w(V'_{cl}(s_1),v'')|^2 dv''dv'ds_1ds\Bigg\}^{\frac12}\nonumber\\
&\leq \frac{C_\fN}{\v^{2}}\Bigg\{\int_{0}^{t} \int_{D}\int_{0}^{s-\lambda\v^2}\exp{\left\{-\frac{\nu_0(t-s_{1})}{C\v^2}\right\}} |f^\v_R(s_{1},X'_{cl}(s_1),v'')|^2 dv''dv'ds_1ds\Bigg\}^{\frac12}\nonumber\\
&=\frac{C_\fN}{\v^{2}}\left\{\int_{0}^{t}\int_{D}\int_{\max\{t_1',0\}}^{s-\lambda\v^2}\exp{\left\{-\frac{\nu_0(t-s_1)}{\v^2}\right\}}|f^\v_R(s_{1},X'_{cl}(s_1),v'')|^2ds_{1}dv'dv''ds\right\}^{\frac12}\nonumber\\
&+\frac{C_\fN}{\v^{2}} \left\{\int_{0}^{t}\int_{D}\int_{0}^{t_1'} \mathbf{I}_{\{t_1'>0\}}\exp{\left\{-\frac{\nu_0(t-s_1)}{\v^2}\right\}}|f^\v_R(s_{1},X'_{cl}(s_1),v'')|^2ds_{1}dv'dv''ds\right\}^{\frac12},
\end{align}
where $t_1':=s-t_{\mathbf{b}}(X_{cl}'(s),v')$. To estimate \eqref{4.33}, we integrate over $v'$,  and make a change of variable $v'\longmapsto z:=X'_{cl}(s_1)$. From the explicit formula $\eqref{4.6}_{3,4}$, we have that
\begin{align}\label{4.34}
\frac{\partial z}{\partial v'}=\frac{\partial X'_{cl}(s_1)}{\partial v'}
= \left(
\begin{array}{ccc}
-(s-s_1) & 0 & 0  \\
0 & -(s-s_1)&0 \\
0 & 0 & -(s-s_1)
\end{array}
\right),\  \mbox{if} \ \max\{0,t_1'\}\leq s_1\leq s-\lambda \v^2,
\end{align}
and
\begin{align}\label{4.35}
\frac{\partial z}{\partial v'}=\frac{\partial X'_{cl}(s_1)}{\partial v'}
= \left(
\begin{array}{ccc}
-(s-s_1) & 0 & 0  \\
0 & -(s-s_1)&0 \\
0 & 0 & s-s_1
\end{array}
\right),\  \mbox{if} \ 0\leq s_1\leq t_1'.
\end{align}
From \eqref{4.34} and \eqref{4.35}, for both cases, it holds  that
\begin{equation*}%\label{4.36}
\left|\mbox{det}\left(\frac{\partial z}{\partial v'}\right)(s_1)\right|=(s-s_1)^3\geq (\lambda \v^2)^3>0, \ \mbox{for}\ s_1\in[0,s-\lambda \v^2],
\end{equation*}
which yields that
\begin{align}\label{4.37}
&\int_{|v'|\leq 2\fN}|f^\v_R(s_{1},X'_{cl}(s_1),v'')|^2dv'\leq \frac{C}{\lambda^{3}\v^{6}} \int_{\mathbb{R}^3_+}|\fep(s_{1},z,v'')|^{2}dz, \  \mbox{for} \  s_1\in[0,s-\lambda \v^2].
\end{align}
 Using \eqref{4.37}, we can further bound the two terms in RHS of  \eqref{4.33} by
\begin{equation*}%\label{4.38}
\frac{C_{\fN,\lambda}}{\v^{3}}\sup_{0\leq s\leq t}\|f^\v_R(s)\|_{L^2}.
\end{equation*}

Collecting all the above terms and multiplying them with $\v^{3}$,  for any small $\lambda>0$ and large $\fN>0$,  one obtains that
\begin{align}
\sup_{0\leq s\leq t}\{\|\v^{3} h^\v_R(s)\|_{L^\infty}\}
&\leq C_R(t)\Big\{\|\v^{3} h^\v_R(0)\|_{L^\infty}+\v^{N-1}+\v^{\fb}\Big\}+C(t) \v^{2}\sup_{0\leq s\leq t}\|\v^3 h^\v_R(s)\|_{L^\infty}^{2} \nonumber\\
&\hspace{-2mm}+C\left\{\lambda+\frac{1}{\fN}+C_R(t)\v\right\}\sup_{0\leq s\leq t}\|\v^3 h^{\v}_R(s)\|_{L^\infty} +C_{\fN,\lambda}\sup_{0\leq s\leq t}\|f^\v_R(s)\|_{L^2}.\nonumber
\end{align}
Noting $t\in[0,\tau]$, first taking $\fN\gg1$ large enough and  $\lambda>0$ small, and finally choosing $0<\v\leq \v_0$ with $\v_0$ small enough,
%\begin{align}\nonumber
%C\left\{\kappa +\frac{}N+C(t)\v\right\}&\leq \frac14,
%\end{align}
we deduce
\begin{align}\nonumber
\sup_{0\leq s\leq t}\{\|\v^{3} h^\v_R(s)\|_{L^\infty}\}
\leq C_R(t)\Big\{\|\v^{3} h^\v_R(0)\|_{L^\infty}+ \v^{N-1}+\v^{\fb}\Big\}+C \sup_{0\leq s\leq t}\|f^\v_R(s)\|_{L^2}.\nonumber
\end{align}
Therefore the proof of 	Lemma \ref{lem5.2} is completed. $\hfill\Box$

\subsection{Proof of Theorem \ref{theorem}}

With Lemma \ref{lem5.1} and \ref{lem5.2} in hand, the rest proof is the same as \cite{Guo Jang Jiang}.  We omit the details for simplicity of presentation. Therefore we complete the proof of Theorem \ref{theorem}.

\

%%%%%%%%%%%%%%%%%%%%%%%%%%%%%%%%%%%%%%%%%%%%%%%%%%%%%%%%%%%%%%%%%%%%%%%%%%%%%%%%%%%%%

\section{Acoustic Limit: Proof of  Theorem \ref{thm5.1}} \label{section7}

To prove Theorem \ref{thm5.1},  we firstly  derive the estimate for two solutions to compressible Euler equations \eqref{1.7}-\eqref{1.12-2} and acoustic systems \eqref{1.9}-\eqref{1.9-1}.
We define $(\varphi_d^\d, \Phi_d^\d, \vartheta_d^\d)$ as
\begin{align} \nonumber
\varphi_d^\d:=\frac{1}{\d^2}(\rho^{\d}-1-\d \varphi),\quad \Phi_d^\d:=\frac{1}{\d^2} (\fu^\d-\d \Phi),\quad \vartheta_d^\d:= \frac{1}{\d^2}(T^{\d}-1-\d \vartheta).
\end{align}
As in \cite{Guo Jang Jiang}, a direct calculation shows
\begin{align}\label{8.1}
\begin{cases}
\dis \pa_t\varphi_d^\d+(\fu^\d\cdot \nabla) \varphi_d^\d+ \rho^\d \, \mbox{\rm div} \Phi_d^{\d}+\d[\nabla\varphi\cdot \Phi_d^\d+\mbox{\rm div}\Phi \, \varphi_d^\d]=-\mbox{\rm div}(\varphi \Phi),\\[2mm]
\dis \rho^{\d} \pa_t \Phi_d^\d +\rho^\d (\fu^\d\cdot \nabla) \Phi_d^\d
+\nabla\big(\rho^\d  \vartheta_d^\d+T^\d \varphi_d^\d\big) \\[2mm]
\qquad-  \vartheta_d^\d\nabla \rho^\d- \varphi_d^\d\nabla T^\d+\d [\pa_t \Phi \, \varphi_d^\d+ \rho^\d (\Phi^\d_d\cdot \nabla) \Phi+\vartheta_d^\d \nabla \varphi+\varphi_d^\d \nabla \vartheta]\\[2mm]
\qquad=-\varphi \pa_t\Phi-\rho^\d (\Phi\cdot\nabla)\Phi-\nabla(\varphi\vartheta),\\[2mm]
\pa_t \vartheta_d^\d+(\fu^\d\cdot\nabla) \vartheta_d^\d+\frac23 T^\d\, \mbox{\rm div} \Phi_d^\d+\d[\nabla\vartheta\cdot \Phi_d^\d+\frac23\mbox{\rm div}\Phi \, \vartheta_d^\d]
=-\Phi\cdot \nabla \vartheta-\frac23 \vartheta\, \mbox{\rm div}\Phi,
\end{cases}
\end{align}
with $(t,x)\in [0,\tau]\times \R^3_+$,  and  the initial and  boundary conditions
\begin{equation}\label{8.2}
(\varphi_d^\d, \Phi_d^\d, \vartheta_d^\d)|_{t=0}=0\quad \mbox{and}\quad \Phi_{d,3}^\d(t,x_{\sp} ,0)\equiv0.
\end{equation}
It is direct to know that \eqref{8.1} is a linear hyperbolic system with characteristic boundary. Then we apply Lemma \ref{lem3.1} to \eqref{8.1}-\eqref{8.2} (the coefficients of \eqref{8.1} is slightly different, but Lemma \ref{lem3.1} is still applicable) to obtain
\begin{equation}\label{8.3}
\sup_{t\in[0,\tau]} \|(\varphi_d^\d, \Phi_d^\d, \vartheta_d^\d)(t)\|_{\MH^{k}(\R_+^3)}\leq C(\tau, \|(\varphi_0, \Phi_0, \vartheta_0)\|_{H^{s_0}(\R^3_+)})
\end{equation}
provided $k\leq s_0-2$ where $\MH^{s}(\R_+^3)$ is the notation defined \eqref{5.4}.

With \eqref{8.3} in hand, we can prove Theorem \ref{thm5.1} by using the same arguments as in section 3.2 of \cite{Guo Jang Jiang}. The details are omitted for simplicity of presentation. Therefore this  completes the proof of Theorem \ref{thm5.1}. $\hfill\Box$

%%%%%%%%%%%%%%%%%%%%%%%%%%%%%%%%%%%%%%%%%%%%%%%%%%%%%%%%%%%%%%%%%%%%%%%%%%%%%%%%%%%%%
%\appendix
\section{Appendix}

\begin{lemma}\label{lemA.1}
1). Let $\Omega_b:=\{(x_{\sp} , x_3): x_{\sp} \in\R^2,\ x_3\in[0,b] \}$ with $b\geq 1$. We assume $f,g\in H^1(\Omega_b)$,  it holds, for any $x_3\in [0,b]$,  that
\begin{align}\label{A.1}
\left|\int_{\R^2} (f g)(x_{\sp} , x_3) dx_{\sp} \right|\leq\|\pa_{x_3}(f,g)\|_{L^2(\Omega_b)}  \|(f,g)\|_{L^2(\Omega_b)} +\frac{1}{b}\|f\|_{L^2(\Omega_b)}
\|g\|_{L^2(\Omega_b)}.
\end{align}
For  $i=1,2$, we have
\begin{align}\label{A.2}
\left|\int_{\R^2} (\pa_{x_i}f \cdot g)(x_{\sp} , x_3) dx_{\sp} \right|
&\leq \|\pa_{x_3}(f,g)\|_{L^2(\Omega_b)}\|\pa_{x_i}(f,g)\|_{L^2(\Omega_b)}\nonumber\\
&\qquad+\frac{1}{b}\|\pa_{x_i}f\|_{L^2(\Omega_b)}\|g\|_{L^2(\Omega_b)}.
\end{align}

2). Let $f,g\in H^{1}(\R^3_+)$, and  $x_3\in \R_+$, it holds that
\begin{align}\label{A.3}
\left|\int_{\R^2} (\pa_{x_i} f \cdot g) (x_{\sp} ,x_3) dx_{\sp} \right|\leq \|\pa_{x_3}(f,g)\|_{L^2(\R^3_+)}\|\pa_{x_i}(f,g)\|_{L^2(\R^3_+)},\ \mbox{for}\ i=1,2.
\end{align}
\end{lemma}

\noindent{\bf Proof.}  We only prove \eqref{A.1} and \eqref{A.2} since  \eqref{A.3} can be obtained by taking the limit $b\rightarrow\infty$ in \eqref{A.2}.

Without loss of generality, we assume that $f,g \in C^2(\Omega_b) \cap H^1(\Omega_b)$. There exists a point $z_b\in [0,b]$ so that
\begin{equation}\label{A.4}
\int_{\R^2} (fg)(x_{\sp} ,z_b) dx_{\sp}  =\frac1{b} \int_0^b\int_{\R^2} (fg)(x_{\sp} ,x_3) dx_{\sp}  dx_3.
\end{equation}
For any given $x_3\in [0,b]$, by using \eqref{A.4}, we have
\begin{align}
\left|\int_{\R^2} (fg)(x_{\sp} ,x_3) dx_{\sp} \right| &\leq \left|\int_{z_b}^{x_3}\int_{\R^2} \pa_{z} (fg) dx_{\sp}  dz\right|+
\frac1{b} \left|\int_0^b\int_{\R^2} (fg)(x_{\sp} ,x_3) dx_{\sp}  dx_3\right|\nonumber\\
&\leq \int_{0}^{b}\int_{\R^2} |\pa_{z}f\cdot g|+|f\cdot \pa_z g| dx_{\sp}  dz+
\frac{1}{b}\|f\|_{L^2(\Omega_b)}
\|g\|_{L^2(\Omega_b)}\nonumber\\
 &\leq  \|\pa_{x_3}(f,g)\|_{L^2(\Omega_b)}\|(f,g)\|_{L^2(\Omega_b)}+\frac{1}{b}\|f\|_{L^2(\Omega_b)}
\|g\|_{L^2(\Omega_b)}.\nonumber
\end{align}
Hence we conclude \eqref{A.1}.

\vspace{2mm}

Similar as \eqref{A.4}, for $i=1,2$, there exists a point $z_b\in [0,b]$ so that
\begin{equation*}%\label{A.5}
\int_{\R^2} (\pa_{x_i}f\cdot g)(x_{\sp} ,z_b) dx_{\sp}  =\frac1{b} \int_0^b\int_{\R^2} (\pa_{x_i}f\cdot g)(x_{\sp} ,x_3) dx_{\sp}  dx_3.
\end{equation*}
For $x_3\in [0,b]$, integrating by parts w.r.t $x_i$, we obtain
\begin{align}
&\left|\int_{\R^2} (\pa_{x_i}f\cdot g)(x_{\sp} ,x_3) dx_{\sp}  \right| \nonumber\\
&\leq \left|\int_{z_b}^{x_3}\int_{\R^2} \pa_{z} (\pa_{x_i}f\cdot g) dx_{\sp}  dz\right|+
\frac1{b} \left|\int_0^b\int_{\R^2} (\pa_{x_i}f\cdot g)(x_{\sp} ,x_3) dx_{\sp}  dx_3\right|\nonumber\\
&\leq \left|\int_{0}^{b}\int_{\R^2} \pa_{z}f\cdot \pa_{x_i}g dx_{\sp}  dz\right| +\left|\int_{0}^{b}\int_{\R^2} \pa_{x_i}f\cdot \pa_z g dx_{\sp}  dz\right| +
\frac{1}{b}\|\pa_{x_i}f\|_{L^2(\Omega_b)}\|g\|_{L^2(\Omega_b)}\nonumber\\
&\leq  \|\pa_{x_3}(f,g)\|_{L^2(\Omega_b)}\|\pa_{x_i}(f,g)\|_{L^2(\Omega_b)}+\frac{1}{b}\|\pa_{x_i}f\|_{L^2(\Omega_b)}\|g\|_{L^2(\Omega_b)}.\nonumber
\end{align}
This completes \eqref{A.2}. $\hfill\Box$

\
%%%%%%%%%%%%%%%%%%%%%%%%%%%%%%%%%%%%%%%%%%%%%%%%%%%%%%%%%%%%%%%%%%%%%%%%%%%%%%%%%%%%%	
\vspace{1.5mm}	
	
\noindent{\bf Acknowledgments.} Y. Guo is supported by an NSF Grant (DMS \#1810868). Feimin Huang's research is partially supported by National Natural Sciences Foundation of China No. 11688101. Yong Wang's research  is partially supported by National Natural Sciences Foundation of China No. 11688101 and 11771429.

\end{document}